\documentclass[final]{siamltex}
\def\FigurePath{./fig}

\usepackage[dvips]{graphicx}
\usepackage{amsmath}
\usepackage{amsfonts}
\usepackage{amssymb}
\usepackage{subfigure}

\usepackage{color}
\usepackage{epsfig}

\def\mm#1{\begin{align}#1\end{align}}
\def\mmn#1{\begin{align*}#1\end{align*}}

\def\dt{\Delta t}

\def\no{{m-1}}% old time
\def\nn{{m}}% new time
\def\dtn{\Delta t_\nn}

\def\R{\mathbb{R}}
\def\Rd1{\R^d\times\R_+^0}
\def\mdot{\!\cdot\!}

\def\ff{\mathbf{f}}
\def\fn{{\ff_\nu}}%{f_\nu}
\def\FF{\mathcal{F}}
\def\FFn{{\FF_\nu^\nn}}

\def\gijn{{g_{ij}^\nn}}

\def\Ab{{\bar A}}
\def\At{{\tilde A}}
\def\Pb{{\bar P}}
\def\Pt{{\tilde P}}
\def\inn{{i,\nn}}

\def\O{{\Omega}}
\def\OT{{\Omega_T}}
\def\Ob{{\overline\Omega}}
\def\OTb{{\overline\OT}}
\def\GT{{\Gamma_T}}

\def\Oh{{\Omega_h}}
\def\OTh{{\Omega_{T,h}}}
\def\Oin{{V_i^\nn}}
\def\Ojn{{V_j^\nn}}

\def\Eint{{\mathcal E_{T,h}^{\textnormal{int}}}}
\def\Eext{{\mathcal E_{T,h}^{\textnormal{ext}}}}
\def\Gin{{\partial V_i^\nn}}%{{\Gamma_{i}^\nn}}
\def\Gij{{\Gamma_{ij}}}
\def\nij{{n_{ij}}}
\def\Gijn{{\Gamma_{ij}^\nn}}

\def\In{{I_\nn}}
\def\Vi{{V_i}}
\def\Vin{{V_i^\nn}}
\def\UU{\mathcal U}
\def\VV{\mathcal V}
\def\VVh{{{\mathcal V}_h}}
\def\myVVh{{{\mathcal V}_h}}
\def\VVh0{{\mathcal V_h^0}}

\def\eJ{{\varepsilon_J}}
\def\eO{{\varepsilon_\Omega}}
\def\egm{{\varepsilon_\Gamma^-}}
\def\egp{{\varepsilon_\Gamma^+}}
\def\egpm{{\varepsilon_\Gamma^\pm}}
\def\vp{\varphi}

\def\ebkn{\bar\eta_k^\nn}

\def\tol{T\!ol}
\def\cfl{C\!F\!L}
\def\myN{{\mathcal N}}%{\mathbb N}

\newtheorem{remark}{Remark}[section]

\def\LMMSvect#1{\mbox{\boldmath$ #1$}}

\newcommand{\replace}[2]{{\textcolor{blue}{ }}{\textcolor{black}{#2}}} 
\title{Adaptive timestep control for nonstationary solutions of the Euler equations\thanks{This work has been 
       performed with funding by the Deutsche
       Forschungsgemeinschaft in the Collaborative Research Center SFB 401
       "Flow Modulation and Fluid-Structure Interaction at Airplane Wings" of
       the RWTH Aachen, University of Technology, Aachen, Germany.}}

\author{Christina Steiner\thanks{Institut f\"ur Geometrie und Praktische Mathematik,
        RWTH Aachen University, D-Templergraben 55, D-52056 Aachen, Germany
        ({\tt steiner@igpm.rwth-aachen.de}).}
        \and Siegfried M\"uller\thanks{Institut f\"ur Geometrie und Praktische Mathematik,
        RWTH Aachen University, D-Templergraben 55, D-52056 Aachen, Germany
        ({\tt mueller@igpm.rwth-aachen.de}).}
        \and Sebastian Noelle\thanks{Institut f\"ur Geometrie und Praktische Mathematik,
        RWTH Aachen University, D-Templergraben 55, D-52056 Aachen, Germany
        ({\tt noelle@igpm.rwth-aachen.de}).}}

\begin{document}

%\bibliographystyle{abbrv}

%\title{Adaptive timestep control for nonstationary solutions of the Euler equations}
%\author{Christina Steiner, Siegfried M\"uller, Sebastian Noelle}
%\date{\today}

\maketitle

\begin{abstract} 
In this paper we continue our work on adaptive timestep control
for weakly nonstationary problems \cite{SteinerNoelle}. The core of the method is
a space-time splitting of adjoint error representations for target functionals
due to S\"uli \cite{Sueli} and Hartmann \cite{Hartmann1998}. The main new
ingredients are (i) the extension from scalar, 1D, conservation laws to the 2D
Euler equations of gas dynamics, (ii) the derivation of boundary conditions for
a new formulation of the  adjoint problem and (iii) the coupling of the adaptive
time-stepping with spatial adaptation. For the spatial adaptation, we use a
multi\replace{scale}{resolution}-based strategy developed by M\"uller \cite{Mueller:02}, and we
combine this with an implicit time discretization. The combined space-time
adaptive method provides an efficient choice of timesteps for implicit
computations of weakly nonstationary flows. The timestep will be very large in
\replace{regions}{time intervalls} of stationary flow, and becomes small when a perturbation enters the
flow field. The efficiency of the solver is investigated by means of an unsteady
inviscid 2D flow over a bump. 
\end{abstract}

%\textbf{Key words:} compressible Euler equations, weakly nonstationary flows,
%Finite Volume and Discontinuous Galerkin methods, adaptive time-stepping,
%adjoint error analysis, multi\replace{scale}{resolution} \replace{analysis}{decomposition}.

%\vspace*{1ex}

%\textbf{AMS subject classifications.}  (MSC2000)
%35L65, 76N15, 65M12, 65M15, 65M50.



\pagestyle{myheadings}
\thispagestyle{plain}
\markboth{CH. STEINER, S. M\"ULLER AND S. NOELLE}{ADAPTIVE TIMESTEP CONTROL}

\section{Introduction}

Today, there is broad consensus that the numerical solution of compressible flow
equations requires a highly resolved mesh to simulate accurately the different
scales of the flow field and its boundaries. Adaptive grid methods can
significantly improve the efficiency by concentrating cells only where they are
most required, thus reducing storage requirements as well as the computational
time. There has been a tremendous amount of research designing, analyzing and
implementing codes which are adaptive in space, see, e.g.,
\cite{clawpack,Mueller:02,KroenerOhlberger2000, Ohlberger2001} and references
therein. 

Here our interest is in timestep adaptation. For stationary problems, local
time\-steps which are linked to the spatial gridsize are commonplace, and they
are heavily built upon the fact that time-accuracy, or time synchronization is
not needed. On the other hand, for fully nonstationary flows, explicit algorithms
whose timestep is governed by the $\cfl$ restriction of at most unity are the
method of choice. In \cite{SteinerNoelle}, we began to explore one of the
remaining gaps, namely weakly nonstationary flows on which we will focus in the
following. Many real world applications, like transonic flight, are
perturbations of stationary flows. While time accuracy is still needed to study
phenomena like aero-elastic interactions, large timesteps may be possible when
the perturbations have passed.  For explicit calculations of nonstationary
solutions to hyperbolic conservation laws, the timestep is dictated by the
$\cfl$ condition due to Courant, Friedrichs and Lewy
\cite{CourantFriedrichsLewy1928}, which requires that the numerical speed of
propagation should be at least as large as the physical one. For implicit
schemes, the $\cfl$ condition does not provide a restriction, since the
numerical speed of propagation is infinite.  Depending on the equations and the
scheme, restrictions may come in via the stiffness of the resulting nonlinear
problem. These restrictions are usually not as strict as in the explicit case,
where the $\cfl$ number should be below unity. For implicit calculations, $\cfl$
numbers of much larger than 1 may well be possible. Therefore, it is a serious
question how large the timestep, i.e., the $\cfl$ number, should be chosen. 

A possible strategy has been investigated by Ferm and L\"otstedt
\cite{Loetstedt} based on  timestep control strategies for ODEs. Here a
Runge-Kutta-Fehlberg method is applied to the semi-discretized flow equations by
which the local spatial and temporal errors are estimated. These errors
determine the local stepsize in time and space. Later on, this idea was also
embedded in fully adaptive multiresolution finite volume schemes, see \cite{
Domingues2007,Domingues2008}. Alternatively  Kr\"oner and Ohlberger
\cite{KroenerOhlberger2000, Ohlberger2001} based their space-time adaptivity 
upon Kuznetsov-type a posteriori $L^1$ error-estimates for scalar conservation
laws.

In this paper we will use a space-time-split adjoint error representation to
control the timestep adaptation. For this purpose, let us briefly summarize the
space-time splitting of the adjoint error representation, see
\cite{Johnson4,BeckerRannacher1996,BeckerRannacher2000,Sueli,Hartmann1998}
for details. The error representation expresses the error in a target functional
as a scalar product of the finite element residual with the dual solution. This error
representation is decomposed into separate spatial and temporal components. The
spatial part will decrease under refinement of the spatial grid, and the
temporal part under refinement of the timestep. Technically, this decomposition
is achieved by inserting an additional projection. Usually, in the error
representation, one subtracts from the dual solution its projection onto
space-time polynomials. Now, we also insert the projection of the dual solution
onto polynomials in time having values which are $H^1$ functions with respect to
space.

This splitting can be used to develop a strategy \replace{for a local choice of
timestep.}{for an adaptive choice of timesteps $\dtn$. For simplicity, we work with
timesteps $\dtn$ which are uniform in space at each time $t^\nn$, so time adaptation only means that
in general, $\dtn\neq\dt^{\nn+1}$. However, note that the adaptive indicator developed here would permit us to
implement adaptive timesteps which are also local in space. See \cite{BergerOliger1984, DawsonKirby2000,MuellerStiriba2007} for more information on local timesteps.}\label{page_intro} 

In contrast to the results reported in \cite{Sueli,Hartmann1998} for 
scalar conservation laws we now investigate weakly nonstationary
solution to the 2D Euler equations. The timestep will be very large in \replace{regions}{time intervalls} of stationary flow, and becomes small when a perturbation enters the flow field. 

Besides applying well-established adjoint techniques to a new test problem, we
further develop a new technique (first proposed by the authors in
\cite{SteinerNoelle}) which simplifies and accelerates the computation of the
dual problem. Due to Galerkin orthogonality, the dual solution $\vp$ does not
enter the error representation as such. Instead, the relevant term is the
difference between the dual solution and its projection to the finite element
space, $\vp-\vp_h$. In \cite{SteinerNoelle} we showed that it is therefore
sufficient to compute the spatial gradient of the dual solution, $w=\nabla \vp$.
This gradient satisfies a conservation law instead of a transport equation, and
it can therefore be computed with the same conservative algorithm as the forward
problem \cite{SteinerNoelle}. The great advantage is that the conservative
backward algorithm can handle possible discontinuities in the coefficients
robustly.

A key step is to formulate boundary conditions for the gradient  $w = \nabla
\vp$ instead of $v$. Generally the boundary conditions for the dual problem
come from the weighting functions of the target functional, e.g., lift or drag.
To formulate boundary conditions for $w$ which are compatible with the target
functional, one has to lift the well-established techniques of characteristic
decompositions from the dual solution to its gradient. We will present details
on that in Section~\ref{sec.boundary_conditions}.

\replace{Starting with a very coarse, but adaptive spatial mesh and $\cfl$ below unity}{
We start with a very coarse, but adaptive spatial mesh. In time we first compute local timesteps
according to a given $\cfl$ number below unity. Then, for simplicity, we use a global timestep which is the minimum over these local timesteps. Using our error indicator,}
we establish timesteps which are well adapted to the physical problem at hand. The
scheme detects stationary time \replace{regions}{intervals}, where it switches to very high $\cfl$
numbers, but reduces the timesteps appropriately as soon as a perturbation
enters the flow field.

We combine our time-adaptation with the spatial adaptive multiresolution
technique \cite{Mueller:02}. In the early 90's Harten \cite{SM-Harten:96}
proposed to use  {\em multiresolution techniques} in the context of finite
volume schemes applied to hyperbolic conservation laws. He employed these
techniques to transform the arrays of cell averages associated with any given
finite volume discretization into a different format that reveals  insight into
the local behavior of the solution.  The cell averages on a given highest level
of resolution ({\em reference mesh}) are represented as cell averages on some
coarse level, where the fine scale information is encoded in arrays of {\em
detail coefficients} of ascending resolution.

In Harten's original approach, the multiresolution \replace{analysis}{decomposition} \cite{SM-Harten:96}
is used to control  a hybrid flux computation by which computational time for
the flux  computation can be saved, whereas the overall  computational
complexity is not reduced but still stays proportional to the number of  cells
on the uniformly fine reference mesh.   Opposite to this strategy,  threshold
techniques are applied to the multiresolution decomposition  in
\cite{Gottschlich-Mueller,Mueller:02,Cohen-Kaber-Mueller-Postel:02} where detail
coefficients below a threshold value are discarded. By means of the remaining
significant details a locally refined mesh is determined whose complexity is
significantly reduced in comparison to the underlying reference mesh. We will
use this approach for the local adaptation in space. These techniques have been
applied successfully to the \replace{}{stationary and instationary} Euler equations
\cite{BramkampLambyMueller2004}.

In the present work, we are interested to combine the  multi\replace{scale}{resolution}-based grid
adaptation with  adjoint techniques to solve efficiently nonstationary problems.
The advantage of this space adaptive method is that it also provides an
efficient break condition for the Newton iteration in the implicit time
integration, see Section~\ref{sec.Newton}.

The paper is organized as follows. We start with a brief description of the
fluid equations and their discretization by implicit finite volume schemes, see
Section~\ref{FAMS}.  The adjoint error control is presented in Section
\ref{aec}, followed by the space-time splitting and the error estimates in
Section~\ref{sec.space_time_splitting}. Section~\ref{sec.boundary_conditions} is
about the boundary conditions of the forward problem, the dual problem and the
conservative dual problem. In Section~\ref{section.numexp} we will present some
details on the numerical realization: The adaptive method in time and grid
generation. To improve the efficiency of the scheme we employ multiresolution
techniques. In Section~\ref{a4-ms} we give a short review of the multi\replace{scale}{resolution}
\replace{analysis}{decomposition}, upon which the adaptation in space is based. In Section
\ref{section.numexp_setup} we present the nonstationary test case, a 2D Euler
transonic flow around a circular arc bump in a channel. In
Sections~\ref{section.numexp_fully_implicit} and
\ref{section.numexp_explicit_implicit} results of the fully implicit and a
mixed explicit-implicit time adaptive strategy are presented to illustrate the
efficiency of the scheme.  In Section~\ref{section:conclusion} we summarize
our results.

\section{Governing equations and finite volume scheme}
\label{FAMS}

For the numerical simulation of nonstationary inviscid compressible fluids 
we solve the time-dependent Euler equations in $\Rd1$. 
These lead to a  system of conservation equations
\begin{align}
\label{eq.euler_class}
U_t + \nabla \cdot f(U) =0 \quad \textnormal{in } \Omega_T,\\
P_-(U^+)\,(\fn(U^+) -g)=0\quad \textnormal{on } \Gamma_T.
\label{eq.euler_bc}
\end{align}
Here $\O\subset\R^d$ is the spatial domain with boundary $\Gamma\!:=\!\partial
\O\subset\R^d$ and $\Omega_T = \Omega\times[0,T)\subset\Rd1$ is the space-time
domain with boundary $\Gamma_T\!:=\!\partial\OT\subset\Rd1$. ${U}=(\varrho,
\varrho\,{v},\varrho\,E)^{\mathrm{T}}$ is the vector of conservative variables
(density of mass, momentum, specific total energy)  and  ${f}=(\varrho\, {v}, 
\varrho\,vv^\mathrm{T}+pI, v\,(\varrho\,E+p))^{\mathrm{T}} = (f_1,\dots,f_{d\replace{}{+2}})^\mathrm{T}$
the array of the corresponding convective fluxes $f_i$, $i= 1,\dots,d\replace{}{+2}$,
in the $i$th coordinate direction.  $p$ is the
pressure and ${v}$ the fluid velocity.  The system of equations is closed by the
perfect gas equation of state $p=\varrho\,(E-0.5\,{v}^2)(\gamma-1)$ with
$\gamma=1.4$ (air).

We denote by $\ff(U):= (f(U),U)$ the space-time flux, by $\nu$ the space-time
outward normal to $\Omega_T$, and by  $\fn(U) := \ff(U)\cdot\nu$ the space-time
normal flux. $U^+$ is the interior trace of $U$ at the boundary $\GT$ (or any
other interface used later on). 
Given the boundary value $U^+$ and the corresponding Jacobian matrix $\fn'(U^+)$,
let $P_-(U^+)$ be the $(d+2)\times(d+2)$-matrix which realizes the
projection onto the \replace{eigenvectors}{eigenspace} of $\fn'(U^+)$ corresponding to negative
eigenvalues (see Section~\ref{sec.boundary_conditions}).
Then the matrix-vector product $P_-(U^+)\, \fn(U^+)$ is the incoming component
of the normal flux at the boundary, and it is prescribed in \eqref{eq.euler_bc}.
\replace{}{Accordingly we define $P_+(U^+)$, with respect to the positive
eigenvalues.}
See Section~\ref{sec.boundary_conditions} for details.

Since it is well-known that solutions $U$ will develop singularities in finite
time \replace{}{\cite{Oleinik63,Lax1964,Sideris1984}}, we pass to the weak formulation of \eqref{eq.euler_class},
\eqref{eq.euler_bc}. It is not fully understood to which space the weak solution
should belong, but loosely based upon the recent work
\cite{DelellisSzekelyhidi_2008} we assume that the solution space is
$$\UU := \textnormal{BV}(\OT).$$
This implies that $U, f(U) \in L^1_{\textnormal{loc}}(\OT)$ and that $d$-dimensional traces
exist. As the space of test functions we choose
$$\VV := \{\vp\in W^{1,\infty}(\OT)\,|\,\text{supp}\,\vp\subset\subset\OTb\},$$
which is consistent with the regularity theory in \cite{Tadmor1991}. Note that
the test functions may take non-zero boundary values. We call
$U\in\UU$ a weak solution of \eqref{eq.euler_class}, \eqref{eq.euler_bc} if for
all $\vp\in\VV$
\begin{align}
  & \int_\OT (U\vp_t + f(U)\nabla\vp)\,dV(x,t) 
     - \int_{\GT} P_+(U^+)\,\fn(U^+)\vp\,dS(x,t) \nonumber \\
  & \;\;\;\; = \; \int_{\GT} P_-(U^+)\,g\vp\,dS(x,t).
  \label{eq.euler_weak}
\end{align}  
We approximate \eqref{eq.euler_weak} by a first or second order finite volume
scheme with implicit Euler time discretization.   The computational spatial grid
$\Oh$ is a set of open cells $V_i$ such that
$$
\bigcup_{i}\;\overline{V_i} = \Ob.
$$
The intersection of the closures of two different cells is either empty or a
union of common faces and vertices. Furthermore let $\myN(i)$ be the set of
cells that have a common face with the cell $i$, $\partial V_i$ the boundary of
the cell $V_i$ and for $j\in\myN(i)$ let $\Gamma_{ij}:= \partial V_i  \cap
\partial V_j $ be the interface between the cells $i$ and $j$ and $n_{ij}$ the
outer spatial normal to $\Gamma_{ij}$ corresponding to cell $i$. Since we will
work on curvilinear grids, we require that the geometric consistency condition 
\mm{\sum_{j \in N ( i )} |\Gamma_{ij}| n_{ij} = 0
\label{eq.geometric_consistency}}
holds for all cells.

Let us define a partition of our time interval $I:=(0,T)$ into subintervals
$\In = [t_\no, t_\nn]$, $1\leq \nn \leq N$, where
\begin{align*}
0=t_0<t_1<\ldots<t_\nn<\ldots<t_N=T.
\end{align*}
The timestep size is denoted by  $\dtn:= t_\nn-t_\no$.
Later on this partition will be defined automatically by the
adaptive algorithm. 
We also denote the space-time cells and faces by $\Vin:=\Vi\times\In$ and 
$\Gijn:=\Gij\times\In$, respectively.
Given this space-time grid the implicit finite volume discretization of
\eqref{eq.euler_weak} can be written as
\mm{
  U_i^\nn+\frac{\dtn}{|V_i|}\,\sum_{j\in\myN(i)}
  |\Gamma_{ij}|\, F_{ij}^\nn
  = U_i^\no \;\;\;\text{for}\;\;\nn\geq 1.
  \label{eq.euler_fv}}
It computes the approximate cell averages $U^\nn_i$ 
of the conserved variables on the new time level. 
For interior faces $\Gij$, the canonical choice for the numerical flux is a
Riemann solver,
\replace{
\begin{align*}
F_{ij}^\nn := F_{\textnormal{riem}}(U_{ij}^\nn,U_{ji}^\nn,n_{ij})
\end{align*}
}
{%
\mm{\label{eq.euler_fv_fint} F_{ij}^\nn := 
F_{\textnormal{riem}}(U_{i}^\nn,U_{j}^\nn,n_{ij})}
}consistent with the normal flux $f_n(U) = f(U)\cdot n_{ij}$.
In the numerical experiments in Sections \ref{section.numexp_fully_implicit}
and \ref{section.numexp_explicit_implicit} we choose Roe's solver \cite{Roe1981}.
If $\Gij\subset\partial\OT=:\GT$, then we follow the definition \eqref{eq.euler_weak} of a weak solution
and define the numerical flux at the boundary by
\replace{
\begin{align*}
  F_{ij}^\nn := P_+(U_{ij}^\nn)\,\ff_{\nu_{ij}}(U_{ij}^\nn)
  + P_-(U_{ij}^\nn)\,g_{ij}^\nn,
  \end{align*}
  }
{
\begin{align}\label{eq.euler_fv_fext}
  F_{ij}^\nn := P_+(U_{i}^\nn)\,\ff_{\nu_{ij}}(U_{i}^\nn)
  + P_-(U_{i}^\nn)\,g_{ij}^\nn,
  \end{align}
  }
where $g_{ij}^\nn$ is the average of $g$ over $\Gijn$. 

For simplicity of presentation we neglect in our notation that due to higher
order reconstruction the numerical flux usually depends on an enlarged stencil
of cell averages.

\section{Adjoint error control - adaptation in time} \label{aec}

In order to adapt the timestep sizes we use a method which involves adjoint
error techniques. We have applied this approach successfully to Burgers'
equation in \cite{SteinerNoelle}. Now we present an extension of this approach
to systems of conservation laws.

Since a finite volume discretization in space and a backward Euler step in time
are a special case of a Discontinuous Galerkin discretization, techniques based
on a variational formulation can be transferred to finite volume methods.

The key tool for the time adaptive method is a space-time splitting of adjoint
error representations for target functionals due to S\"uli \cite{Sueli} and
Hartmann \cite{Hartmann1998}. It provides an efficient choice of timesteps for
implicit computations of weakly nonstationary flows. The timestep will be very
large in time \replace{regions}{intervalls} of stationary flow, and become small when a perturbation
enters the flow field.

\subsection{Variational Formulation}

In this section we rewrite the finite volume method as a Galerkin method, 
which makes it easier to apply the adjoint error control techniques.

Let us first introduce the space-time numerical fluxes. Let
$\Oin=V_i\times\In\in\OTh$ be a space-time cell, and let
$\gamma\subset\partial\Oin$ be one of its faces, with outward unit normal $\nu$.
There are two cases: if $\nu$ points into the spatial direction, then
$\gamma=\Gij\times\In$ and $\nu=(n,0)$. If it points into the positive time
direction, then $\gamma=V_i\times\{t_\nn\}$, and $\nu=(0,1)$. Now we define the
space-time flux by
\mm{
  \FFn(U_h) = \left\{ \begin{array}{cl}
  F_{ij}^\nn \;\; \text{from} \;\; \eqref{eq.euler_fv_fint} & 
  ~~\text{if}~~\nu=\nij \;\; \text{and}~~\gamma\in\Eint \\
  F_{ij}^\nn \;\; \text{from} \;\; \eqref{eq.euler_fv_fext} & 
  ~~\text{if}~~\nu=\nij \;\; \text{and}~~\gamma\in\Eext \\
  (1-\theta)U_i^\no + \theta U_i^\nn & 
  ~~\text{if}~~\nu=(0,1) \;\; \text{and}~~ \nn \geq 1 \\
  U_i^0 & ~~\text{if}~~\nu=(0,1) \;\; \text{and}~~ \nn=0
  \end{array}\right.
  \label{eq.space-time-flux}}
where $\Eint$ are the interior faces and $\Eext$ the boundary faces. In the
third case, $\theta\in[0,1]$, so the numerical flux in time direction is a
convex combination of the cell averages at the beginning and the end of the
timestep. Different values of $\theta$ will yield different time
discretizations, e.g., explicit Euler for $\theta =0$, implicit Euler for $\theta
=1$.

Let $\myVVh:=W^{1,\infty}(\OTh)$ be the space of piecewise Lipschitz-continuous
functions. Now we introduce the semi-linear form $N$ by
\mm{
  N &: \;\;\;\;\;\; \UU \times \myVVh \to \mathbf R \nonumber \\
  N(U,\vp) &:= \sum_\inn(\FFn(U),\vp)_{\Gin} \nonumber \\
  & \;\; - \sum_\inn\left( (U,\vp_{h,t})_{\Oin} + 
  (f(U),\nabla\vp)_{\Oin}. \right)
  \label{eq.var_form}}

Here and below the sum is over the set $\{(i,\nn)\,|\,\Oin\in\OTh\}$, i.e., all
gridcells. Now we rewrite the finite volume method \eqref{eq.euler_fv} as a
first order Discontinuous Galerkin method (DG0): 

\mm{
  \text{Find}\;\; U_h\in\VVh0 \;\;\;\;\text{such that}\;\;\;\;
  N(U_h,\vp_h) = 0  \quad\quad\forall \vp_h\in\VVh0,
  \label{eq.DG0}}
where $\VVh0$ is the space of piecewise constant functions over $\OTh$.
\begin{remark}\label{remark.dg0_galerkin_orthognal}
(i) For the DG0 method, $U_h,\vp\in\VVh0$ are piecewise constant, so the last two
terms in \eqref{eq.DG0}, containing derivatives of $\vp$, disappear. Moreover,
due to the geometric condition \eqref{eq.geometric_consistency}
\mmn{
\sum\limits_{\{j\,|\,\Gijn\subset\Gin\}}\ff(U_i^\nn)\cdot\nu_{ij}^\nn = 0}
holds for all cells $\Oin$. Therefore, the DG0 solution may be characterized by:
Find $U_h\in\VVh0$ such that
\mm{
  \sum_\inn(\FFn(U_h)-\fn(U^+_h),\vp_h)_{\Gin} = 0
  \quad\quad\forall \vp_h\in\VVh0. \label{eq.DG0_b}}
This form is convenient to localize our error representation later on.  

\replace{}{
(ii) The Galerkin orthogonality \eqref{eq.DG0_b} also holds for the {\bf first and higher
order finite volume schemes}. These schemes satisfy \eqref{eq.DG0_b} due to their conservation property for local cell averages, which correspond to piecewise constant test functions $\vp_h$.}
\replace{}{As we go along, we will see that all the error representations developed in this paper hold for higher order finite volume schemes, as well. The reason is that we  will only test with $\vp_h\in\VVh0$.}
\end{remark}

\subsection{Adjoint error representation for target functionals}
\label{sec.error_representations}

In this section we define the class of target functionals $J(U)$ treated in
this paper, state the corresponding adjoint problem and prove the
error representation which we will use later for adaptive timestep control.

Before we derive the main theorems, we would like to give a preview of an
important difference between error representations for linear and nonlinear
hyperbolic conservation laws.
For linear conservation laws (and many other linear PDE's), it is possible to
express the error in a user specified functional,
\mm{\label{eq:error_representation_a}
  \eJ &:= J(U)-J(U_h),}
as a computable quantity $\eta$, so
\mm{\label{eq:eJ}
  \eJ &= \eta}
(see, e.g., \cite{BeckerRannacher1996, BeckerRannacher2000, Sueli,
BarthLarson2002, HartmannHoustonSISC2002, HartmannHoustonJCP2002,
SueliHouston2003} and the references therein).
In general, $\eta$ will be an inner product of the numerical residual
with the solution of an adjoint problem. Below we will see that such a 
representation does not hold for nonlinear hyperbolic conservation laws.
The nonlinearity will give rise to an additional error $\egm$ on the
inflow boundary, an error $\egp$ on the outflow boundary, and a linearization
error $\eO$ in the interior domain. Altogether, the error representation
in Theorem~\ref{theorem:error_representation_h} will be
\mm{\label{eq:error_representation_b}
  \eJ+\egm+\egp+\eO &= \eta.}
Our adaptation is based on equidistributing this $\eta$.

Typical examples for the functional $J$ are the lift or the drag of a
body immersed into a fluid. To simplify matters we consider
functionals of the following form:
\mm{
  J(U) = (U,\psi)_{\OT} - (P_+(U^+) \fn(U^+),\psi_\Gamma)_{\GT},
  \label{eq.functional}}
where $\psi$ and $\psi_\Gamma$ are weighting functions in the interior
of the space-time domain $\Omega_T$ and at the boundary $\Gamma_T$.
\replace{We will give an example in Section~\ref{subsection:example_functional}.}
{This class of functionals  includes lift and drag and also the functional used
in the numerical eperiments in Section~\ref{section.numexp} 
(see Section~\ref{subsection:example_functional} for details).}

\replace{}{Let $U_h$ be a solution of the first or higher order finite volume scheme \eqref{eq.euler_fv}. As pointed out in Remark~\ref{remark.dg0_galerkin_orthognal}, these schemes are Galerkin orthogonal when tested with piecewise constant functions.}
As a first step towards deriving the identity \eqref{eq:error_representation_b},
we generalize the well-known error representation (see, e.g., Tadmor 
\cite[(2.16)]{Tadmor1991}) from
initial value problems ($\Omega=\R^d$) to the initial boundary value
problem~\eqref{eq.euler_class}--\eqref{eq.euler_bc}. For this, let
$\Ab=\Ab(U,U_h)$ be the averaged Jacobian
\replace{\begin{align}
  \Ab &:= \int\limits_0^1 \frac{d}{d\tau}
  f(U_h+\tau (U-U_h))d\tau. \label{eq.Ab}
\end{align}}
{
\begin{align}
  \Ab &:= \int\limits_0^1 \textnormal{D}
  f(U_h+\tau (U-U_h))d\tau. \label{eq.Ab}
\end{align}
}
Note that $\Ab$ is in general discontinuous \replace{}{with respect to $x$}, and it is conservative in the sense
that
\begin{align}\label{eq.AB conservative}
f(U)-f(U_h) &= \Ab\,(U-U_h).
\end{align}
Let $\Pb_\pm$ be the corresponding projection matrices. Then
\begin{theorem}\label{theorem:error_representation}
Suppose $\vp \in \VV$ solves the adjoint problem
\mm{
  \partial_t \vp + \Ab^\mathrm{T}\nabla\vp &= \psi \quad \; \; \mathrm{in} \; \OT,
  \label{eq.adjoint} \\
  \bar{P}^\mathrm{T}_+ (\vp -\psi_\Gamma) &=0 
  \quad \; \;\mathrm{on} \; \Gamma_T
  \label{eq.adjoint_bc}}
\replace{}{where $\psi$ defines the functional $J(U)$ in \eqref{eq.functional}.}  
Let $\eJ$ be defined by \eqref{eq:error_representation_a} and let
\mmn{
\egm &:= -((P_-(U^+)-P_-(U_h^+))g,\vp)_\GT \\
\egp &:= -(P_+(U^+)\fn(U^+)-P_+(U_h^+)\fn(U_h^+),\vp-\psi_\Gamma)_\GT.}
Then
\mm{\label{eq:error_representation}
    \eJ+\egm+\egp = \eta}
for all $\vp_h \in {\mathcal V}_h^0$
\mm{\eta &:= N(U_h,\vp) = N(U_h,\vp-\vp_h). \label{eq:eta}} 
\end{theorem}
Here $U^+$ and $U_h^+$ are the traces of $U$ and $U_h$ at the boundary  $\GT =
\partial \OT$. In particular, $U_h\equiv U_i^m$ for a boundary face $\Gijn$.
\begin{remark}
For the initial line $\Omega\times\{t=0\}\subset\GT$, the projections become
trivial,
\mmn{P_-(U)=P_-(U_h)=I,\quad P_+(U)=P_+(U_h)=0,}
%$$P_-(U)=P_-(U_h)=I,\quad P_+(U)=P_+(U_h)=0,$$
so 
\mmn{\egm=\egp=0.}
%$$\egm=\egp=0.$$
Similarly the boundary errors vanish at time $t=T$ and for supersonic spatial
boundaries. For subsonic spatial parts of the boundary, the boundary errors
cannot be computed a posteriori. Together with the error in the functional they
will be estimated by the approximate error 
representations~\eqref{eq:error_representation} 
and~\eqref{eq:error_representation_h}.
\end{remark}

For the adjoint problem \eqref{eq.adjoint} and \eqref{eq.adjoint_bc} the role
of time is reversed and hence $\bar{P}^\mathrm{T}_+$ plays the role of $P_-$ in
\eqref{eq.euler_bc}. Here $\psi_\Gamma$ comes from the weighting function in
the functional \eqref{eq.functional}. We will present details on the boundary
conditions for the dual problem in Section~\ref{sec.boundary_conditions}. Note
that the right-hand side in \eqref{eq:error_representation} depends on the
solution $U$ not only due to the boundary term $P_-(U^+)g$, but mainly because
$\vp$ is the solution of \eqref{eq.adjoint}.

We would like to give a short proof of
Theorem~\ref{theorem:error_representation}, since there are some subtleties due 
to the boundary conditions \eqref{eq.euler_bc} and \eqref{eq.adjoint_bc}.
\\[1ex]
{\em Proof of Theorem~\ref{theorem:error_representation}.}
By definitions \eqref{eq:error_representation_a} and \eqref{eq.functional} of $\eJ$ and $J$, resp.,
\mm{
  \eJ &= J(U)-J(U_h) \nonumber \\
  &= (U-U_h,\psi)_\OT - (P_+(U^+)\fn(U^+)-P_+(U_h^+)\fn(U_h^+),\psi_\Gamma)_\GT.
  \label{eq.proof1}}
Using \eqref{eq.adjoint}, \eqref{eq.AB conservative}, the definitions 
\eqref{eq.euler_weak} of a weak solution and \eqref{eq.var_form} of the
variational form, we obtain
\mm{&\phantom{==} (U-U_h,\psi)_\OT \nonumber \\
  &= (U-U_h,\vp_t+\Ab^\mathrm{T}\nabla\vp)_\OT \nonumber \\
  &= (U-U_h,\vp_t)_\OT+(\Ab(U-U_h),\nabla\vp)_\OT \nonumber \\
  &= \left((U,\vp_t)_\OT+(f(U),\nabla\vp)_\OT\right)
    - \left((U_h,\vp_t)_\OT+(f(U_h),\nabla\vp)_\OT\right) \nonumber \\
  &= (P_-(U^+)g+P_+(U^+)\fn(U^+),\vp)_\GT \nonumber \\
  &\phantom{==} + N(U_h,\vp)-\sum_\inn(\FFn(U_h),\vp)_{\Gin}
  \label{eq.proof2}}
Since $\vp$ is continuous, the fluxes across interior faces cancel each other.
Using the definition
\eqref{eq.euler_fv_fext} of the boundary fluxes, we obtain
\mm{\sum_\inn(\FFn(U_h),\vp)_{\Gin}
  &= \sum_\inn(\FFn(U_h),\vp)_{\Gin\cap\GT} \nonumber \\
  &= \sum_{\{(i,j,\nn)\,|\,\Gijn\subset\GT\}}
     \left(P_-(U_i^\nn)\gijn+P_+(U_i^\nn)\fn(U_i^\nn),
     \vp\right)_\Gijn \nonumber \\
  &= \left(P_-(U_h^+)g+P_+(U_h^+)\fn(U_h^+),\vp\right)_\GT
     \label{eq.proof3}}
where $j$ in the second line is chosen such that $\Gijn\subset\GT$.
Combining \eqref{eq.proof1}--\eqref{eq.proof3} yields
\mm{
  \eJ &= N(U_h,\vp)-\left(P_-(U_h^+)g
     +P_+(U_h^+)\fn(U_h^+),\vp\right)_\GT \nonumber \\
  &\phantom{=.}+\,\left(P_-(U^+)g+P_+(U^+)\fn(U^+),\vp\right)_\GT \nonumber \\
  &\phantom{=.}-\,\left(P_+(U^+)\fn(U^+)-P_+(U_h^+)\fn(U_h^+),\psi_\Gamma\right)_\GT
    \nonumber \\
  &= N(U_h,\vp)-\egm-\egp
    \nonumber \\
  &= N(U_h,\vp-\vp_h)-\egm-\egp.}
In the last step we have used the definition of the DG0 scheme \eqref{eq.DG0}.
This completes the proof. $\square$

\begin{remark}
(i) In \cite{Tadmor1991} Tadmor proves the well-posedness of the
Cauchy problem for the adjoint equation \eqref{eq.adjoint_bc} -- \eqref{eq.adjoint}
for scalar, convex, one-dimensional conservation laws. 
The key observation is that, if the forward solution $U$ has jump
discontinuities, then due to the entropy condition the jump of the
transport coefficient $\bar{A}^\mathrm{T}$ has a distinct sign. This makes it
possible to follow the characteristics of the adjoint problem
backwards in time.
\\[0.5ex]
(ii) Equivalently one can define the adjoint solution via a
variational formulation, see, e.g., \cite{BarthLarson2002,
HartmannHoustonSISC2002, HartmannHoustonJCP2002, SueliHouston2003}.
\end{remark}

Besides the unknown boundary error terms $\egpm$ a fundamental difficulty
remains if one tries to design an adaptive algorithm based on
Theorem~\ref{theorem:error_representation}: Since the exact solution $U$ is not
known, we cannot compute $\Ab$. Therefore we cannot approximate the solution
$\vp$ of the adjoint problem \eqref{eq.adjoint_bc}--\eqref{eq.adjoint} as well
as the error indicator $\eta$. 

The following theorem, in which we replace $\Ab$ and $\Pb_\pm$ by
\mm{
  \At := \replace{A(U_h)}{\textnormal{D}
  f(U_h)} \quad\text{and}\quad \Pt_\pm := P_\pm(\At),
  \label{eq:At}}
overcomes this difficulty.
\begin{theorem}\label{theorem:error_representation_h}
Suppose $\vp \in \VV$ solves the approximate adjoint problem
\mm{
  \partial_t \vp + \At^\mathrm{T}\nabla\vp &= \psi \quad \; \; \mathrm{in} \; \OT,
  \label{eq.adjoint_h} \\
  \Pt^\mathrm{T}_+ (\vp -\psi_\Gamma) &=0  \quad \; \;\mathrm{on} \; \Gamma_T,
  \label{eq.adjoint_h_bc}}
Let $\eJ$, $\egm$, $\egp$ and $\eta$ be as in 
Theorem~\ref{theorem:error_representation}, and let
\mmn{
  \eO &= (f(U)-f(U_h)-\At(U-U_h),\nabla\vp)_\OT.}
Then
\mm{\label{eq:error_representation_h}
  \eJ+\egm+\egp+\eO = \eta.}
\end{theorem}

\proof
The proof is almost the same as the one of
Theorem~\ref{theorem:error_representation}, except that we have to replace $\Ab$
in the third line of \eqref{eq.proof2}, $$(\Ab(U-U_h),\nabla\vp)_\OT,$$ by
$\At$. This yields the additional term
\mmn{
  ((\Ab-\At)(U-U_h),\nabla\vp)_\OT
  = (f(U)-f(U_h)-\At(U-U_h),\nabla\vp)_\OT
  = \eO.}
This completes the proof. $\square$
\\[1ex]
\replace{}{A formal Taylor series expansion suggests that $\eO$ is 
quadratic in the error $U-U_h$.}

\section{Space-time splitting and the error estimate}
\label{sec.space_time_splitting}

The error representation \eqref{eq:error_representation} is not yet suitable for
time adaptivity, since it combines space and time components of the residual and
of the difference $\vp-\vp_h$ of the dual solution and the test function. The
main result of this section is an error estimate whose components depend either
on the spatial grid size $h$ or the timestep $k$, but never on both. The key
ingredient is a space-time splitting of \eqref{eq:error_representation_h} based on
$L^2$-projections. Similar space-time projections were introduced previously in
\cite{Hartmann1998,Sueli} \replace{}{for the linear heat equation and linear advection}.
In \cite{SteinerNoelle} we adapted them to the \replace{}{nonlinear Burgers equation.
The projections can be generalized directly} to finite
element spaces and space-time Discontinuous Galerkin methods of arbitrary order.

Let $\Ojn= V_j \times \In$ and $P_{s,r}(\Ojn) = P_s(V_j)\times P_r(\In)$ be the
space of polynomials of degree $s$ on ${V_j}$ and $r$ on $\In$. Furthermore let
$\hat P_{\In}^r(\Oin)=\{w\in L^2(\Oin)|w(x,\cdot)\in P_r(\In),\forall x\in
V_j\}$, and $\hat P_{V_j}^s(\Oin)=\{w\in L^2(\Oin)|w(\cdot,t)\in
P_s(V_j),\forall t\in \In\}$. For $r \geq 0$ define the $L^2$-projection 
$\Pi_{\In}^r: L^2(\Oin) \to \hat P_{\In}^r(\Oin)$  onto piecewise polynomials in
time via
\begin{align}
  (U(x, \cdot)-\Pi_{\In}^r U(x, \cdot), \vp(x, \cdot))_{\In} = 0 
  \quad \forall \vp\in \hat P_{\In}^r(\Oin), \forall x\in V_j,
\end{align}
and for $s \geq 0$ define the $L^2$-projection
$\Pi_{V_j}^s: L^2(\Oin) \to \hat P_{V_j}^s(\Oin)$ 
onto piecewise polynomials in space  via
\begin{align}
  (U(\cdot, t)-\Pi_{V_j}^s U(\cdot, t), \vp(\cdot, t))_{V_j} = 0 
  \quad \forall \vp\in \hat P_{V_j}^s(\Oin), \forall t\in \In.
\end{align}
Similarly let the $L^2$-projection onto space-time polynomials
$\Pi_{\Oin}^{s,r}: = L^2(\Oin) \to  P_{s,r}(\Oin)$
be defined via
\begin{align}
  (U-\Pi_{\Oin}^{s,r}U,\vp)_{\Oin}=0\quad \forall \vp\in 
  P_{s,r}(\Oin).
\end{align}
Note that $\Pi_{\Oin}^{s,r} = \Pi_{{V_j}}^s \Pi_{\In}^r = \Pi_{\In}^r
\Pi_{{V_j}}^s$. First we choose $\vp_h$ in the error representation
\eqref{eq:error_representation_h} to be $\vp_h =\Pi_{h,k}^{s,r} \vp$, i.e.,
$\vp_h\mid_{\Oin} =\Pi_{{V_j}}^s \Pi_{\In}^r \vp= \Pi_{\In}^r \Pi_{{V_j}}^s
\vp$. This leads to the identity
\mm{
  \vp-\Pi_{h,k}^{s,r}\vp =\vp-\Pi_{\In}^r \vp + \Pi_{\In}^r
  \vp -\Pi_{h,k}^{s,r}\vp = (\textnormal{id}-\Pi_{\In}^r)
  \vp+(\textnormal{id}-\Pi_{V_j}^s )\Pi_{\In}^r \vp.
  \label{eq:Pi_hksr}}
Now we restrict ourselves to finite volume methods. These are based on
space-time cell averages, and therefore the corresponding order in the DG
context would be $r=s=0$, even for higher order FV schemes. Using
\eqref{eq:Pi_hksr}, we obtain the following splitting of the error representation
\eqref{eq:error_representation_h} in Theorem~\ref{theorem:error_representation_h}:
\mm{
  \eta &= N(U_h,\vp) \;= N(U_h,\vp-\Pi_{h,k}^{0,0}\vp) \nonumber\\
  &= N(U_h,(\textnormal{id}-\Pi_{\In}^0 )\vp
    +(\textnormal{id}-\Pi_{V_j}^0 )\Pi_{\In}^0 \vp)\nonumber\\ 
  &= N(U_h,(\textnormal{id}-\Pi_{\In}^0 )\vp)
    +N(U_h,(\textnormal{id}-\Pi_{V_j}^0 )\Pi_{\In}^0 \vp)\nonumber\\ 
  &= \sum_\inn(\FFn(U_h)-\fn(U_h^+),(\textnormal{id}-\Pi_{\In}^0 )\vp )_\Gin \nonumber \\
  &\phantom{=} + \sum_\inn(\FFn(U_h)-\fn(U_h^+),
    (\textnormal{id}-\Pi_{V_j}^0 )\Pi_{\In}^0 \vp )_\Gin \nonumber \\
  &=:\eta_k +\eta_h
  \label{eq:error_representation2}}
where $\eta_k$ is the time-component and $\eta_h$ the space-component of the
error representation $\eta$. To summarize, we have shown the following corollary
of Theorem~\ref{theorem:error_representation_h}:
\begin{corollary}\label{cor:error_representation_hk}
Under the assumptions of Theorem~\ref{theorem:error_representation_h}, the
following error representation holds:
\mm{\label{eq:error_representation_hk}
  \eJ+\egm+\egp+\eO = \eta_k+\eta_h.}
\end{corollary}
In \cite{SteinerNoelle} we showed by numerical experiments that $\eta_k$ depends
only an $k$ and that $\eta_h$ depends only on $h$, and that they both decrease
with first order. We will use the asymptotic behavior of the error term $\eta_k$
to derive an adaptation strategy in time.

\section{Boundary conditions, the conservative dual problem and functionals
at the boundary}
\label{sec.boundary_conditions}

In this section we present details of the boundary conditions for the forward
\eqref{eq.euler_class}--\eqref{eq.euler_bc} and the dual problem
\eqref{eq.adjoint_h}--\eqref{eq.adjoint_h_bc}. Then we will
recall the conservative approach to the dual problem, which we introduced in
\cite{SteinerNoelle} and derive boundary conditions for the  gradient of the
dual problem. \replace{}{Note, that we assume that the Jacobian matrix is diagonalizable. This assumption holds for the Euler equations.}

\subsection{Boundary conditions for the forward and the dual problem}
\label{sec.boundary_fp}

First we introduce some notation and state boundary conditions for the forward
problem. For simplicity of notation, we restrict ourselves to the spatial domain
$\Omega$ with boundary $\Gamma$, normal $n(x)$, and normal flux $f_n = f\cdot
n$. As in~\eqref{eq:At} let $\At = f_n'(U_h)$ be the Jacobian of $f_n$ evaluated
at the approximate solution $U_h$, and let $\Pt_\pm$ be the projection matrices
which map vectors onto the eigenspaces of $\At$ corresponding to positive and
negative eigenvalues, respectively. They are defined in detail in
\eqref{eq.P_pm} below.

In order to explain the boundary conditions in \eqref{eq.euler_bc},
\eqref{eq.euler_fv_fext}, \eqref{eq.adjoint_bc}, and \eqref{eq.adjoint_h_bc}, we recall the theory of
boundary value problems for hyperbolic systems, see, e.g., \cite{lax,
Kreiss-Lorenz}. Boundary values have to be prescribed along characteristics
entering the domain. Therefore the solution, or the fluxes, have to be
decomposed into in- and outgoing components.

Let $L=L(\At)$ and $R=R(\At)$ denote the $(d+2)\times(d+2)$ matrices of the left and
right eigenvectors of ${\At}$, and $\Lambda=\Lambda(\At) = \textnormal{diag}(\lambda_1(\At),
\dots,\lambda_{d+2}(\At))$ the diagonal matrix of the eigenvalues of $\At$, so
\mmn{
  \At = R \Lambda L.}
As usual, the positive and negative parts of $\At$ are
\mmn{
  {\At_\pm} = R\Lambda^\pm L.}
We now introduce the notations
\mm{
  \Pt_{\pm}  := R D_\pm L,
  \label{eq.P_pm}}
where $D_\pm$ is the diagonal matrix $D_\pm :=\textnormal{diag}(\chi^\pm(\lambda_i))$
with 
\mmn{
  \chi^\pm (\lambda)= \max(0, \textnormal{sign}(\pm \lambda)).}
Then we observe the identities
\mmn{
  \Pt_+\Pt_-= \Pt_-\Pt_+=0, \quad\quad \Pt_\pm^2=\Pt_\pm, \quad\quad \Pt_\pm 
  \At = R D_\pm L R \Lambda L = \At_\pm.}
Note that $\Pt_\pm$ and $\At$ commute:
\mmn{
  \At \Pt_\pm = R\Lambda LR D_\pm L 
  = R\Lambda D_\pm L= R D_\pm\Lambda L 
  =  R D_\pm LR\Lambda L = \Pt_\pm {\At}.}
  
This specifies the boundary fluxes \eqref{eq.euler_fv_fext} for the forward
finite volume solver and \eqref{eq.adjoint_h_bc} for the linearized adjoint
problem. The boundary conditions \eqref{eq.euler_bc} and \eqref{eq.adjoint_bc}
are derived analogously.

\subsection{The conservative dual problem}
\label{subsection:conservative_dual_problem}

The adjoint equation \eqref{eq.adjoint_h} is a system of linear transport
equations with discontinuous coefficients. Therefore, numerical approximations
may easily become unstable. Another inconvenience is that in order to obtain a
meaningful error representation in \eqref{eq:error_representation_h}, the
approximate adjoint solution $\vp$ should not be contained in $\VVh0$.
Therefore, $\vp$ is often computed in the more costly space $\VV_h^1$.

In \cite{SteinerNoelle} we have proposed a simple alternative which helps to
avoid both difficulties. Instead of computing the dual solution $\vp$ we will
compute its gradient
\mmn{w := \nabla\vp,}
which is the solution of the conservative dual problem
\mm{
  w_t+\nabla\cdot(\At^\mathrm{T}w) = \nabla \psi \;\;\;\text{in}\;\OT.
  \label{eq.adjoint_w}}
This system is in conservation form, and therefore it can be solved by any
finite volume or Discontinuous Galerkin scheme. Moreover, \eqref{eq.adjoint_w}
may be solved in $\VVh0$, since a piecewise constant solution $w$ already
contains crucial information on the gradient of $\vp$. 

The scalar problem treated in \cite{SteinerNoelle} was set up in such a way that
the characteristic boundary conditions for the dual problem became trivial. In
the following, we develop the boundary conditions in the more general case
needed in the present paper.

Denoting the flux in \eqref{eq.adjoint_w} by $H:=\At^\mathrm{T}w$, the boundary condition
\eqref{eq.adjoint_h_bc} becomes
\mm{
  \Pt_+^\mathrm{T}\,(H-H_\Gamma) = 0 \;\;\;\text{on}\;\GT,
  \label{eq.adjoint_w_bc}}
i.e., we prescribe the incoming component $\Pt_+^\mathrm{T}H$. Here $H_\Gamma$ is a given
real-valued vector function, which depends on $\psi_\Gamma$.  However, this
characteristic boundary condition needs to be interpreted carefully. Using
\eqref{eq.P_pm} and \eqref{eq.adjoint_w} and denoting the interior trace at the
flux by $H_{\textnormal{int}}$, we may introduce the boundary flux by \mmn{H := \Pt_-^\mathrm{T} \,
H_{\textnormal{int}} + \Pt_+^\mathrm{T} \, H_{\Gamma} \quad \textnormal{on } \GT.} Note that all the
projections $\Pt_\pm$ used below depend on the point $(x,t)\in\GT$ via
the outside normal vector $\nu(x,t)$. The value $\Pt_-^\mathrm{T} H_{\textnormal{int}}$ may be assigned
from the trace $w_{\textnormal{int}}$ at the interior of the computational domain, $$ \Pt_-^\mathrm{T}
\, H_{\textnormal{int}} =\Pt_-^\mathrm{T} \,(\At^\mathrm{T}w)_{\textnormal{int}}. $$ The boundary values $\Pt_+^\mathrm{T} \, H_\Gamma$
are computed using the PDE 
\mm{
  \vp_t = -H +\psi
  \label{eq.v_t}}
with boundary values \eqref{eq.adjoint_h_bc},
\mm{
  \Pt_+^\mathrm{T} \, H_\Gamma 
  &= \Pt_+^\mathrm{T} \,(-\vp_t +
    \psi)|_\Gamma \nonumber \\
  & = -(\Pt_+^\mathrm{T}\,\psi_\Gamma)_t 
    + \Pt_+^\mathrm{T}\,\psi \nonumber \\
  & \approx - \frac1{\dtn} \left(\Pt_+^{\mathrm{T},\nn}
    \,\psi_\Gamma^\nn - \Pt_+^{\mathrm{T},\no}
    \,\psi_\Gamma^\no\right) 
    + \Pt_+^{\mathrm{T},\no}\,\psi.}
This completes the definition of
the numerical boundary conditions for the conservative dual problem.

\subsection{The time component of the error representation}

Let us have another look at the time component of the error representation
\eqref{eq:error_representation2},
\mmn{
  \eta_k &= \sum_\inn(\FFn(U_h)-\fn(U_h^+),
  (\textnormal{id}-\Pi_{\In}^0)\vp)_{\partial\Oin}.}
To compute the leading order part of $\eta_k$ we assume that
$\vp\in\VV_h^1$. In this case, 
\mmn{
  (\textnormal{id}-\Pi_{\In}^0)\vp(x,t) = \left(t-\frac{t^\nn+t^\no}2\right) \vp_t.}
Note that $\vp_t$ is piecewise constant. Using \eqref{eq.v_t} we obtain
\mmn{
  \eta_k &= \sum_\inn\left(\FFn(U_h)-\fn(U_h^+),
  \left(\cdot-\frac{t^\nn+t^\no}2\right)\,(\psi-\At^\mathrm{T}w)\right)_{\partial\Oin}.}
Since $U_h$ is piecewise constant, the integrals over the time-like faces
$\partial V_i\times \In$ drop out, and only those over the
space-like faces $V_i\times\{t_\nn\}$ remain, so
\mmn{
  \eta_k &= \sum_\inn
  \left((\FFn(U_h)-\fn(U_h^+),\frac{\dtn}2\,(\psi-\At^\mathrm{T}w))_{V_i} \right.
  \nonumber \\
  &\phantom{..} - \left.
  (\FF_\nu^\no(U_h)-\fn(U_h^+),\frac{\dtn}2\,(\psi-\At^\mathrm{T}w))_{V_i}\right)
  \nonumber \\
  &= \sum_\inn
  (\FFn(U_h)-\FF_\nu^\no(U_h),\frac{\dtn}2\,(\psi-\At^\mathrm{T}w))_{V_i}.}
From the definition \eqref{eq.space-time-flux} of the flux in time direction
this simplifies further, namely
\mm{
  \eta_k &= \frac{\dtn}2\,\sum_\inn
  \left((1-\theta)(U_i^\no-U_i^{\nn-2})
  +\theta(U_i^\nn-U_i^\no),\psi-\At^\mathrm{T}w\right)_{V_i}
  \label{eq.eta_k}}
(we set $U_i^{-1}=U_i^{0}$ in the first summand). Thus our temporal error
indicator is simply a weighted sum of time-differences of the approximate
solution $U_h$, and the weights can be computed from the data $\psi$ and the
solution $w$ of the conservative dual problem~\eqref{eq.adjoint_w}. 

In our adaptive strategy, we will use the localized indicators
\mm{
  \ebkn := \frac12\,\sum\limits_i
  \left|((1-\theta)(U_i^\no-U_i^{\nn-2})
  +\theta(U_i^\nn-U_i^\no),\psi-\At^\mathrm{T}w)_{V_i}\right|.
  \label{eq.localized_error_indicator}}
In the next section we present an example for the boundary conditions for the
dual problem.

\subsection{Example: Functionals at the
boundary}\label{subsection:example_functional} In the numerical examples in
Section~\ref{section.numexp_fully_implicit} we will consider the 2D Euler
equations. Let $\Gamma_s$ be the solid wall, where we impose the reflecting
boundary condition $v\mdot n= 0$. Thus the flux in normal direction $n$ is given
by
\mmn{
   f_n =p(0, n, 0)^\mathrm{T}.}
The eigenvalues of $\At$ are $\lambda_1=v\cdot n -c =-c$, 
$\lambda_2=\lambda_3=v\cdot n=0$  and $\lambda_4=v\cdot n +c=c$,
and we can compute 
\mmn{
  \Pt_+ &=R\, \mbox{\textnormal{diag}}(0,0,0,1)L,\\
  \Pt  _- &=R\, \mbox{\textnormal{diag}}(1,0,0,0) L.}
and
\mmn{
  \Pt_+\,f_n(U)=\frac{p}2\left(\frac1c, n, \frac{c}{\gamma-1}\right).}
As our functional we choose the space-time integral of the pressure
at the solid wall,
\mmn{
   J(U) = \int_0^T\int_{\Gamma_{s}}p\,dS(x)\,dt.}
If we choose
\mmn{
  \psi_\Gamma =2(0, n,0)^\mathrm{T},}
then $J(U)$ may be rewritten in terms of characteristic projections,
\mmn{
  J(U) =\int_0^T (\Pt_+\,f_n(U),\psi_\Gamma)_{\Gamma_s} dt.}

\replace{}{Note that the functional is in the form suggested in 
Section~\ref{sec.error_representations}. It is a modifcation of lift and drag,
which is given by
\begin{align*}
J_{l/d}(U) = \frac{2}{\rho_\infty |v_\infty|^2 \bar{l}}\int_{\Gamma_{s}}p\, {n}\cdot \psi \,dS(x),
\end{align*}
where $\bar{l}$ denotes a reference length, the subscript $\infty$ indicates
free stream quantities, $\psi$ is given by 
$\psi_d = (\textnormal{cos}(\alpha), \textnormal{sin}(\alpha))^\mathrm{T}$
or $\psi_l = (\textnormal{-sin}(\alpha), \textnormal{cos}(\alpha))^\mathrm{T}$ 
for the drag and lift coefficient, with respect to an angle of attack $\alpha$.}
In the numerical experiments in 
Sections~\ref{section.numexp}--\ref{section.numexp_explicit_implicit}, we will
multiply $\psi_\Gamma$ with an additional weighting function 
(see \eqref{eq.func_p}).

\section{Numerical realization}
\label{section.numexp}

Before we set up our test problem (Section~\ref{section.numexp_setup}) and
present numerical experiments (Sections~\ref{section.numexp_fully_implicit} --
\ref{section.numexp_explicit_implicit}),  we have to specify some details on the
adaptive concept, the grid generation and  the numerical flux evaluation on
locally refined grids with hanging nodes.

\subsection{Adaptive Method in time}

Now we combine the multi\replace{scale}{resolution} based approach introduced in Section~\ref{a4-ms}
and the time adaptive method derived from the space-time splitting of the error
representation to get a space-time adaptive algorithm:
\begin{itemize}
\item solve the primal problem \eqref{eq.euler_weak} on a {\em coarse}
  adaptive spatial grid using uniform \replace{$\cfl$ numbers ($\cfl = 0.8$),} {global timesteps, with maximal ($\cfl = 0.8$),}
\item compute the dual problem \eqref{eq.adjoint} and
  \eqref{eq.adjoint_h_bc} and the space-time-error representation 
  \eqref{eq:error_representation2}. In particular, compute the localized 
  error indicators $\ebkn$ using \eqref{eq.localized_error_indicator}.
\item compute the new adaptive timestep sizes depending on the temporal part
  of the error representation and the $\cfl$ number on the new grid, 
  aiming at an equidistribution of the error,
\item solve the primal problem using the new timestep sizes on a 
  {\em finer} spatial grid.
\end{itemize}

The advantage is, that the first computations of the primal problems and the
dual problem are done on a coarse spatial grid, and therefore have low cost.
These computations provide an initial guess of the timesteps for the computation
on the finer spatial grid. We will restrict the timestep size from below to
$\cfl=0.8$, since smaller timestep sizes only add numerical diffusion to the
scheme and increase the computational cost. Note that all physical effects
already have to be roughly resolved on the coarse grid in order to determine a
reliable guess for the timesteps on the fine grid.

We will deal with some aspects in detail in the numerical examples in Section
\ref{section.numexp_fully_implicit} and  \ref{section.numexp_explicit_implicit}.

\subsection{Multi\replace{scale}{resolution} \replace{analysis}{decomposition} and adaptation in space} \label{a4-ms}

A finite volume discretization is typically working on cell averages. In order
to analyze the local regularity behavior of the data we employ the concept of
biorthogonal wavelets \cite{SM-Carnicer-Dahmen-Pena:96,
Cohen-Daubechies-Feauveau:92}. This approach may be considered as a natural
generalization of Harten's discrete framework  \cite{SM-Harten:96}. The core
ingredients are a hierarchy of nested grids, biorthogonal wavelets and the
multi\replace{scale}{resolution} decomposition. In the following we will only summarize the basic
ideas. For technical details  we refer the reader to the book \cite{Mueller:02}
and \cite{BramkampLambyMueller2004}, respectively.

\begin{figure}[ht]
\vspace*{-5mm}
\begin{minipage}{1.0\columnwidth}
  \begin{minipage}{0.49\columnwidth}
    \vspace*{4mm}
    \begin{center}
      \includegraphics[width=0.9\textwidth]{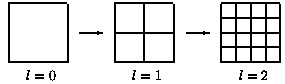}  
    \end{center}
    \vspace*{-3mm}
    \caption{Sequence of nested grids \label{fig.ms1}}
  \end{minipage}
  \begin{minipage}{0.49\columnwidth}
    \vspace*{4mm}
    \begin{center}
      \includegraphics[width=0.9\textwidth]{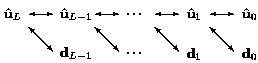}  
    \end{center}
     \caption{multi\replace{scale}{resolution} transformation \label{fig.ms2}}
  \end{minipage}
\end{minipage}
\vspace*{-5mm}
\end{figure}

{\bf Step 1:} {\em Multi\replace{scale}{resolution} \replace{analysis}{decomposition}.}
The fundamental idea is to present the cell averages $\hat{\LMMSvect{u}}_L$
representing the discretized flow field at fixed time level $t^\nn$ on a given
uniform highest level of resolution $l=L$ ({\em reference mesh})  associated
with a given finite volume discretization ({\em reference scheme}) as cell
averages $\hat{\LMMSvect{u}}_0$ on some coarsest level $l=0$. Here the fine
scale information is encoded in arrays of {\em detail coefficients}
$\LMMSvect{d}_l$, $l=0,\ldots,L-1$ of ascending resolution, see Figure
\ref{fig.ms2}.

The multi\replace{scale}{resolution} decomposition is performed on a hierarchy of {\em nested} grids
${\mathcal{G}}_l$ with increasing resolution $l=0,\ldots,L$ determined by dyadic
grid refinement of the logical space,  see Figure \ref{fig.ms1}. 

{\bf Step 2: }{\em Thresholding.}
It can be shown that the detail coefficients become small with increasing
refinement level when the underlying function is locally smooth. This motivates
us to  discard all detail coefficients $d_{l,{k}}$  whose absolute values fall
below a level-dependent  threshold value $\varepsilon_l = 2^{l-L} \varepsilon$ 
in order to compress the original data.  Let ${\mathcal{D}}_{L,\varepsilon}$ be
the set of {\em significant details}. The ideal strategy would be to determine
the threshold value $\varepsilon$ such that the  {\em discretization error} of
the reference scheme,  i.e., the difference between exact solution and reference
scheme, and  the {\em perturbation error}, i.e.,  the difference between the
reference scheme and the adaptive scheme, are balanced, see
\cite{Cohen-Kaber-Mueller-Postel:02}.  

{\bf Step 3: }{\em Prediction and grading.}
Since the flow field evolves in time, grid adaptation is  performed  after each
evolution step to provide  the adaptive grid at the {\em new} time level.  In
order to guarantee the  adaptive scheme to be {\em reliable} in the sense that
no significant future feature of the solution is missed, we have to {\em
predict}  all significant details at the new time level $n+1$ by means of the
details at the {\em old} time level $n$.  Let ${\tilde
{\mathcal{D}}}^\nn_{L,\varepsilon}\supset {\mathcal{D}}^\nn_{L,\varepsilon} \cup
{\mathcal{D}}^\nn_{L,\varepsilon}$ be the  prediction set. The prediction
strategy is detailed in \cite{Cohen-Kaber-Mueller-Postel:02}. In view of the
grid adaptation step this set is additionally inflated  such that it corresponds
to a graded tree, i.e.,  the number of levels between two neighboring cells
differs at most by 1. 

{\bf Step 4: }{\em Grid adaptation.}
By means of the set ${\tilde {\mathcal{D}}}^\nn_{L,\varepsilon}$ a locally
refined grid is determined. For this purpose, we recursively check (proceeding
levelwise from coarse to fine) whether there exists a significant detail on a
cell. If there is one,  then we refine the respective cell. We finally obtain
the  locally refined grid with hanging nodes represented by the index set 
${\mathcal{G}}_{L,\varepsilon}$, see for example Figure \ref{channel-stat2}.

\subsection{Grid generation.}

The computational domain in our test configuration is bounded by curvilinear
boundaries. For this domain we compute a parametric grid mapping 
${x}:[0,1]^2\to\Omega$. Then a hierarchy of Cartesian grids for the parameter
domain is mapped to  a grid hierarchy of curvilinear meshes in  the
computational domain. The grid mapping is realized efficiently by a sparse
B-Spline representation,  cf.~\cite{BramkampLambyMueller2004,SM-Lamby:06}. Then
the locally refined grids are determined by evaluation of this mapping. In our
computations the underlying discretization is always  a hierarchy of curvilinear
grids.

\subsection{Newton method for the nonlinear system.}\label{sec.Newton}

The 2D Euler equations are discretized with the finite volume method on an
adaptive grid. The resulting system of nonlinear equations is discretized by a
Newton method in each timestep. Two delicate aspects of the Newton method are
the choice of the initial value and the choice of the break condition. In our
computations we choose the solution of the old timestep as initial value for the
Newton iterations. The resulting system of linear equations is solved using GMRES \replace{}{\cite{SaadSchultz}}
with an \replace{ILU}{incomplete LU factorization (ILU)} preconditioning, see \replace{}{\cite{BeyWittum}}.  The break condition for the Newton method is
coupled with the threshold value of the multi\replace{scale}{resolution} method: If the defect of the
Newton method is below the threshold of the multi\replace{scale}{resolution} representation, we will
stop the iteration process.

\begin{remark}
There are also other strategies to terminate the Newton iterations.  If we
assume, that in each iteration (both the nonlinear and the linear) the residuals
will drop, we could choose as break condition how many times the residual has to
decrease. Additionally we can set an absolute value of the number of iterations,
if no other breaking condition holds. There are two main problems, which may
occur. The first one is, that in the case where the solution is stationary, no
iteration will decrease the error. The other  problem is, that if the timestep
sizes are large, and the residual is  very large in the beginning and will drop
down fast, the Newton iteration will stop, but will not lead to a good solution,
with a small residuum. In the numerical examples we observe that the coupling
with the multi\replace{scale}{resolution} representation leads to very efficient results. Alternative
strategies to control the timestepping could be based on the residual or the
defect of the Newton-method, the boundary conditions, the $\cfl$ number. For
such a strategy it is obvious that we have many parameters which have to be
chosen and optimized, see \cite{Pollul2008}.
\end{remark}

\subsection{Computation of the dual problem}

The conservative dual problem \eqref{eq.adjoint_w_bc} is a system of $d\mdot m$
conservation laws, where $m$ is the number of equations of the forward problem
and $d$ the number of space dimensions. Since for the backward solver robustness
is more important than accuracy we solve
\eqref{eq.adjoint_w}--\eqref{eq.adjoint_w_bc} with a finite volume method using
Lax-Friedrichs numerical flux and $\cfl<1$.

\section{Setup of the numerical experiment}
\label{section.numexp_setup}

An nonstationary variant of a classical stationary 2D Euler transonic flow,
considered in \cite{RV:81}, is investigated to illustrate the efficiency of the
adaptive method. 

\subsection{Steady state configuration} 
\begin{figure}[ht]
    \begin{center}
          \includegraphics[width=10cm,clip]{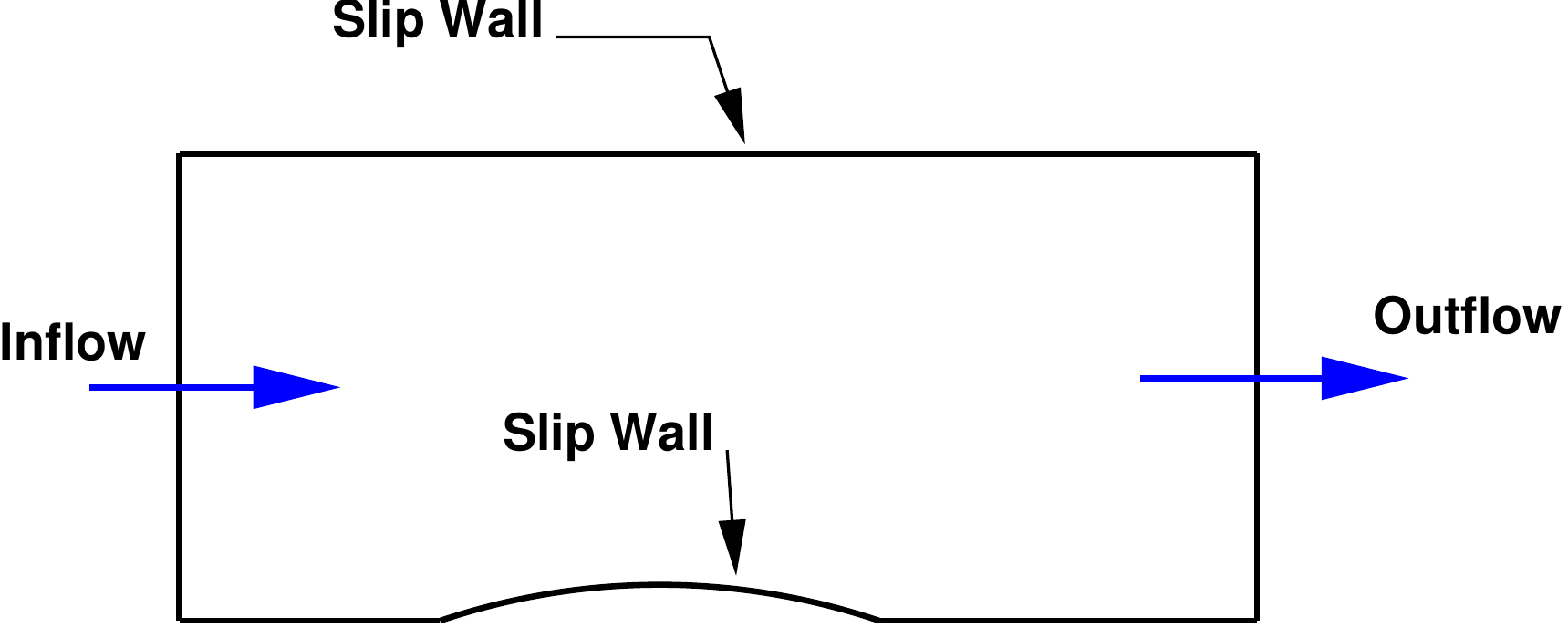}
    \end{center}
    %\vspace*{-7mm}
  \caption{Circular arc bump configuration of the computational domain $\Omega$. \label{channel-config}}
\end{figure}

First we consider the classical setup in the
stationary case. The computational domain is a channel of $3m$ length and $2m$
height with an arc bump of $l = 1m$ secant length and $h = 0.024m$ height cut
out, see Figure \ref{channel-config}.  At the inflow boundary, the Mach number
is 0.85 and a homogeneous flow field characterized by the free-stream quantities
is imposed. At the outflow boundary, characteristic boundary conditions are
used. We apply slip boundary conditions across the solid walls, i.e., the normal
velocity is set to zero. In the numerical examples in Section
\ref{section.numexp_fully_implicit} the height of the channel is $2m$ and the
length $6m$.

\begin{figure}[ht]
    \begin{center}\vspace{-5cm}
          \includegraphics[width=\textwidth]{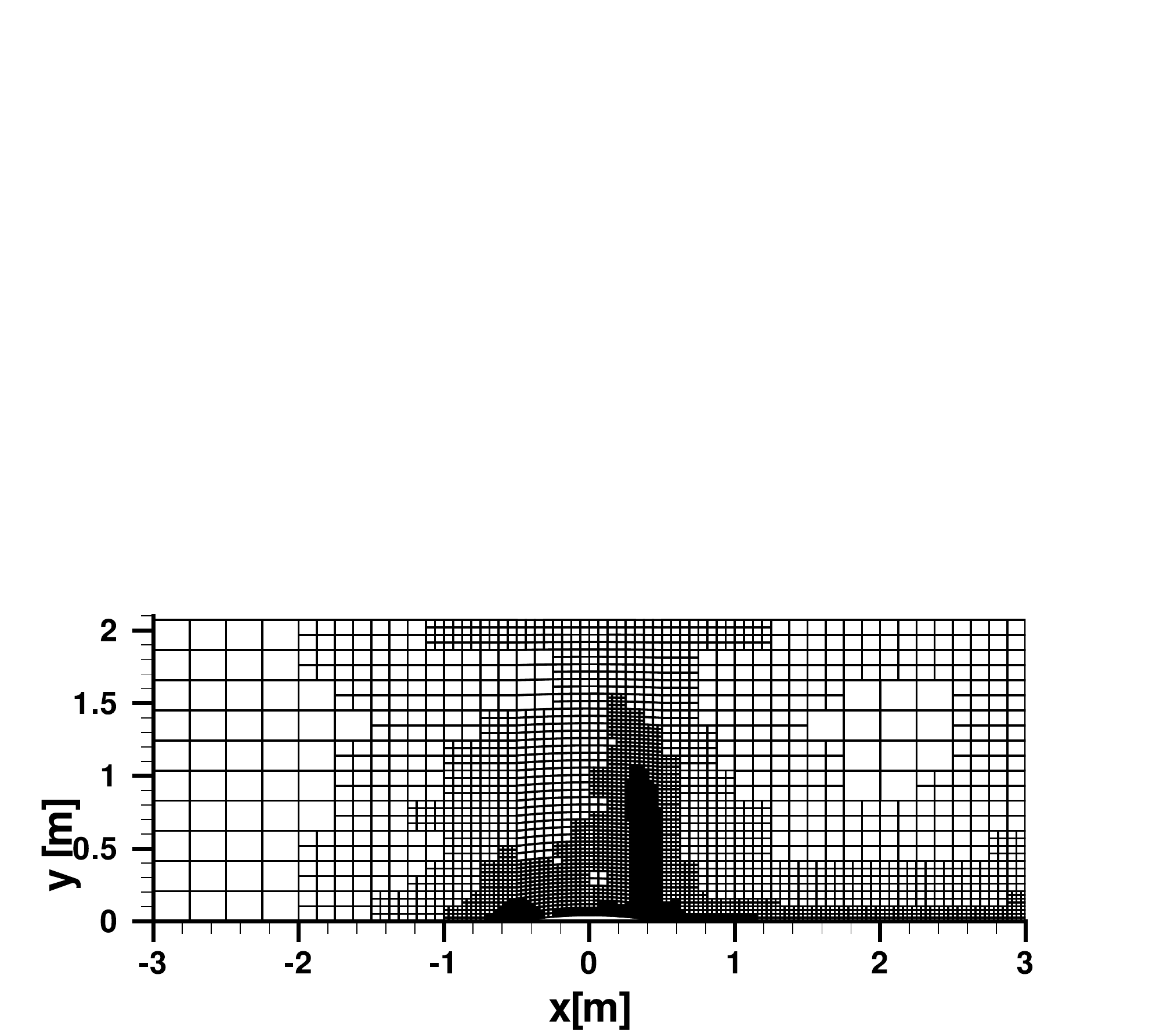}
    \end{center}
    %\vspace*{-7mm}
  \caption{Adaptive grid $L=5$, to steady state solution in Figure \ref{channel-stat1} of the circular arc bump configuration.  \label{channel-statgrid1}}
\end{figure}

\begin{figure}[ht]
    \begin{center}\vspace{-5cm}
          \includegraphics[width=\textwidth]{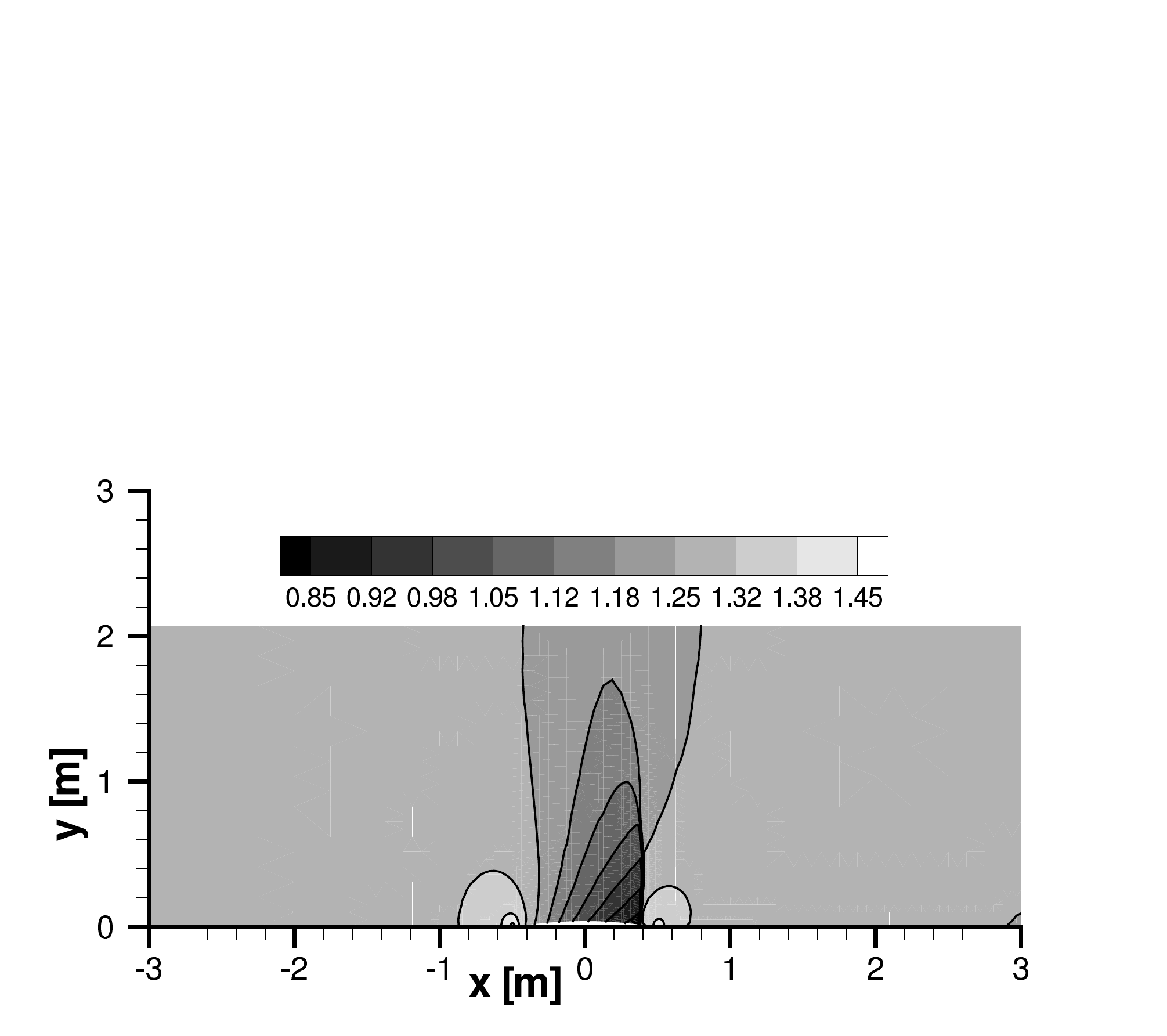}
    \end{center}
    %\vspace*{-7mm}
  \caption{Steady state solution of the circular arc bump configuration: Isolines of the density, $L=5$.  \label{channel-stat1}}
\end{figure}

The threshold value in the grid adaptation step %for the multi\replace{scale}{resolution} \replace{analysis}{decomposition}  is
$\varepsilon = 1 \times 10^{-3}$ and computations are done on adaptive grids
with finest level $L = 2$ and $L=5$ respectively. In general, a smaller
threshold value results in more grid refinement whereas a larger value gives
locally coarser grids. 

In the stationary case at Mach 0.85 there is a compression shock separating a
supersonic and  a subsonic domain. %  The shock wave is sharply captured  and
the stagnation areas are highly resolved, see Figures \ref{channel-stat1} and
\ref{channel-stat2}.  

We will use this steady state solution as initial data for the nonstationary test
case.

\subsection{nonstationary test case} 
Now we define our nonstationary test case
prescribing a time-dependent perturbation coming in at the inflow boundary.
First we keep the boundary conditions fixed, and prescribe the corresponding
stationary solution as initial data. Then we introduce, for a short time period
$[t_b, t_e]$, a perturbation $\alpha$ of the pressure at the left boundary, see
\eqref{eq.pertubation}. We will impose two perturbations of the inflow boundary
conditions, at time $t_b^1= 0.004 s$ until $t_e^1= 0.005 s$ and at $t_b^2= 0.022
s$ until $t_e^2= 0.023 s$. These perturbations increase and decrease in a short
time period of $\tau = 0.00005 s$. The first perturbation is about 20 percent of
the pressure at the inflow boundary and the second 2 percent. The perturbations
imposed move through the domain and leave it at the right boundary. Then the
solution is stationary again. The total time is $t= 0.029s$.

The perturbations are given by:
\mmn{
  w_p^i (t) = \left\{ \begin{array}{ccc} \left(\frac{t-t_b^i}{\tau}\right)^2& 
  \mathrm{for}&t_b^i<t\leq t_b^i+\tau \\ 
  1& \mathrm{for}&t_b^i+\tau<t \leq t_e^i-\tau\\
  \left(\frac{t-t_e^i}{\tau}\right)^2& \mathrm{for}&t_e^i-\tau<t
  \leq t_e^i \end{array}\right.,}
\mm{\label{eq.pertubation}
  p_{in}(t) = p_\infty w_p(t) = p_\infty
  \left\{\begin{array}{ccc}
  1& \mathrm{for}&t \leq t_b^1\\
  1+\alpha^1 w_p^1 (t)& \mathrm{for}&t_b^1<t \leq t_e^1\\
  1& \mathrm{for}&t_e^1<t \leq t_b^2\\
  1+\alpha^2 w_p^2 (t)& \mathrm{for}&t_b^2<t\leq t_e^2\\
  1& \mathrm{for}&t_e^2<t
  \end{array}\right.}
with perturbation parameters listed in Table \ref{table.pert}.
\begin{table}[ht] 
  \begin{center}
  \begin{tabular}{c|c|c|c|c}
    $i$&$\alpha^i$ &$t_b^i [s]$ & $t_e^i [s]$ & $\tau^i [s]$ \\ \hline
    1&0.2	&  0.004&  0.005&    0.00005  \\  \hline
    2&0.02	&  0.022&  0.023&    0.00005  \\ 
  \end{tabular}\\[1ex] 
  \end{center}
  \caption{Parameters of the perturbations $w_p^i, i=1,2$ at the left boundary to
    equation \eqref{eq.pertubation}.}
  \label{table.pert} 
\end{table} 
\begin{figure}
  \begin{center}
  \includegraphics[width=\textwidth]{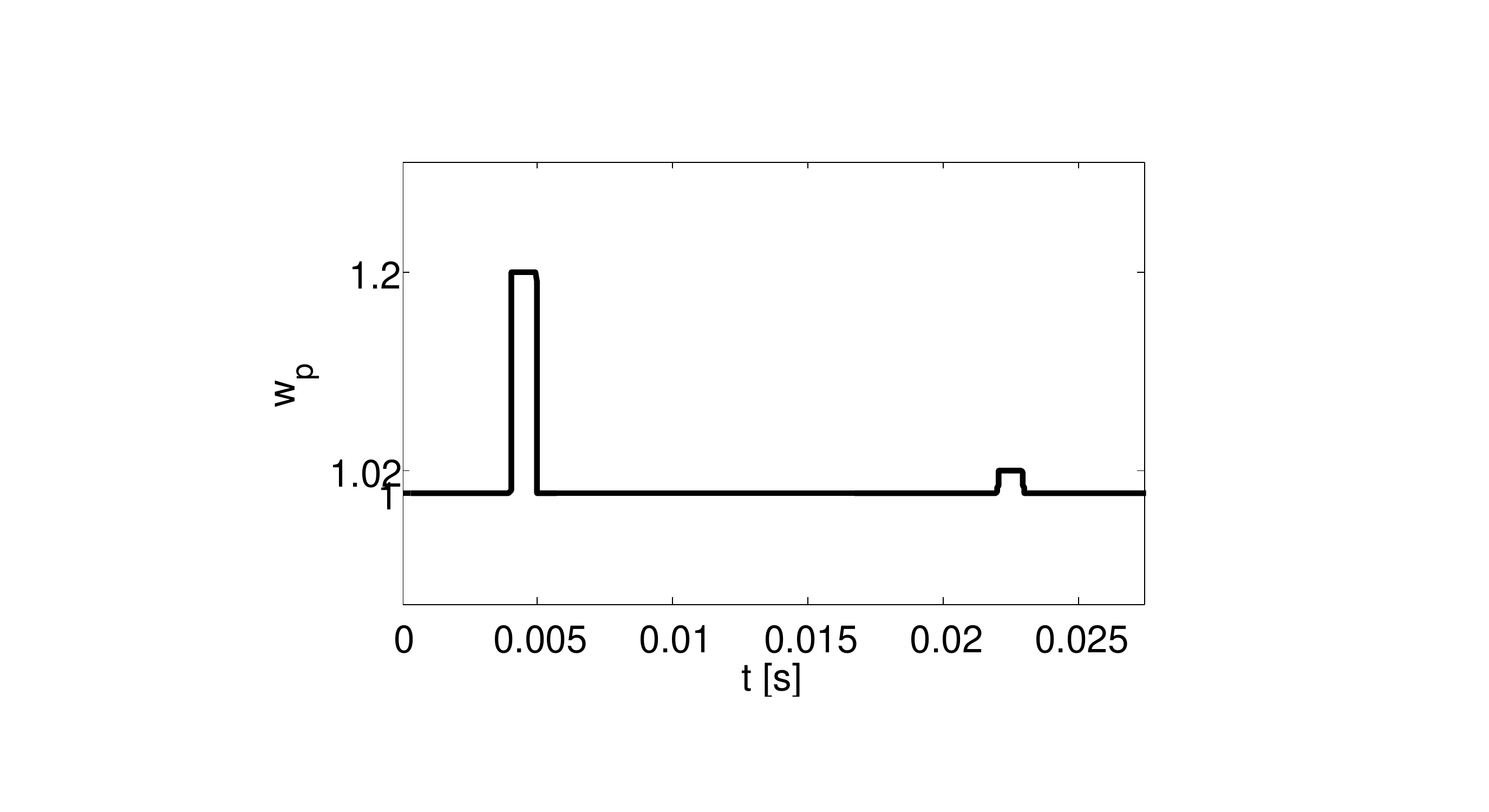}
  \vspace*{-15mm}
  \caption{Weighting function of the perturbation. $p_{in}(t) = p_\infty w_p(t)$} \label{fig.pertub_p}
  \end{center}
\end{figure}
The first computation is done on an adaptive grid with finest level $L=2$. We
also compute the dual solution and the error representation on this level. Using
the time-space-split error representation \eqref{eq:error_representation2} we
derive a new timestep distribution aiming at an equidistribution of the error.
Finally this is modified by imposing a $\cfl$ restriction from below.

We aim to equidistribute the error and prescribe a tolerance
$\;\tol(5)=2^{-3}\bar\eta_k^{ref}\,$, where $\bar\eta_k^{ref}$ is the temporal
error from the computation on level $L=2$.

For this set-up we will show that the adaptive spatial refinement together with
the time-adaptive method will lead to an efficient computation.  The multi\replace{scale}{resolution}
method provides a well-adapted spatial representation of the solution, and the
dual solution will detect time-domains where the solution is stationary. In
these domains, the equidistribution strategy will choose large timesteps.
\begin{figure}
  \begin{center}
  \includegraphics[width=0.49\textwidth]{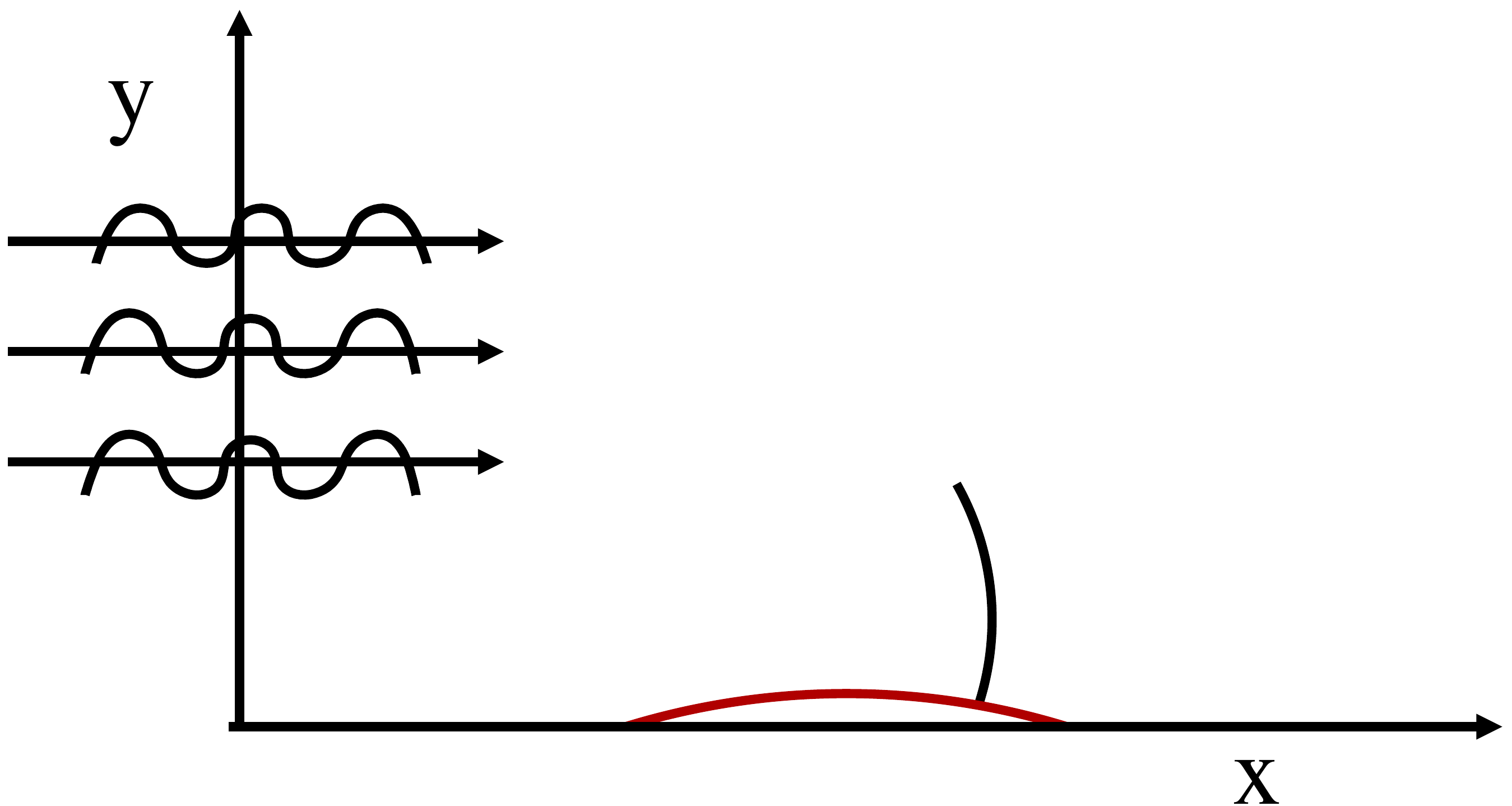}
  \includegraphics[width=0.49\textwidth]{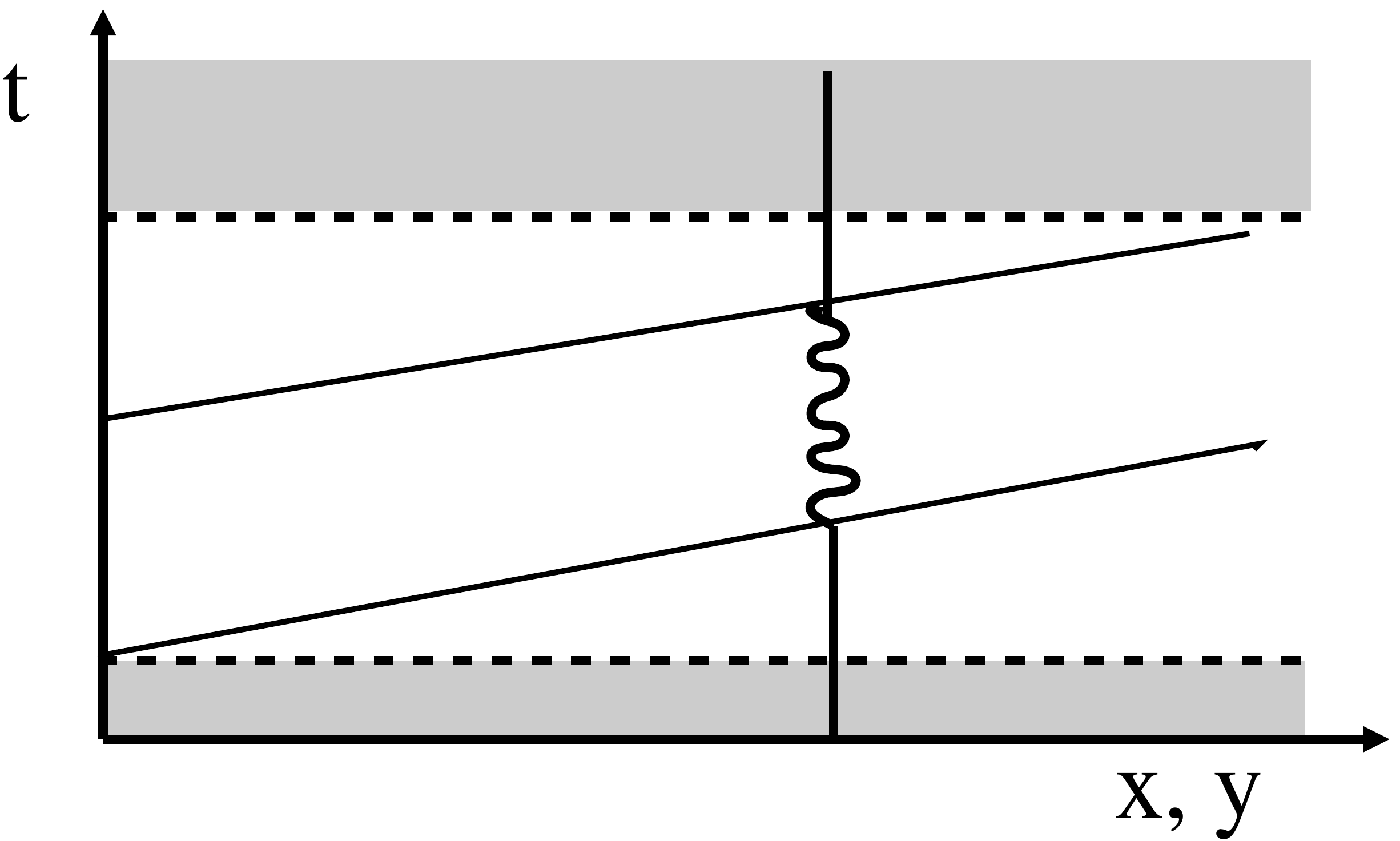}
  \caption{Schematic illustration of the nonstationary test case.
    Left: perturbation at the left boundary. Right: nonstationary (white)
    and stationary (grey) time domains of the solution}
  \end{center}
\end{figure}
\subsection{Target functional}
Now we set up the target functional. The functional $J(U)$ is chosen as a
weighted average of the normal force component exerted on the bump and at the
boundaries before and behind the bump:
\begin{align}\label{eq.func_p}
  J(U) =\sum_{i=1}^7 \int_0 ^T \int_{\kappa_{i}} p \,\psi_{i} (x,y) ds
\end{align}
with
\begin{align*}
  \kappa_{i} &= \left\{
  (x,y)\in \Gamma \,:\, x\in [x_i-0.25, x_i+0.25]\right\}\\
  \psi_{i}(x) &= (x-(x_i-0.25))^2(x+(x_i+0.25))^2/0.25^4,\quad x\in \kappa_i. 
\end{align*}

Here $\Gamma$ is only the bottom part of $\Gamma$. In all computations presented
in this section the functional \eqref{eq.func_p} is chosen, which is the
pressure averaged at several points at the bump in front and behind the bump.
The $x$-coordinates of these points at the bottom are $x_i$ = -3, -2, -1, 0, 1,
2, 3. At each of these points $x_i$ a smooth function $\psi_{i}$ is given with
support $x_i-0.25, x_i+0.25$. This functional measures the pressure locally.

\begin{figure}[ht]
  \begin{center}\vspace*{-35mm}
  \includegraphics[width=0.6\textwidth]{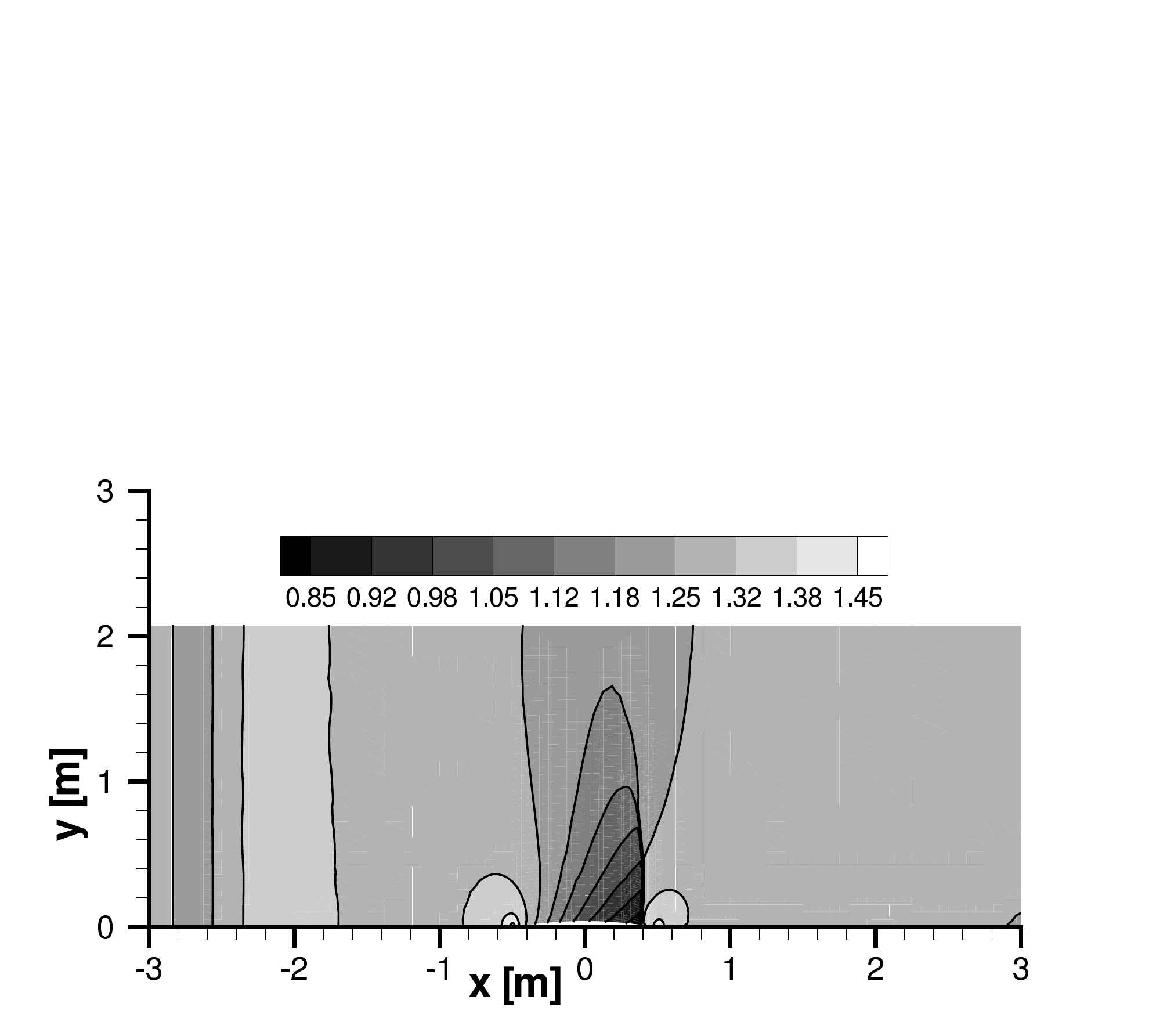}\vspace*{-30mm}
  \includegraphics[width=0.58\textwidth]{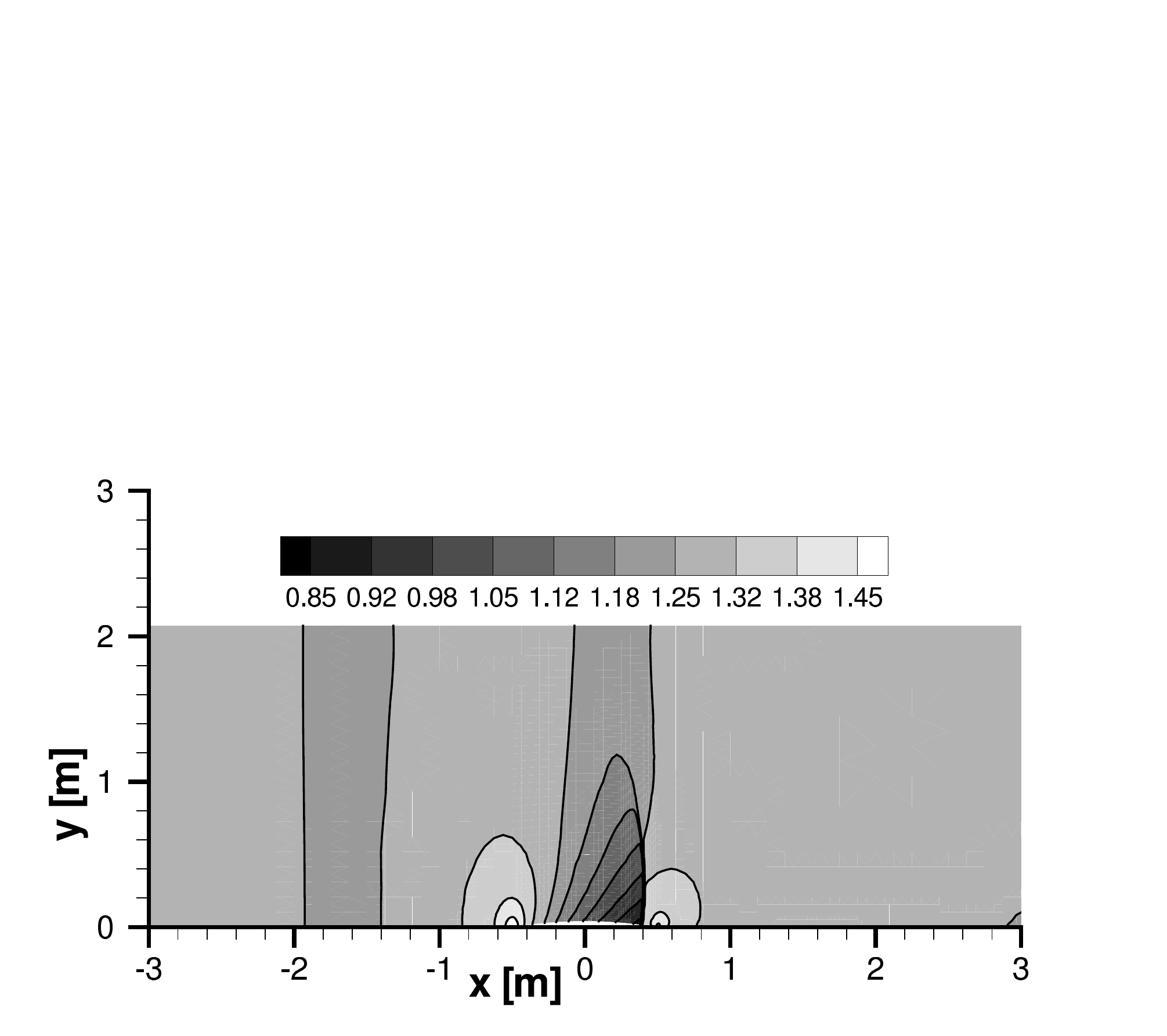}\vspace*{-30mm}
  \includegraphics[width=0.58\textwidth]{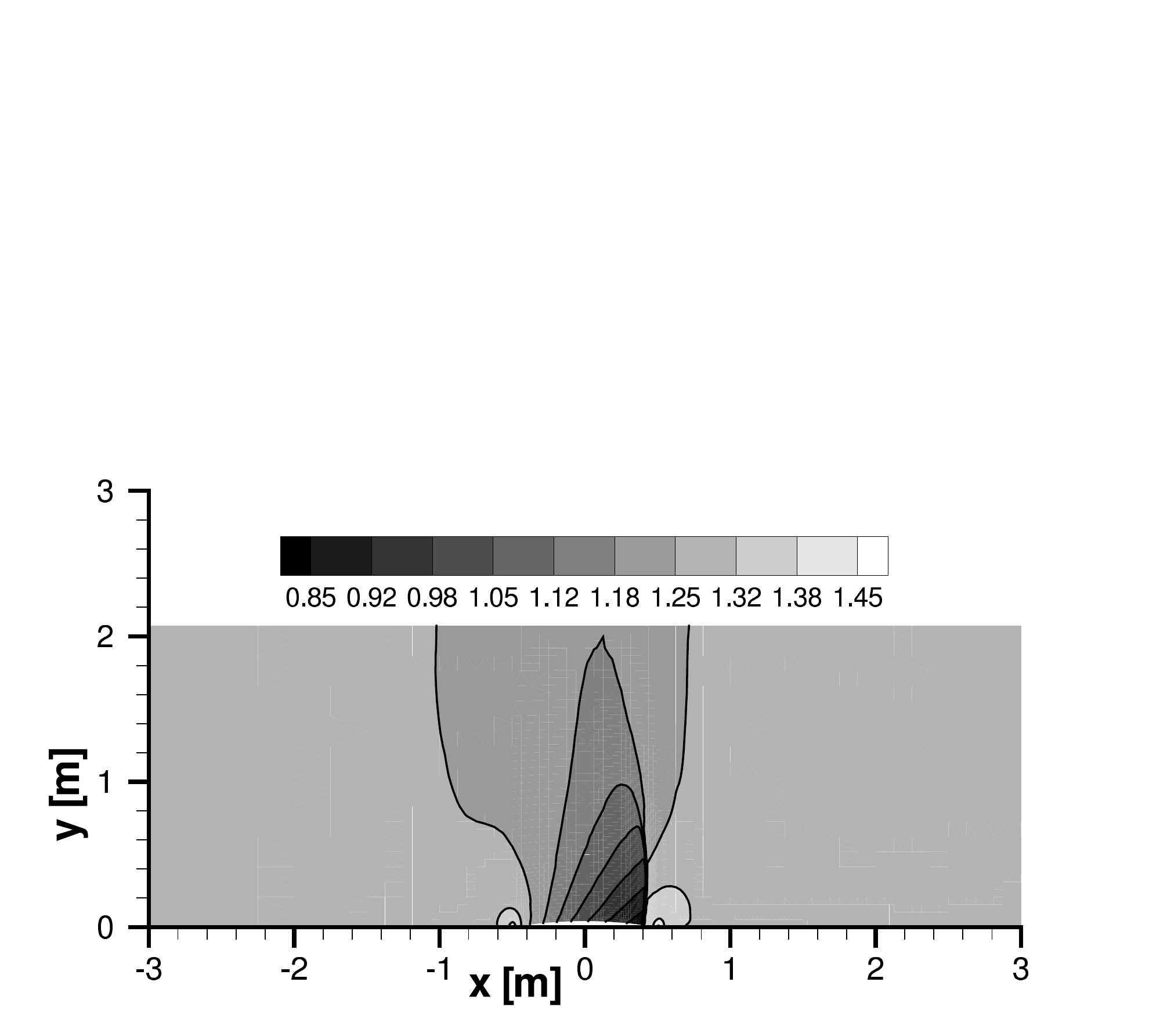}\vspace*{-30mm}
  \includegraphics[width=0.58\textwidth]{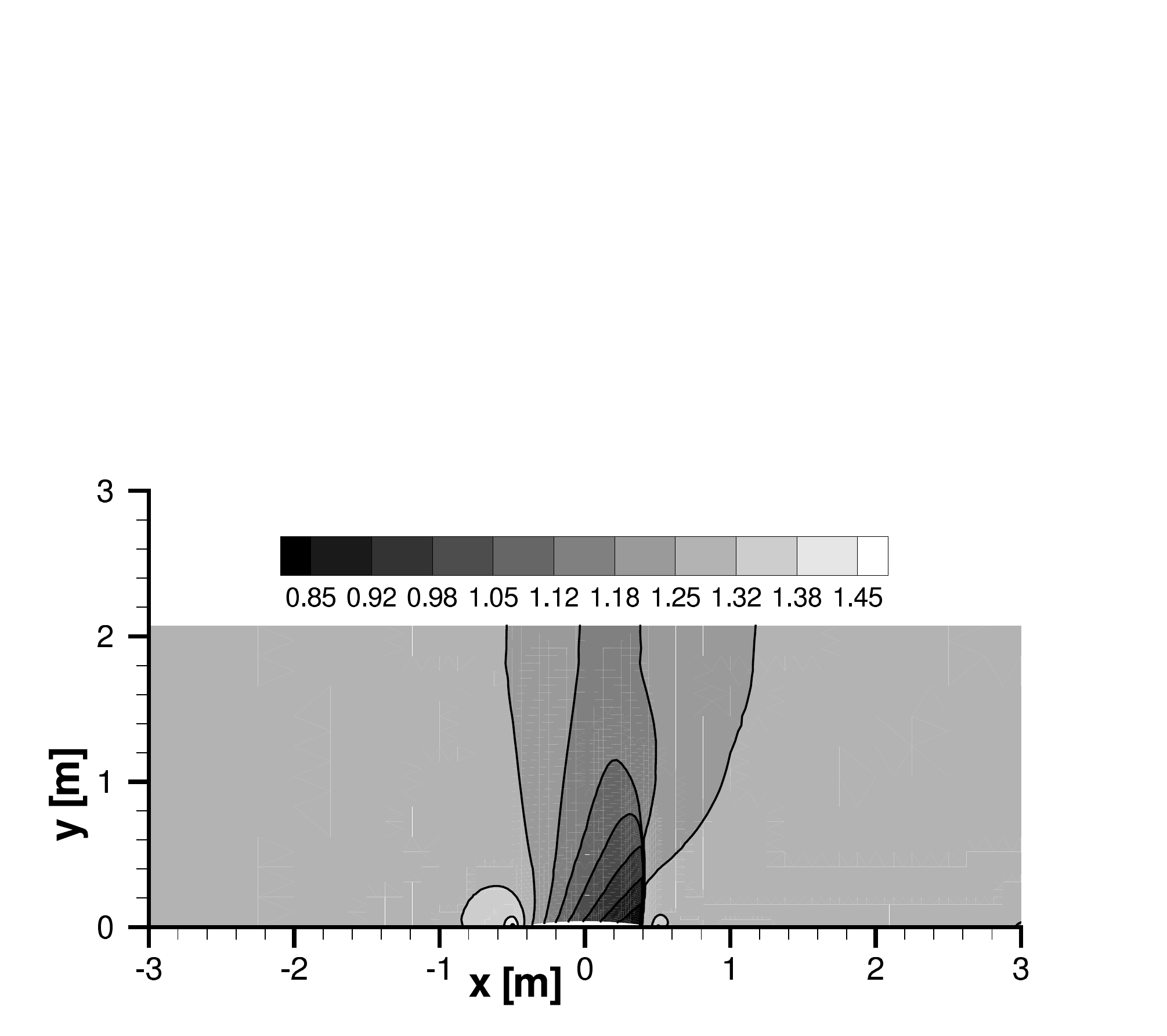}\vspace*{-30mm}
  \includegraphics[width=0.58\textwidth]{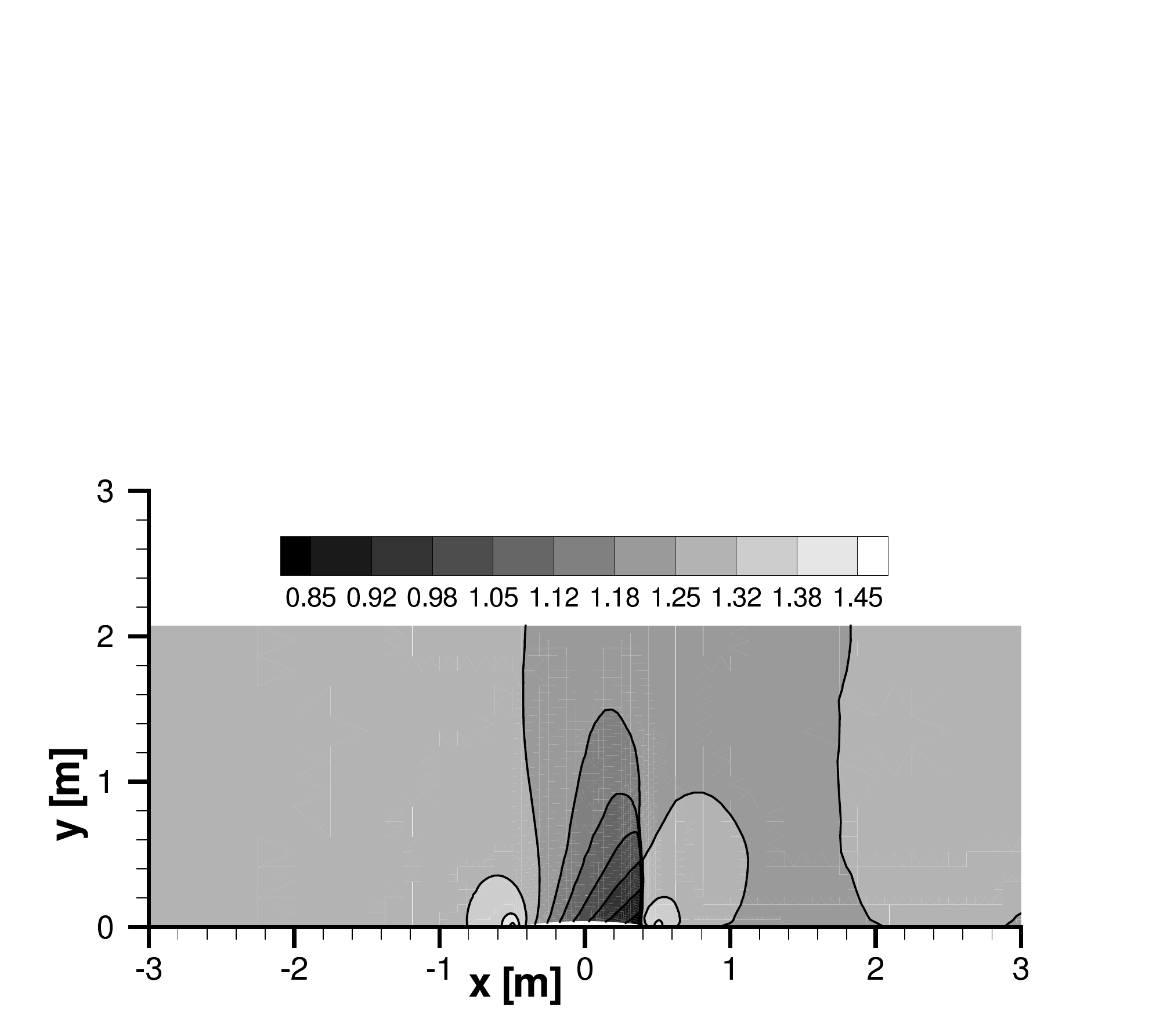}\vspace*{-5mm}
  \end{center}
  %\vspace*{-7mm}
  \caption{nonstationary solution of the circular arc bump configuration, uniform timestep $\cfl =1$, level $L=5$, isolines of the density, perturbation entering on the left and leaving on the right side of the computational domain, from top to bottom: $t$ = 0.0057s, 0.00912s, 0.01254s, 0.01596s, 0.01938s.}\label{channel-stat3}
\vspace*{-20mm}
\end{figure}

\clearpage

\begin{figure}[ht]
  \begin{center}\vspace*{-35mm}
  \includegraphics[width=0.6\textwidth]{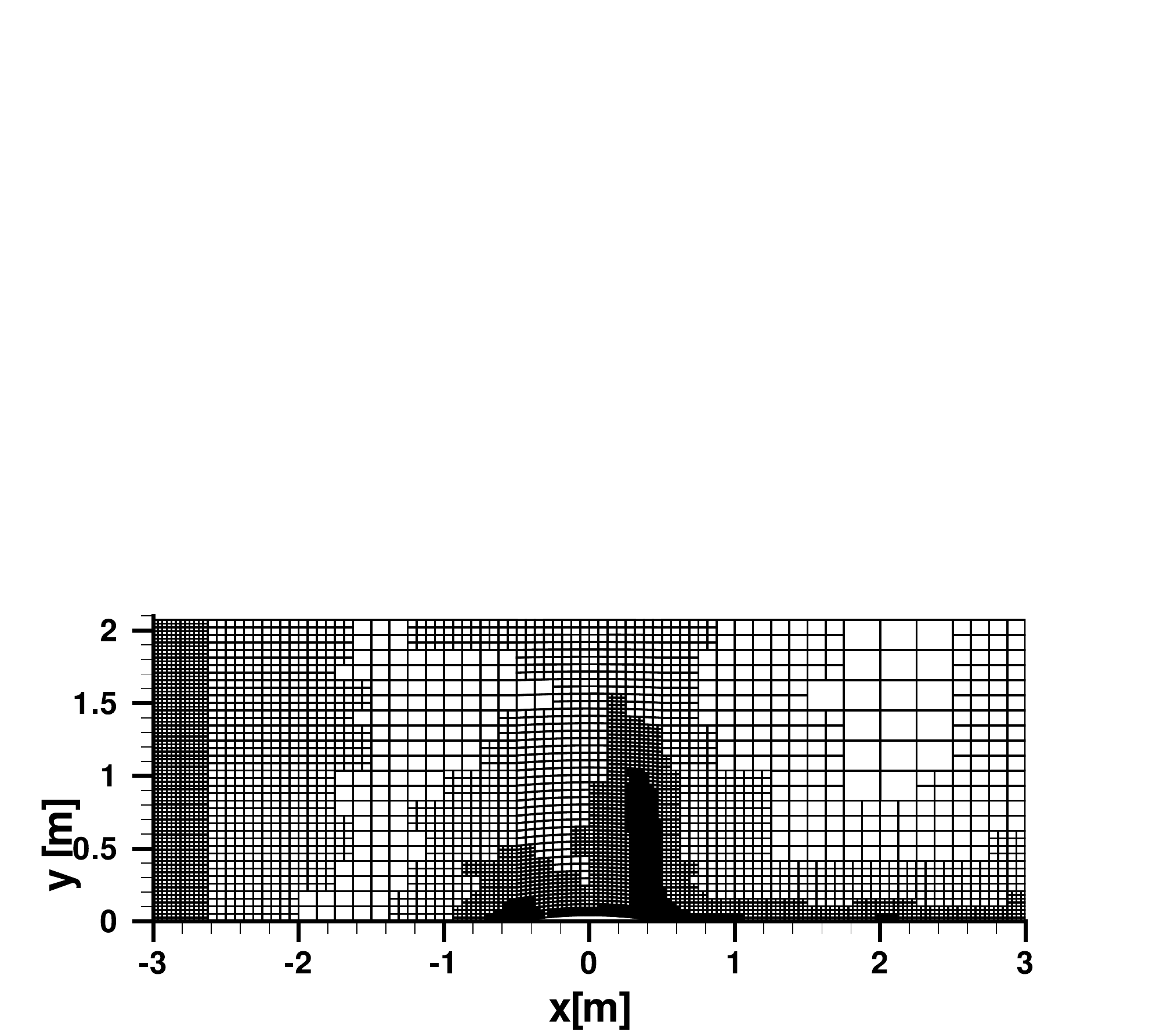}\vspace*{-35mm}
  \includegraphics[width=0.6\textwidth]{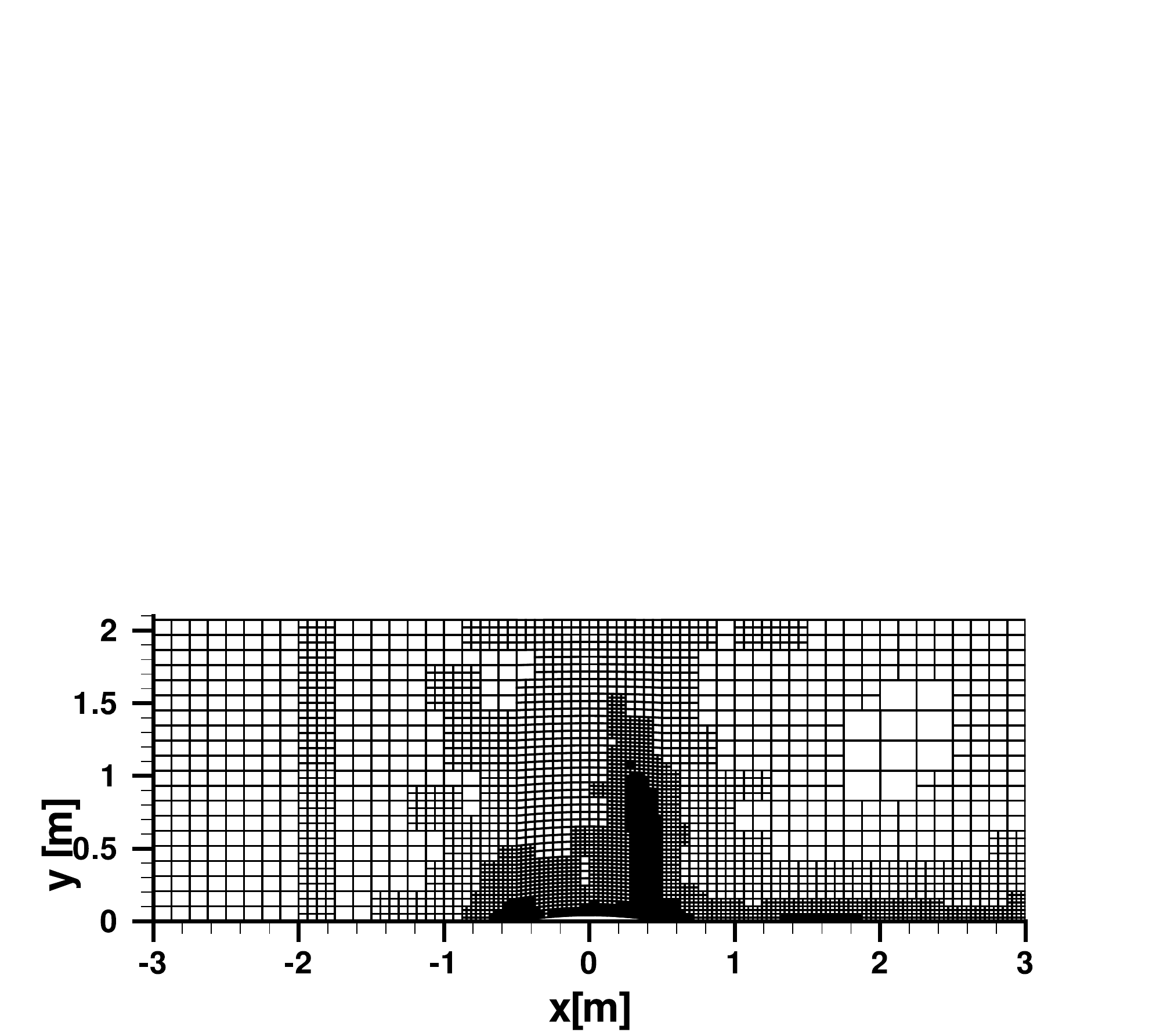}\vspace*{-35mm}
  \includegraphics[width=0.6\textwidth]{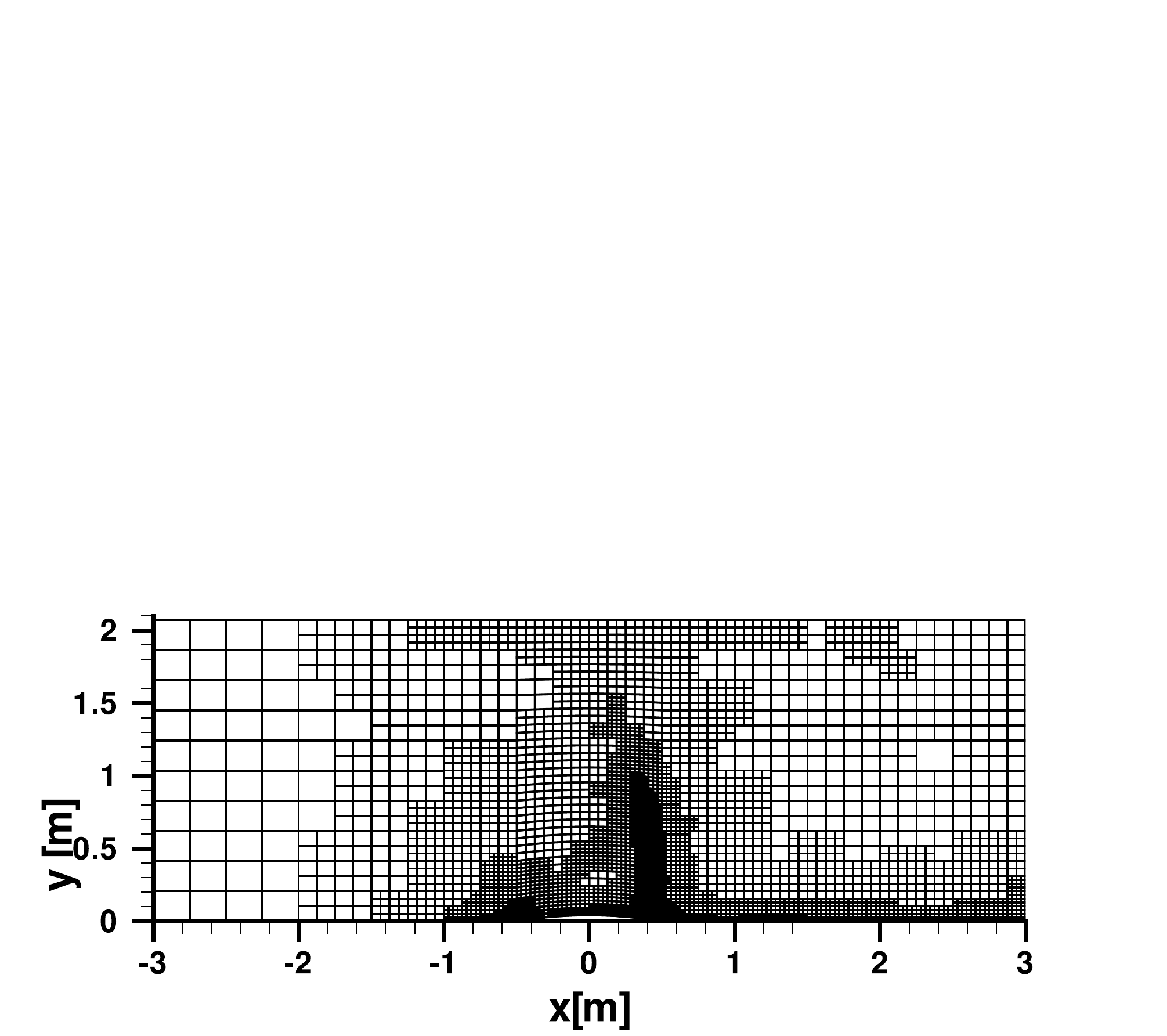}\vspace*{-35mm}
  \includegraphics[width=0.6\textwidth]{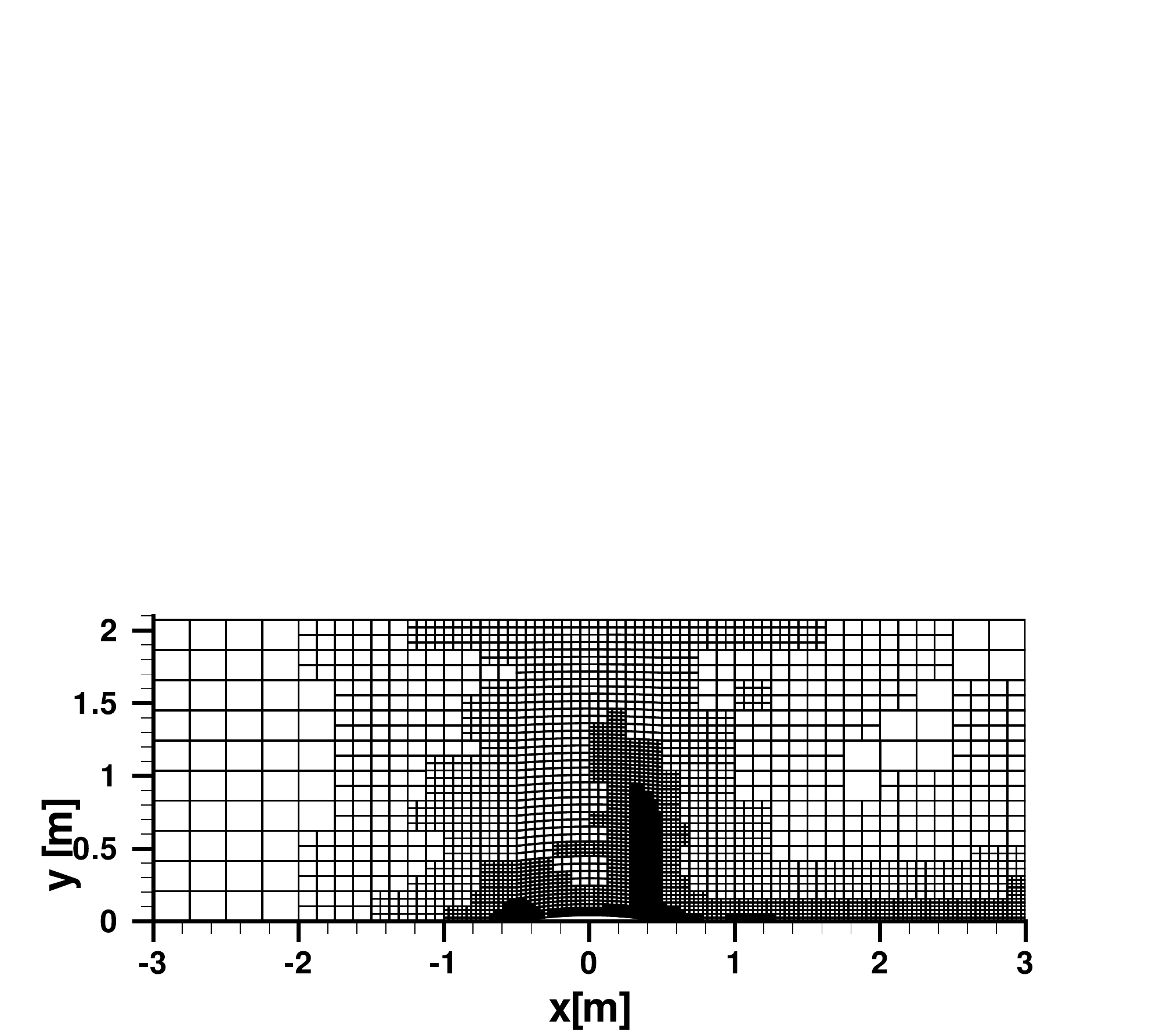}\vspace*{-35mm}
  \includegraphics[width=0.6\textwidth]{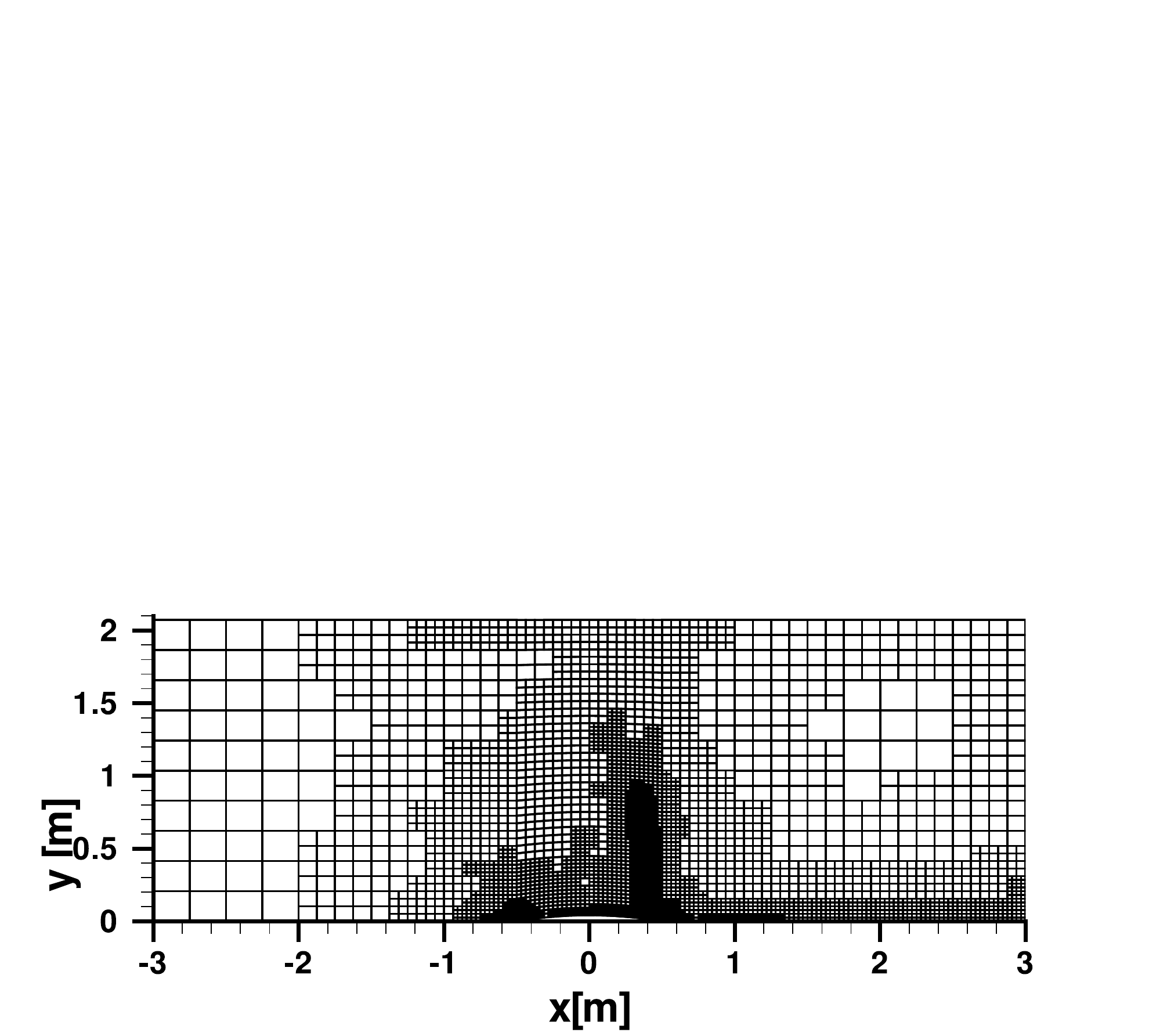}%\vspace*{-5mm}
  \end{center}
  %\vspace*{-7mm}
  \caption{nonstationary solution of the circular arc bump configuration:
    adaptive grid $L=5$ for the computation Figure \ref{channel-stat3},
    perturbation entering on the left side and leaving on the right
    side of the computational domain, from top to bottom:
    $t$ = 0.0057s, 0.00912s, 0.01254s, 0.01596s, 0.01938s.}\label{channel-stat2}
\vspace*{-5mm}\end{figure}

\begin{remark}
(i) Our experience is that the averaged pressure is only computed accurately if
the whole flow field is well resolved. Therefore our functional-based time
adaptation together with the multi\replace{scale}{resolution} spatial adaptation seems to yield a
reliable global accuracy both in space and time. 

(ii) Note that the functional $J(U)$ is never evaluated because $U$ is not known
at all. Instead we compute the localized indicator 
\eqref{eq.localized_error_indicator}, which involves solving the conservative
linearized adjoint problem \eqref{eq.adjoint_w} -- \eqref{eq.adjoint_w_bc}.
\end{remark}

In Figure \ref{channel-stat3} we show a time-sequence of the nonstationary test 
case computed with uniform $\cfl=1$ on an adaptive grid with finest level $L=5$. 
In Figure \ref{channel-stat2} the corresponding adaptive grids are presented. 
Note that the perturbation entering at the left boundary and moving through the 
boundary is resolved very well.

\section{Fully implicit computational results}
\label{section.numexp_fully_implicit}

\subsection{Numerical strategies}

We will present and compare three strategies to demonstrate the efficiency of
our space-time adaptive method.

\begin{enumerate}
\item {\bf Adaptive timesteps via adjoint indicator.}\\
  The first strategy is the one we proposed in \cite{SteinerNoelle}:
  We first compute a forward solution on a coarse grid ($L=2$) and solve 
  the adjoint problem on the adaptive grid of the forward solution.
  Then we use the information of the error representation based on the dual solution to determine a new sequence of timesteps. This sequence is used in the computation of the forward solution on a grid with finest level $L=5$, where we additionally restrict the $\cfl$ number from below. 
\item {\bf Adaptive timesteps via ad hoc indicator.}\\
  In the second approach we compare our indicator with ad hoc indicators,
  which do not require the solution of an adjoint problem. To get these
  indicators we first do a computation on a coarse grid ($L=2$), and
  compute residuals in time of the approximate solution.
\item {\bf Uniform timesteps.}\\
  In a third approach we will set-up uniform timestep distributions
  with the same number of timesteps as in the adaptive case of strategy one.
  We compare the results with our timestep distribution and a uniform in time
  computation with $\cfl =1$ and $\cfl =10$. 
\end{enumerate}

%\clearpage

We want to compare these strategies with respect to the following main aspects:
\begin{itemize}
\item What is the quality and what are the costs determining the adaptive
  timestep sequence from computations on the coarse grid?
\item Is the predicted adaptive timestep sequence well-adapted to the solution
  on the fine grid?
\item Do ad hoc indicators without computing a dual solution lead to comparable
  results?
\item How is the solution affected if we use uniform timesteps larger than
  the predicted adaptive timestep sequence?
\end{itemize}

In order to quantify the results we have to compare with a reference solution.
Since the exact solution is not available we perform a computation with $L=5$
refinement levels using implicit timestepping with $\cfl= 1$. This is a very
expensive approximation for the nonstationary case. For all of the above issues
we will discuss the quality of the solution, the computational costs (time and
memory) and the efficiency.

\subsection {Adaptive timesteps via adjoint indicator}\label{sec.numadjoint}

Now we use the error representation on finest level $L=2$ for a new time
adaptive computation on finest level $L=5$. 

\begin{figure}[ht]  %[hbtp]
  \begin{center} {\vspace{-3cm}
  {\includegraphics[width=0.9\textwidth]{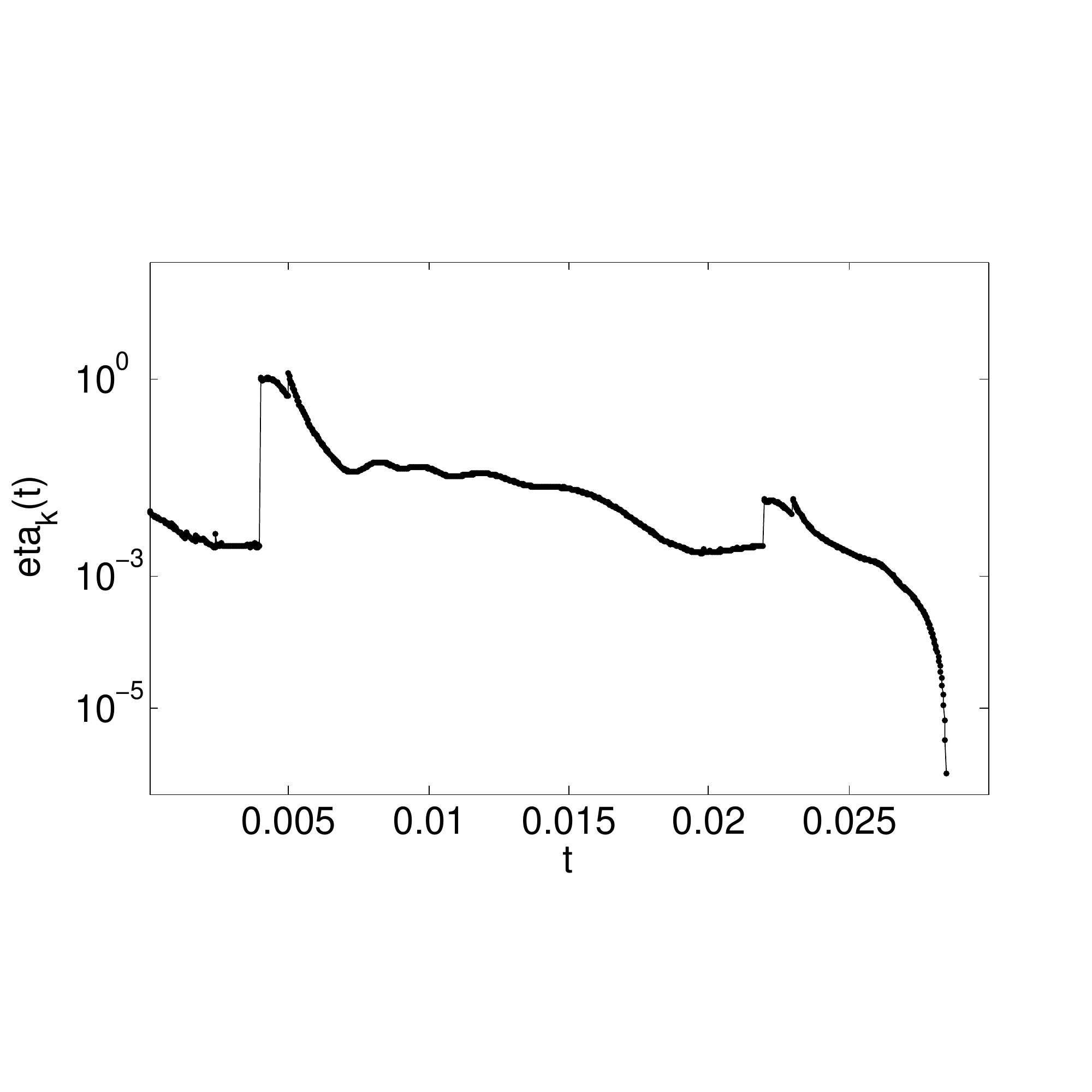}}\vspace{-4cm}
  {\includegraphics[width=0.9\textwidth]{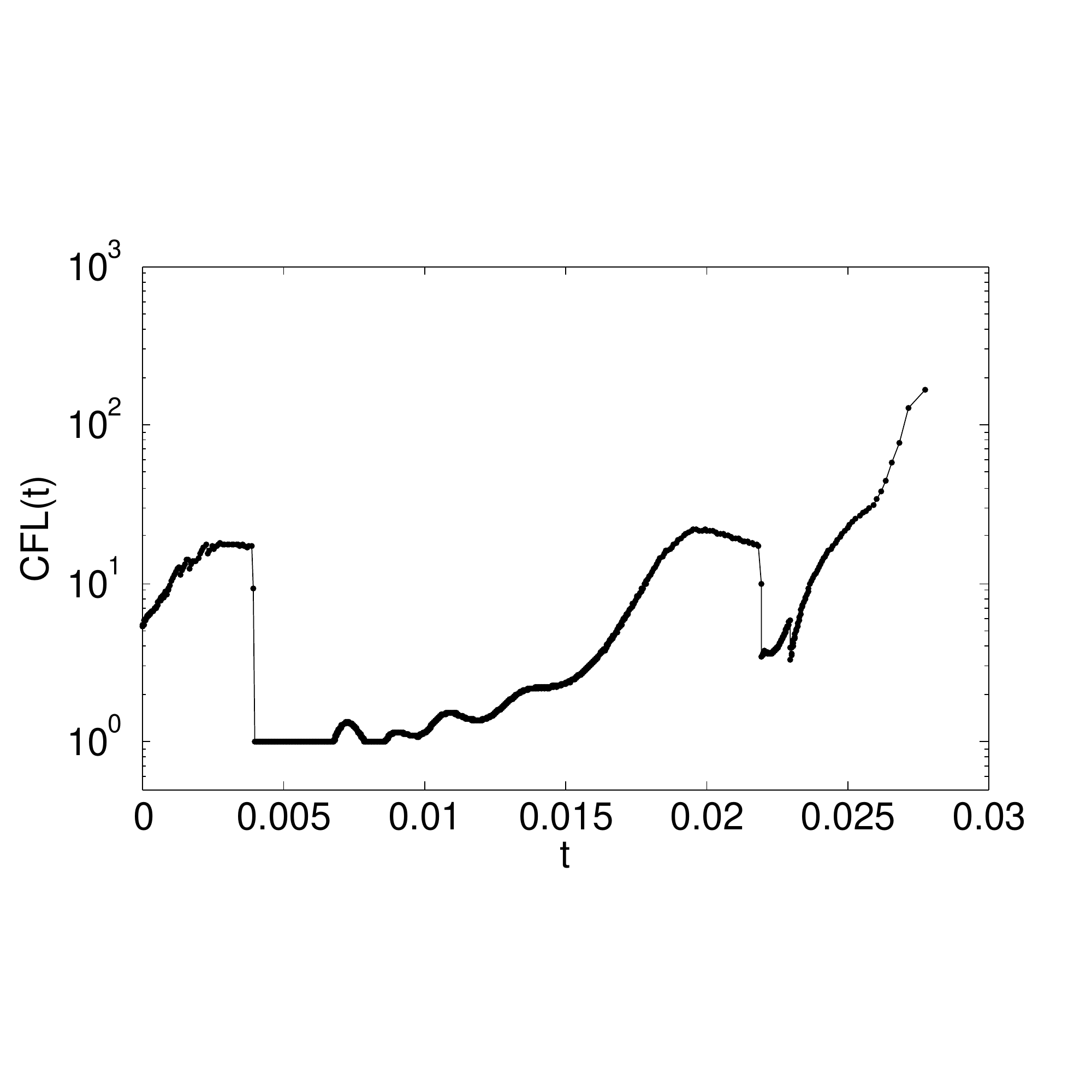}}\vspace{-2.5cm}}
  \end{center}
  \caption{Time component of the error representation $\ebkn$(top) and new
    timesteps with $\cfl$ restriction from below $\cfl(t_n)$ (bottom).
  \label{fig.vgl_euler2d_1}}
%\vspace{-3cm}
\end{figure}

The first computation is done on a mesh with finest level $L=2$. We compute
until time $T= 0.0285 s$, which takes 1000 timesteps with $\cfl=1$. We use the
results of the error representation of this computation to compute a new
timestep distribution. The forward problem takes $329 s$ and the dual problem
including the evaluation of the error representation $619 s$ on an Opteron 8220
processor at 2.86 GHz. The total computational costs are $948 s$, and in memory
we have to save 1000 solutions (each timestep) of the forward problem which
corresponds to 48 MB (total). This gives us a new sequence of adaptive timesteps
for the computation of level $L$ = 5. The error indicator and the new timesteps
are presented in Figure \ref{fig.vgl_euler2d_1}. In time intervals where the
solution is stationary, i.e., at the beginning, and after the perturbations have
left the computational domain, the timesteps are large. In time intervals where
the solution is nonstationary we get well-adapted small timesteps.

Then we use the adaptive timestep sequence for a computation on level $L=5$ and
compare it with a uniform in time computation using $\cfl= 1$. The uniform
computation needs 8000 timesteps and the computational time is $21070 s$. The
time adaptive solution is computed with 2379 timesteps and this computation
takes 9142$s$.

\begin{figure}[t]
\begin{center}\vspace*{-2.5cm}
\includegraphics[width=0.6\textwidth]{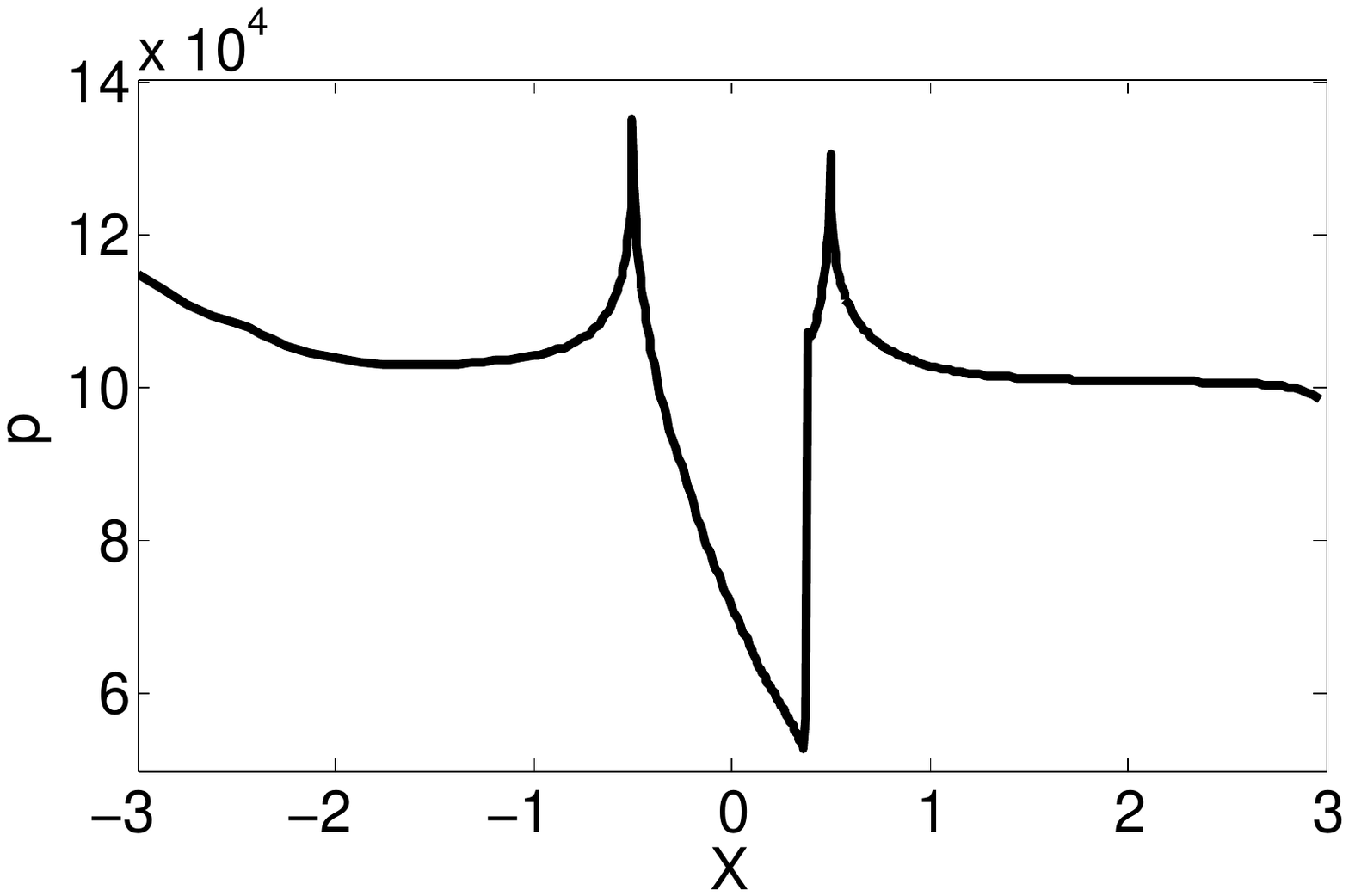}\vspace{-6cm}
\includegraphics[width=0.6\textwidth]{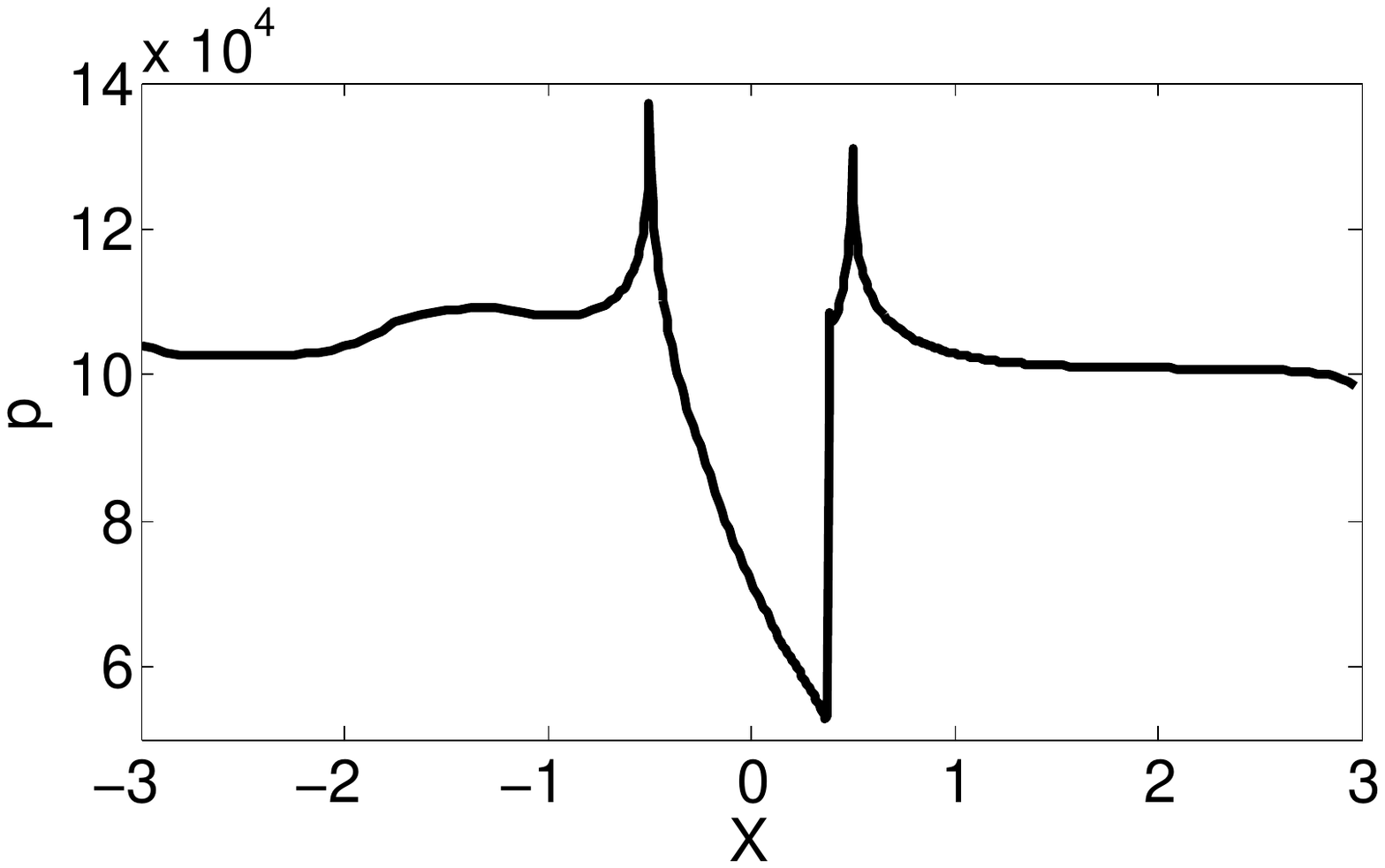}\vspace{-6cm}
\includegraphics[width=0.6\textwidth]{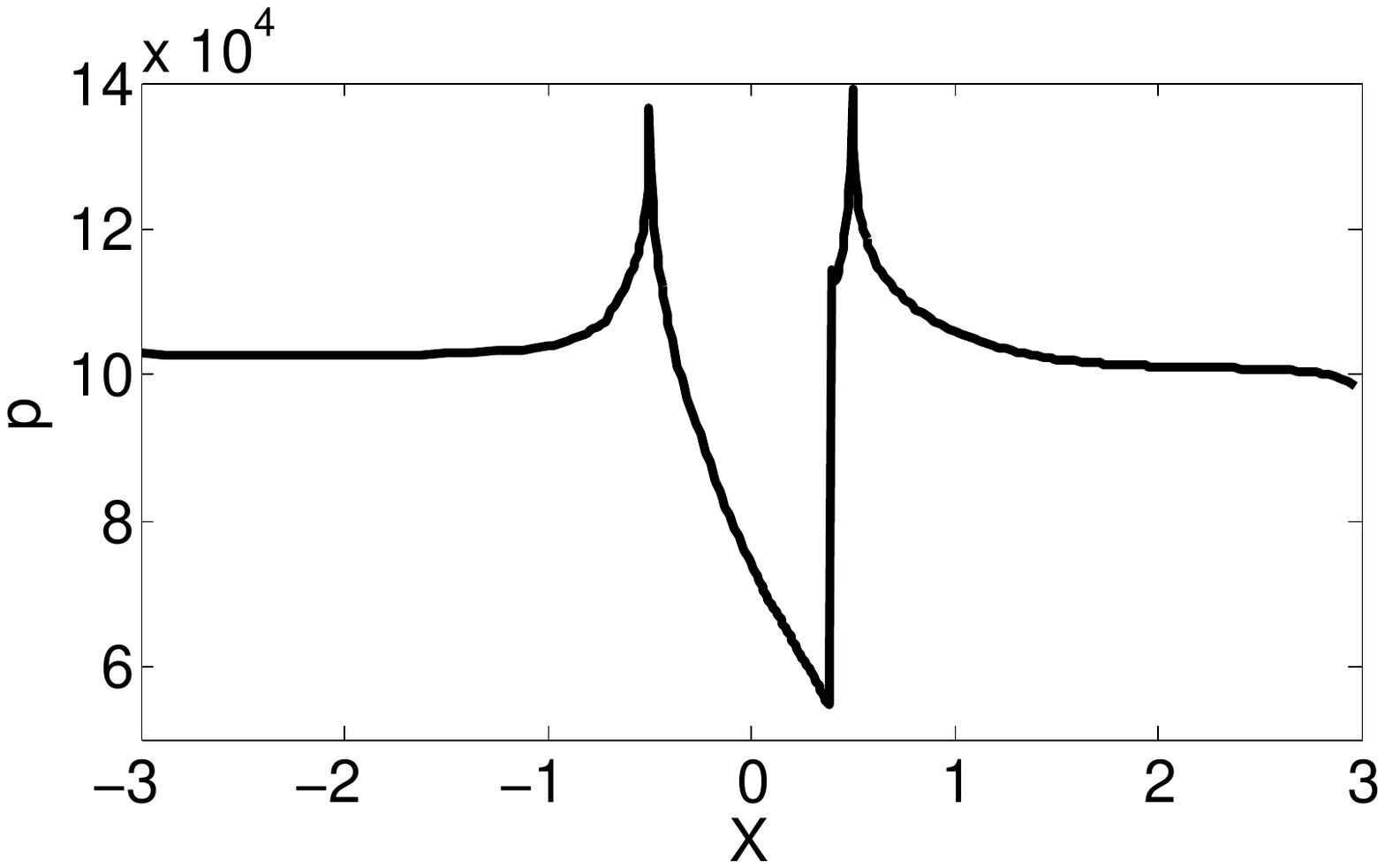}\vspace{-6cm}
\includegraphics[width=0.6\textwidth]{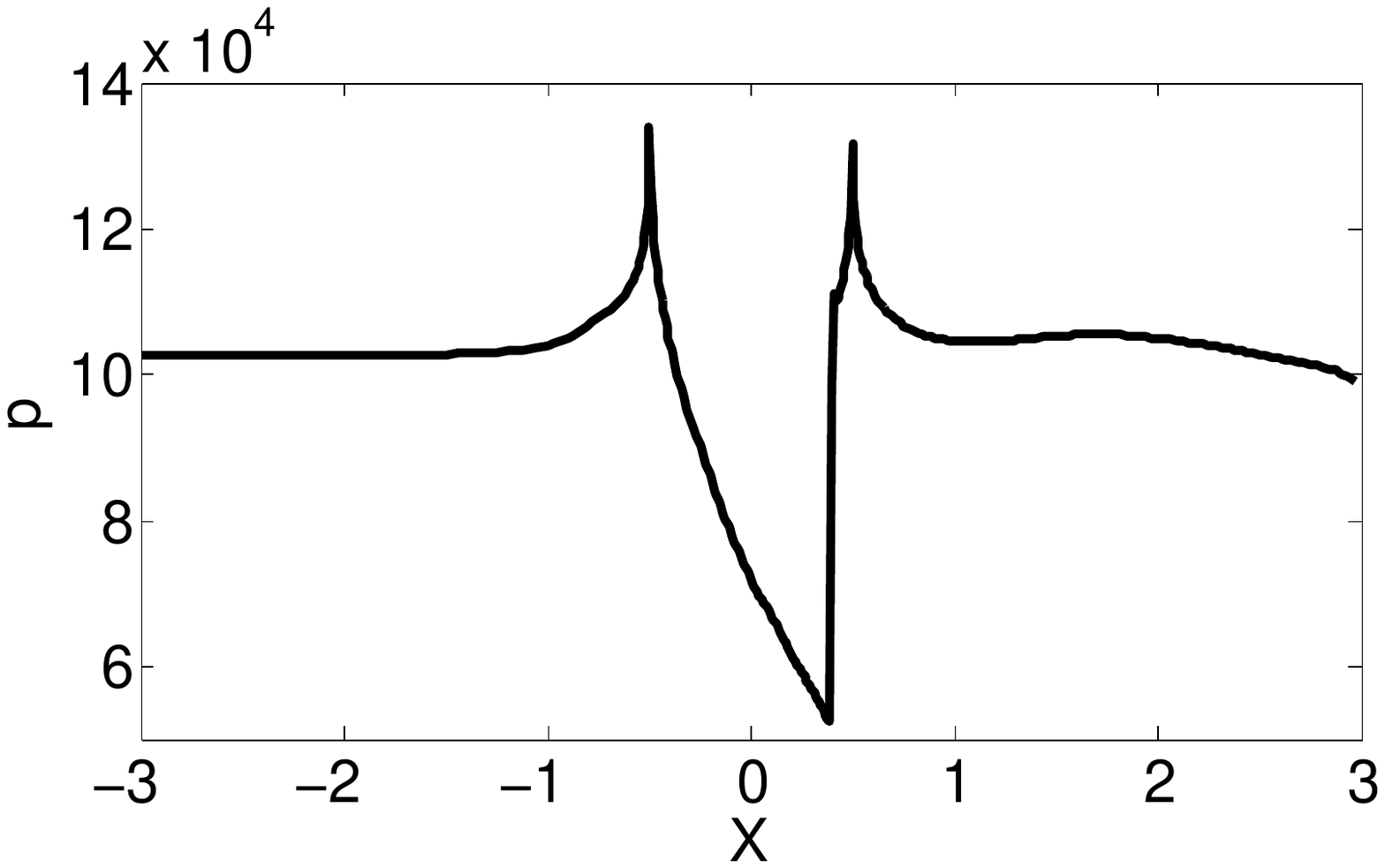}\vspace{-3cm}
\caption{Reference solution for 2D Euler equations: adaptive spatial grid (finest level $L=5$), uniform time steps ($\cfl=1$). Pressure $p$ at bottom boundary at times $t$=0.005002, 0.007125, 0.009990, 0.011975, from top to bottom.} \label{fig.vgl_euler2d_4}
\end{center}
\vspace*{-0.5cm}\end{figure}

In Figure \ref{fig.vgl_euler2d_4} we show a sequence of plots of the uniform, $\cfl=1$ computation on an adaptive spatial grid with finest level $L=5$. In Figure \ref{fig.vgl_euler2d_5} we compare the pressure distribution at the bottom boundary of the uniform solution and the time adaptive solution at several times. The two solutions on level $L=5$ match very well.

\subsection{Adaptive timesteps via ad hoc indicator}\label{sec.numadhoc}

\begin{figure}
  \begin{center}\vspace{-3cm}
  \includegraphics[width=0.8\textwidth]{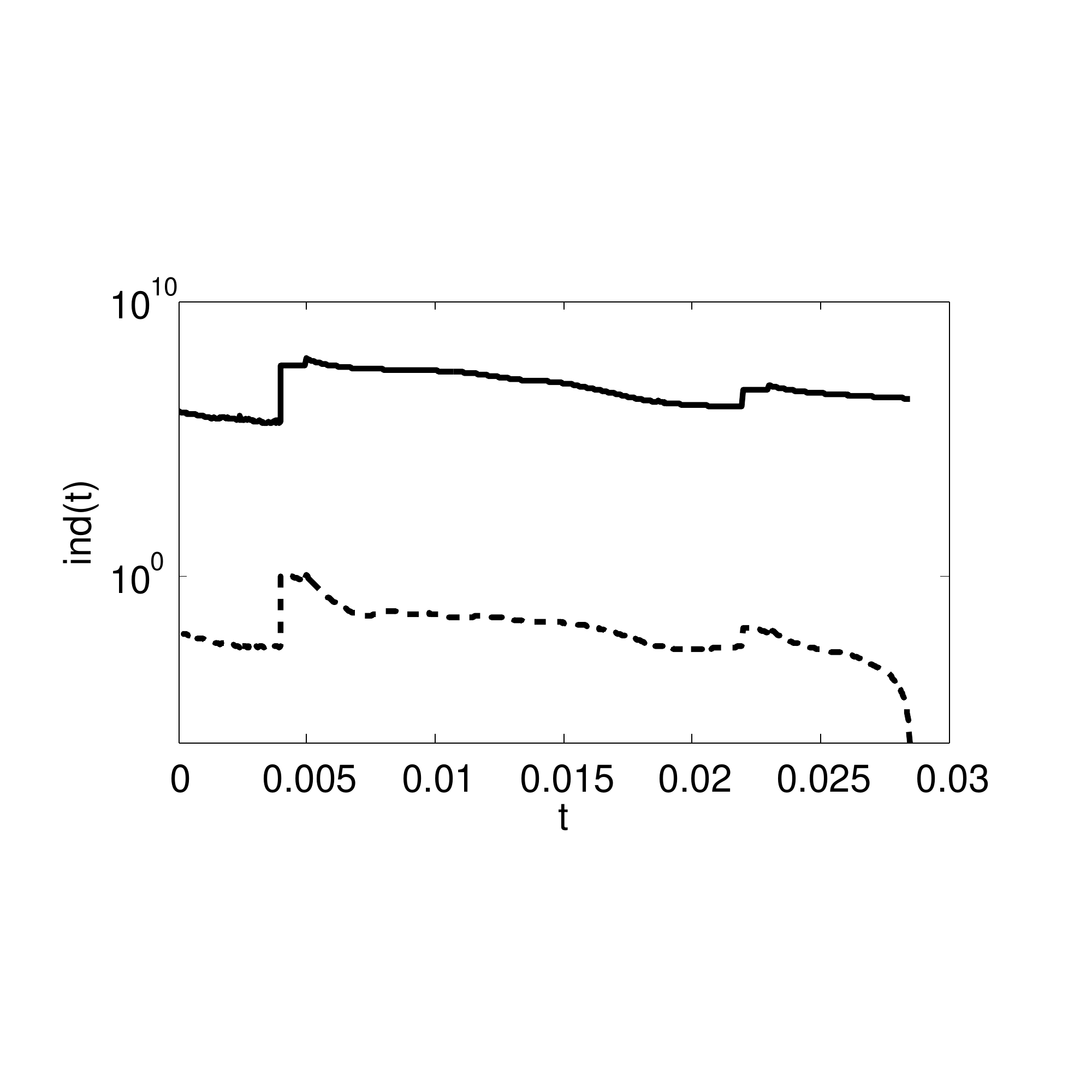}\vspace{-4.5cm}
  \includegraphics[width=0.8\textwidth]{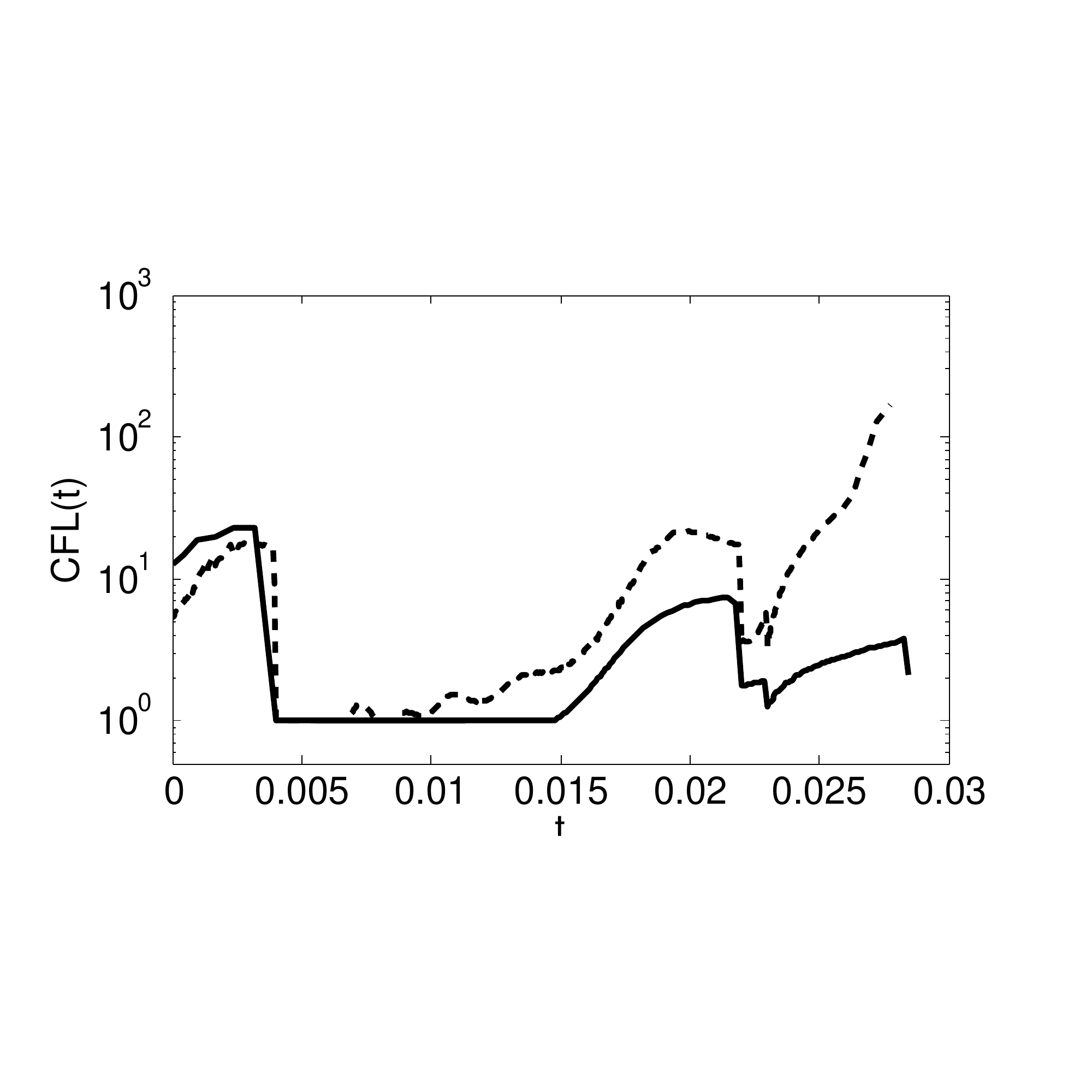}\vspace{-2.5cm}
  \caption{Comparison of error indicators (top) and timestep sequences (bottom)
    derived from error representation via dual problem (dashed line) and from
    variation indicator $ind$ (bold line).} \label{fig.oth_ind_tstps}
  \end{center}
\end{figure}

Here we replace the adjoint indicator by an ad hoc indicator, which estimates
the variation of the solution from one timestep to the following,
\mm{\label{eq.TVL}
  ind(\nn) = \sum_{i} |U_i^\nn -U_i^\no|_1 |\Oin|
  \approx \|U_h(\cdot,t^\nn)-U_h(\cdot,t^\no)\|_{L^1(\O)}.} 
Clearly this indicates whether the solution is stationary or not. Even though we
do not know any theoretically justified global decay rates of $ind$, and much
less of the error, as the timestep is refined, it is reasonable to assume that
the local variation in time decays linearly with the timestep. Similarly as for
the adjoint error control, we do a  first computation on a coarse grid with
finest level $L=2$, where we compute the indicator, see Figure
\ref{fig.oth_ind_tstps}. Then we redistribute the timesteps.

The indicator \eqref{eq.TVL} compares as follows to adjoint error indicator (see
Figure~\ref{fig.oth_ind_tstps}): Both indicators detect stationary and
nonstationary \replace{regions}{time intervalls}. The temporal distribution is very similar, but the 
variation indicator $ind$ leads to considerably more timesteps than the error
indicator from the dual approach (3640 vs. 2379). In particular, most timesteps
are smaller than in the case with adaptation via adjoint problems. Since the
timesteps are restricted from below, it leads to computations which are in
general more expensive but not more accurate. Results are not displayed. One
advantage may be that we do not need to compute a dual solution, which makes the
computation of the variation indicator less expensive. But this is only a small
advantage, since we compute the error indicators on a coarse mesh, which takes
$619 s$ on level $L=2$ for the adjoint approach. If we use the timesteps which
we compute from $ind(\nn)$ and do a time adaptive computation on a grid with
finest level $L=5$ then we will not get equally distributed indicators. This
holds also true if we do not apply the timestep restriction from below.

We have also implemented some variations of the discrete variation indicator,
which lead to similar results. Another approach was to choose the maximum jump
of the solution in one cell, both weighted and not weighted with the size of the
cell. This was an approximation to the $L^\infty$-norm. This indicator is not
very useful, since it turned out to be highly oscillating. Therefore we do not
present results for this indicator.

\subsection{Uniform timesteps}\label{sec.numuniform}

In many nonstationary computations where no a priori information is known, one
reasonable choice is to use uniform $\cfl$ numbers. Therefore we will compare
the computation with adaptive implicit timesteps with implicit computations
using uniform $\cfl$ numbers. In Section~\ref{sec.numadjoint} we have already
done a computation with uniform $\cfl$ number, $\cfl=1$, on a grid with $L=2$,
to get timestep sizes for an adaptive computation on a grid with $L=5$.  As a
reference solution we also computed with uniform $\cfl$ number, $\cfl =1$, a
solution of the problem on a grid with $L=5$. Now we compare these computations
with computations using higher uniform $\cfl$ numbers.

\begin{figure}[ht]
\begin{center}\vspace{-2.7cm}
\includegraphics[width=0.6\textwidth]{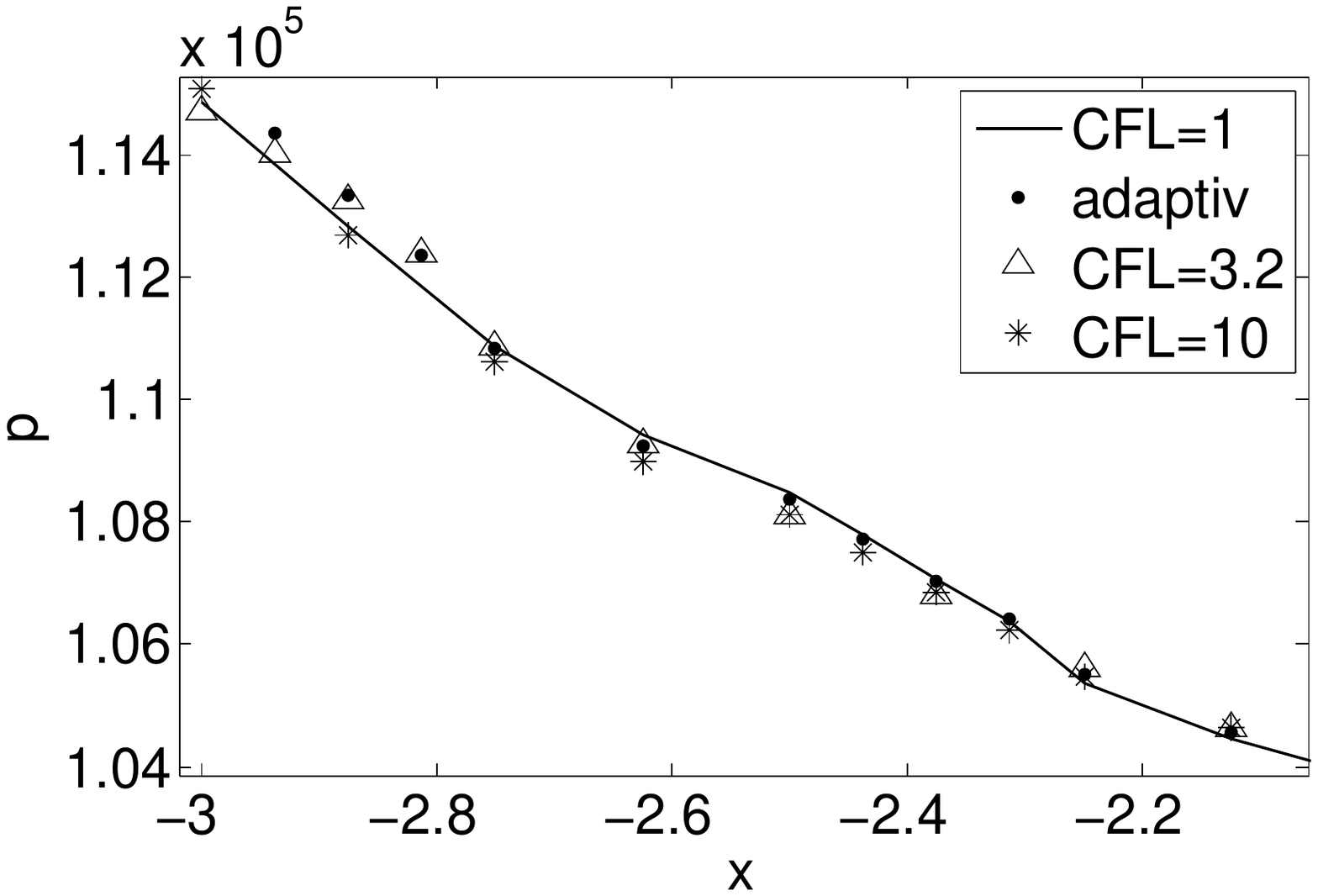}\vspace{-6cm}
\includegraphics[width=0.6\textwidth]{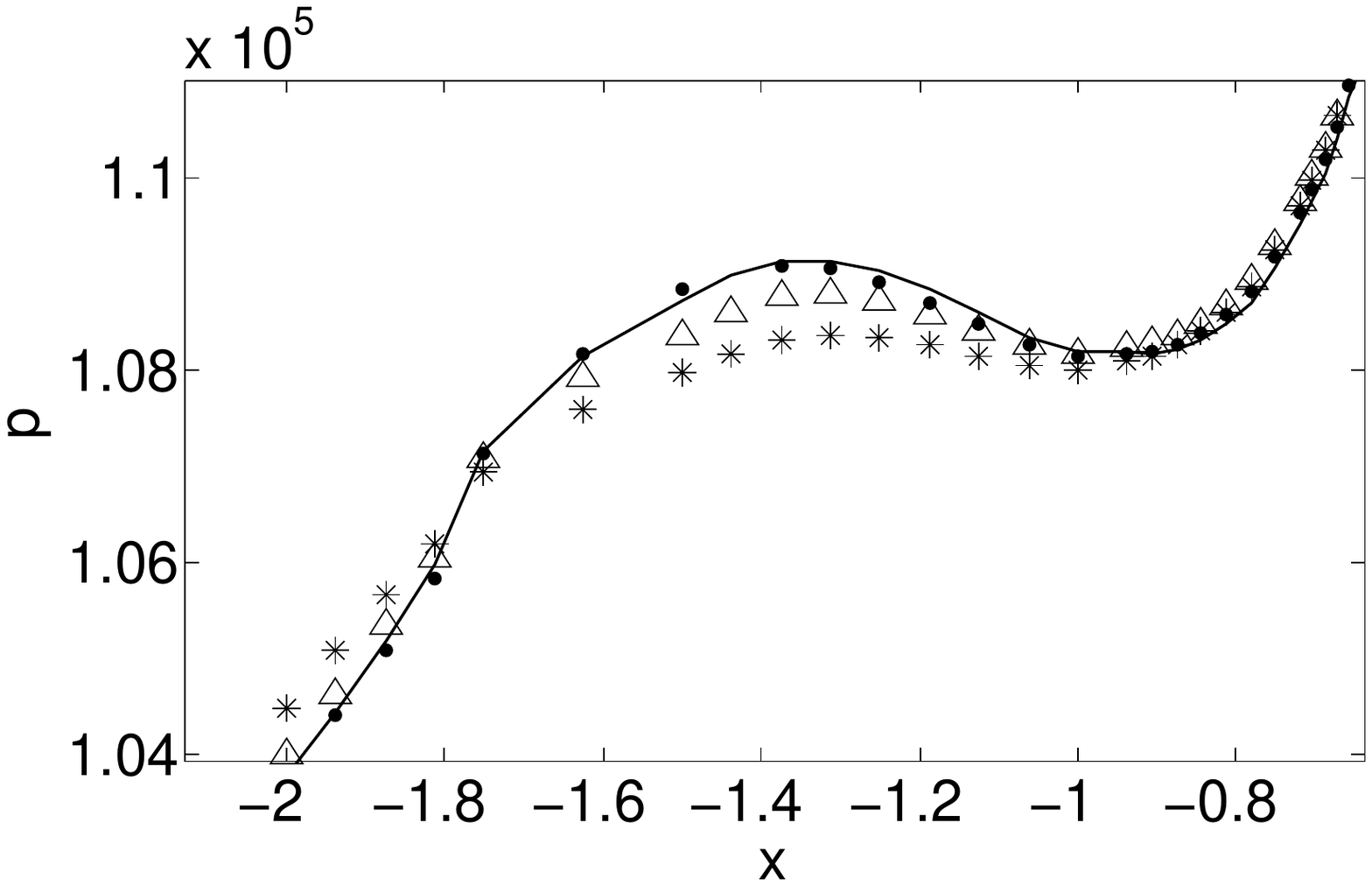}\vspace{-6cm}
\includegraphics[width=0.6\textwidth]{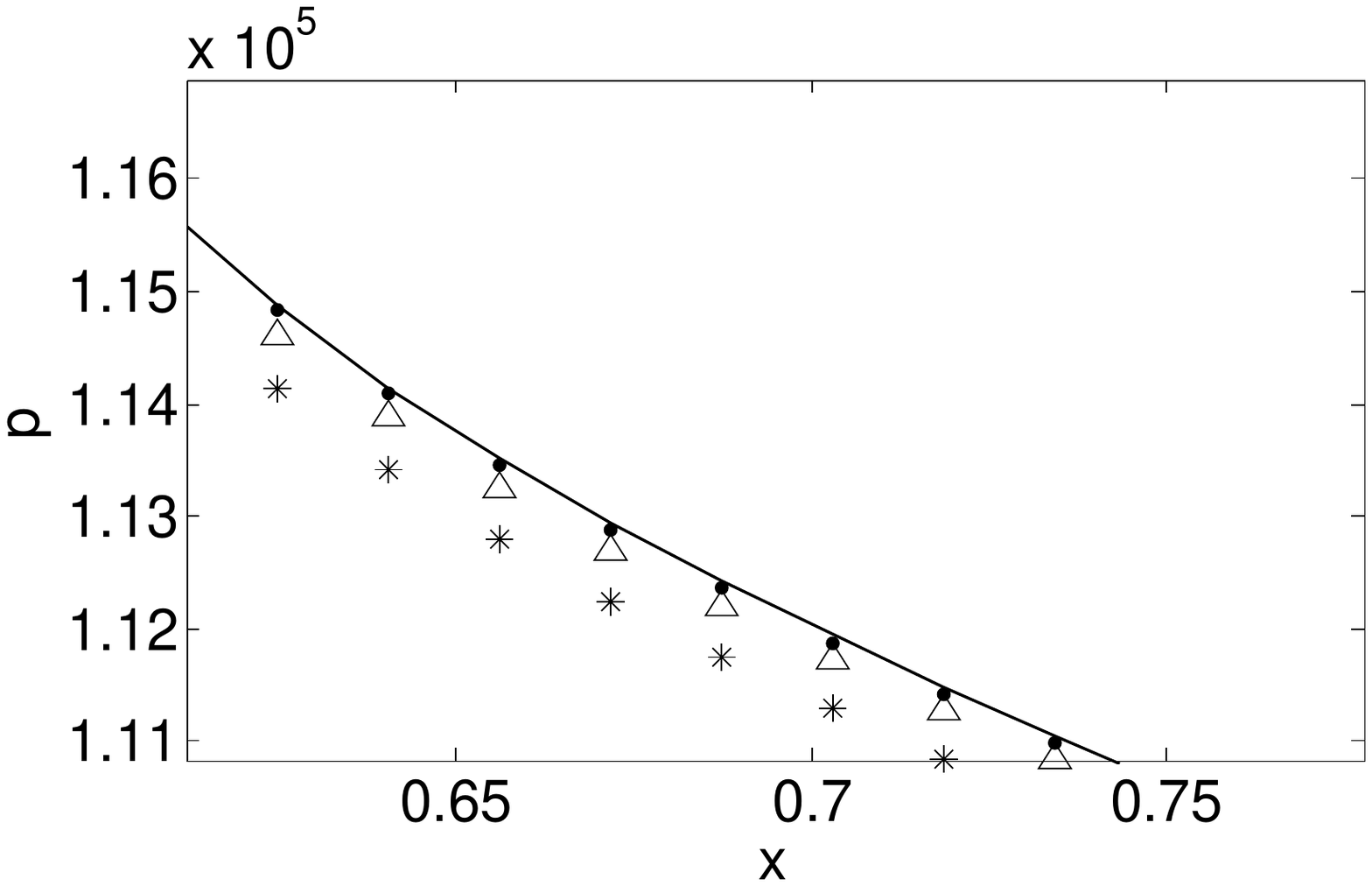}\vspace{-6cm}
\includegraphics[width=0.6\textwidth]{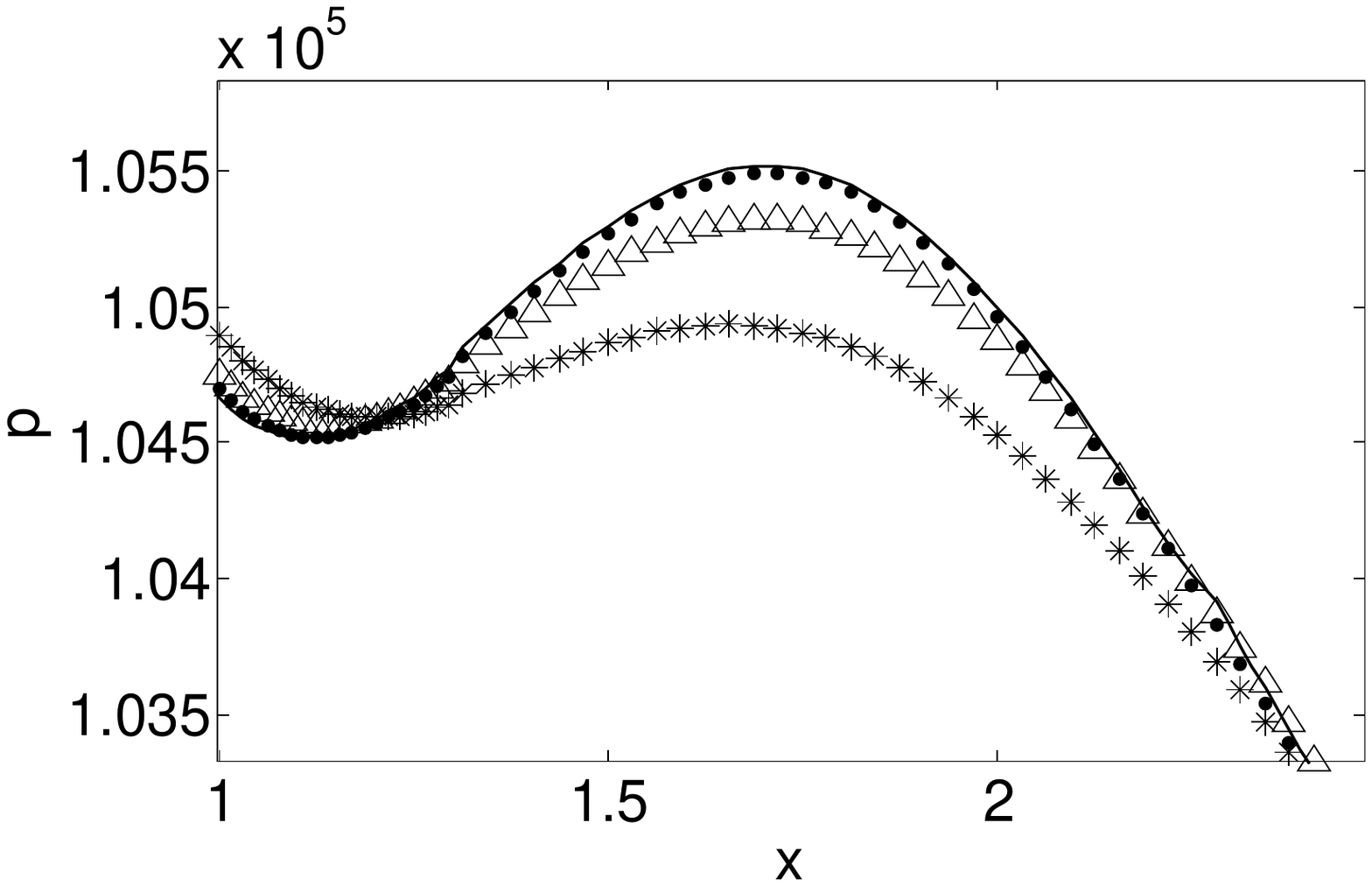}\vspace{-3.3cm}
\caption{Zoom tracing the perturbations in Figure~\ref{fig.vgl_euler2d_4}
at times $t$=0.005002, 0.007125, 0.009990, 0.011975, from top to bottom. Comparison of uniform timesteps with $\cfl = 1,\;3.2,\;10$ and time-adaptive strategy.} \label{fig.vgl_euler2d_5}
\end{center}
\end{figure}

%\begin{figure}%
%\centering
%\subfigure[$t$=0.005002]{\includegraphics[width=0.595\textwidth,trim = 7mm 84mm 8mm 83mm,clip]{\FigurePath/diss/5all_1404_new.pdf}}\\
%\subfigure[$t$=0.007125]{\includegraphics[width=0.595\textwidth,trim = 7mm 87mm 8mm 86mm,clip]{\FigurePath/diss/5all_2000_new.pdf}}\\
%\subfigure[$t$=0.009990]{\includegraphics[width=0.595\textwidth,trim = 7mm 87mm 8mm 86mm,clip]{\FigurePath/diss/5all_2806_new.pdf}}\\
%\subfigure[$t$=0.011975]{\includegraphics[width=0.595\textwidth,trim = 7mm 87mm 8mm 86mm,clip]{\FigurePath/diss/5all_3362_new.pdf}}
% \caption{Zoom tracing the perturbations in Figure~\ref{fig.vgl_euler2d_4}
%at several times. Comparison of uniform timesteps with $\cfl = 1,\;3.2,\;10$ and time-adaptive strategy.}
%\label{fig.vgl_euler2d_5}
%\end{figure}

First we choose a uniform $\cfl$ number of approximately 3.2, which corresponds
to 2500 timesteps. This equals roughly the number of timesteps in the
adaptive method, and hence it should give a fair comparison. The uniform
computation takes about $11509 s$, more than the $9142 s$ of the adaptive
computation (see Table~\ref{table.comp2}).  In the uniform computation most of
the timesteps are more expensive, since they need more Newton steps, and more
steps for solving the linear problems. This shows that the understanding of the
dynamics of the solution pays directly in the nonlinear and linear solvers.
Moreover, it can be seen from Figure~\ref{fig.vgl_euler2d_5} that the quality of
the solution is considerably worse than for the adaptive computation.

Another computation with $\cfl$ number 10 takes only $3290 s$. However, as can
be seen from Figure~\ref{fig.vgl_euler2d_5} the solution is badly approximated:
In the beginning of the computation the solutions of the different methods match
very well, which means that the inflow at the boundary is well-resolved. As time
goes on, the solutions differ more and more. After the perturbation has passed
the bump, the perturbations differ strongly. Only the time-adaptive method
approximates the reference solution ($\cfl=1$) closely.

\subsection{Newton iterations and linear iterations}
\label{sec.Newtonliniterations}

\begin{figure}%
\centering
\subfigure[adaptive $\cfl$]{\includegraphics[width=0.6\textwidth,trim = 7mm 100mm 8mm 80mm,clip]{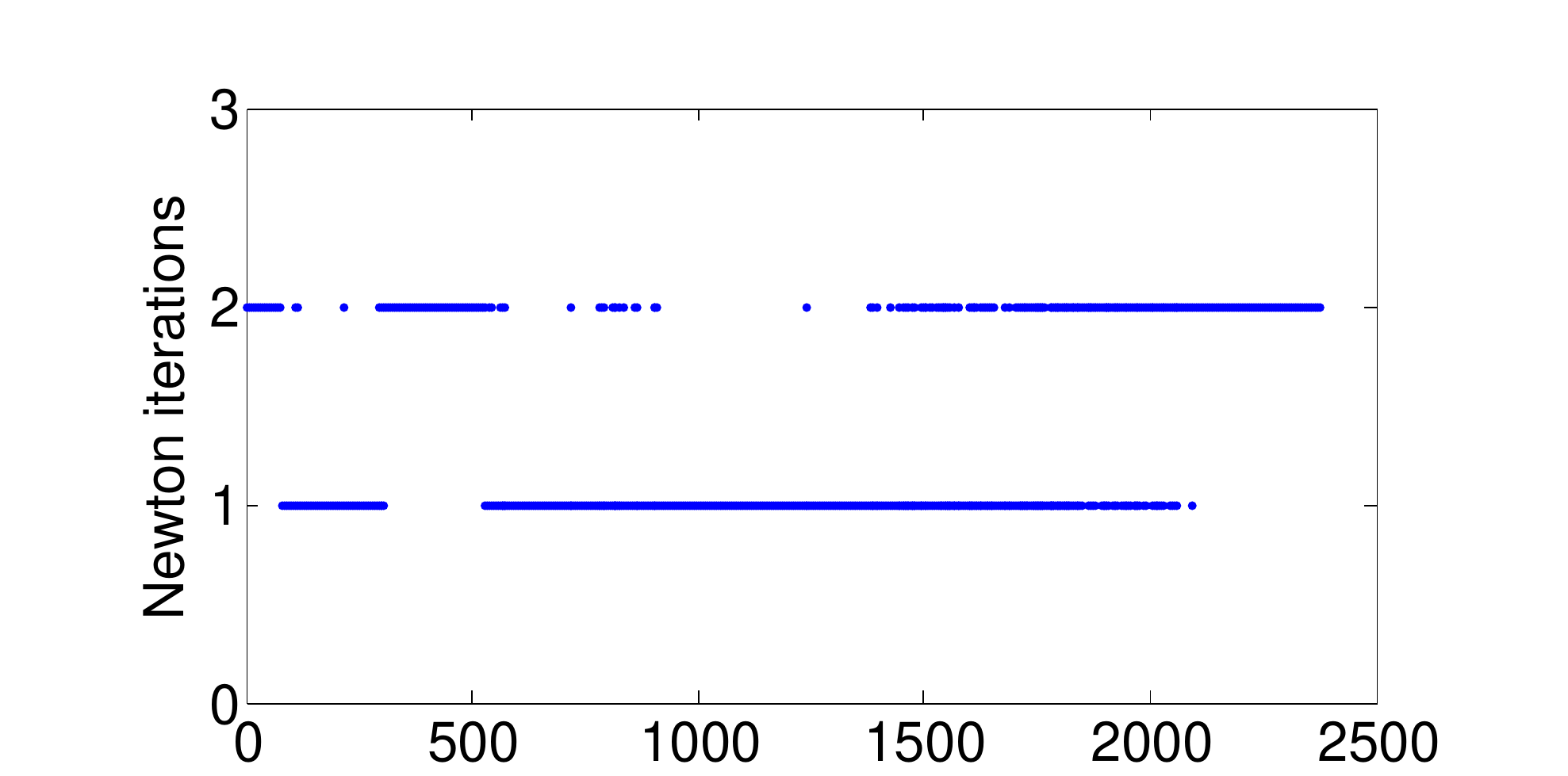}}\\
\subfigure[$\cfl=1$]{\includegraphics[width=0.6\textwidth,trim = 7mm 100mm 8mm 80mm,clip]{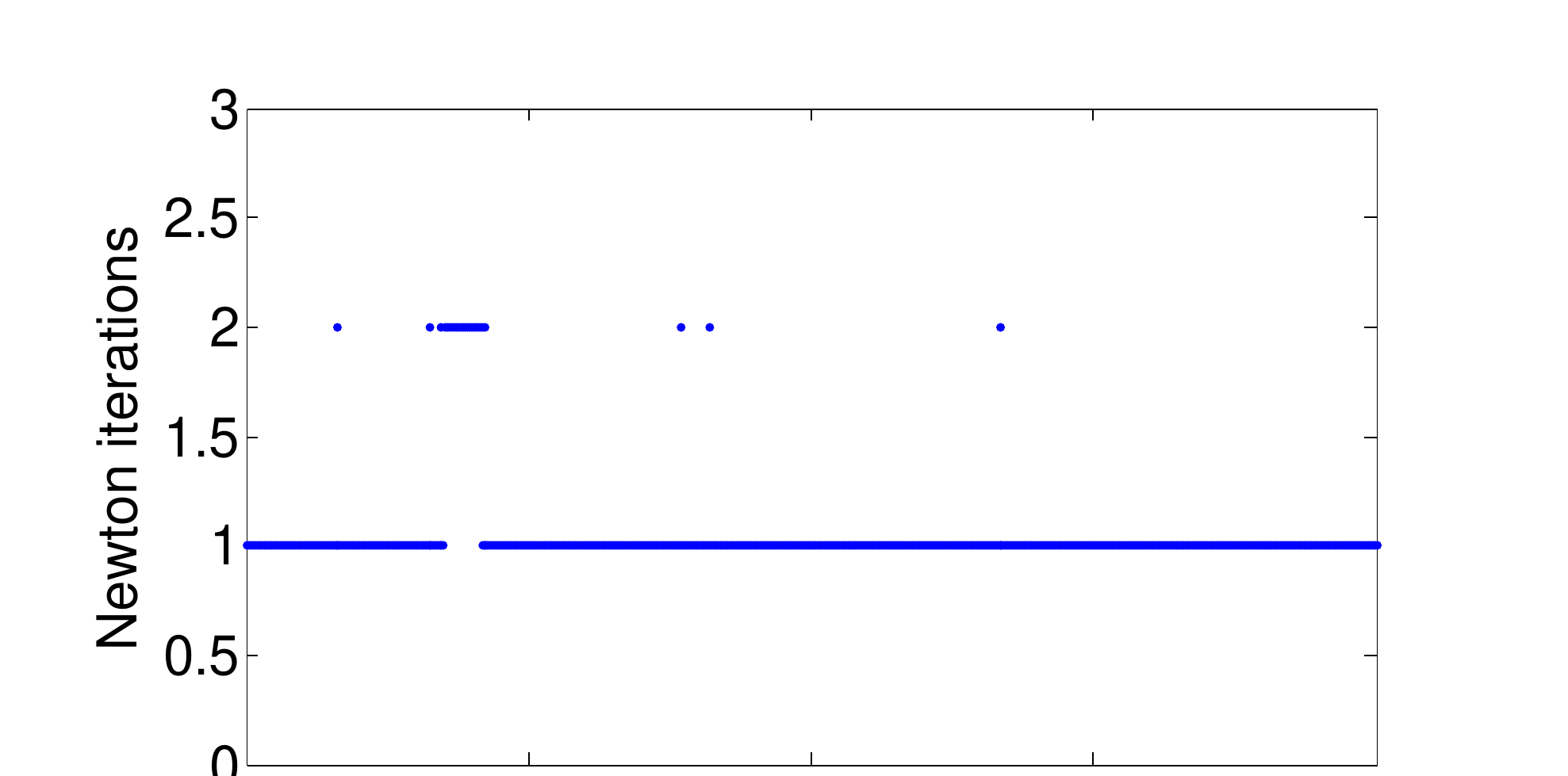}}\\
\subfigure[$\cfl\approx$3.2]{\includegraphics[width=0.6\textwidth,trim = 7mm 100mm 8mm 80mm,clip]{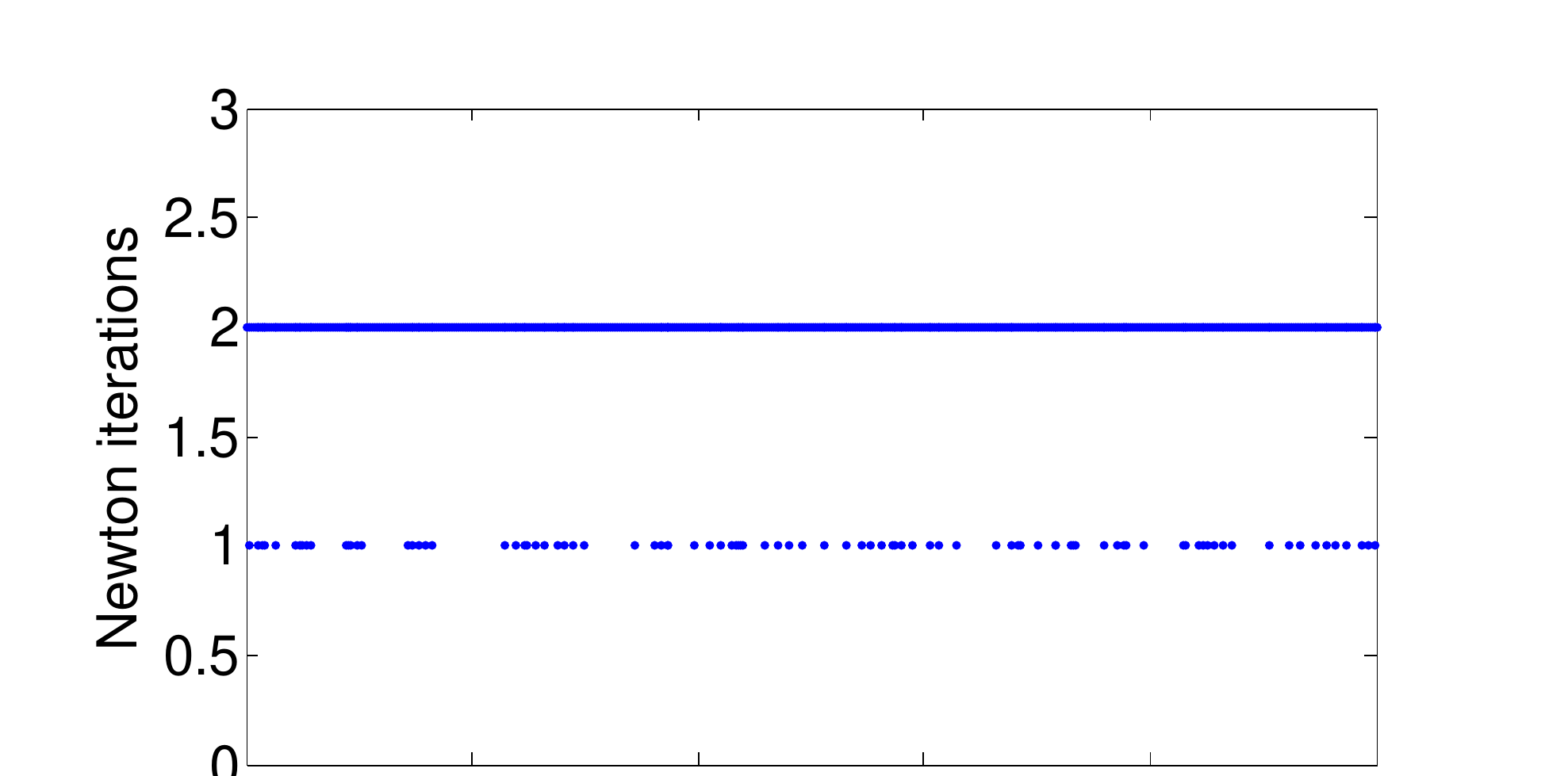}}\\
\subfigure[$\cfl=10$]{\includegraphics[width=0.6\textwidth,trim = 7mm 100mm 8mm 80mm,clip]{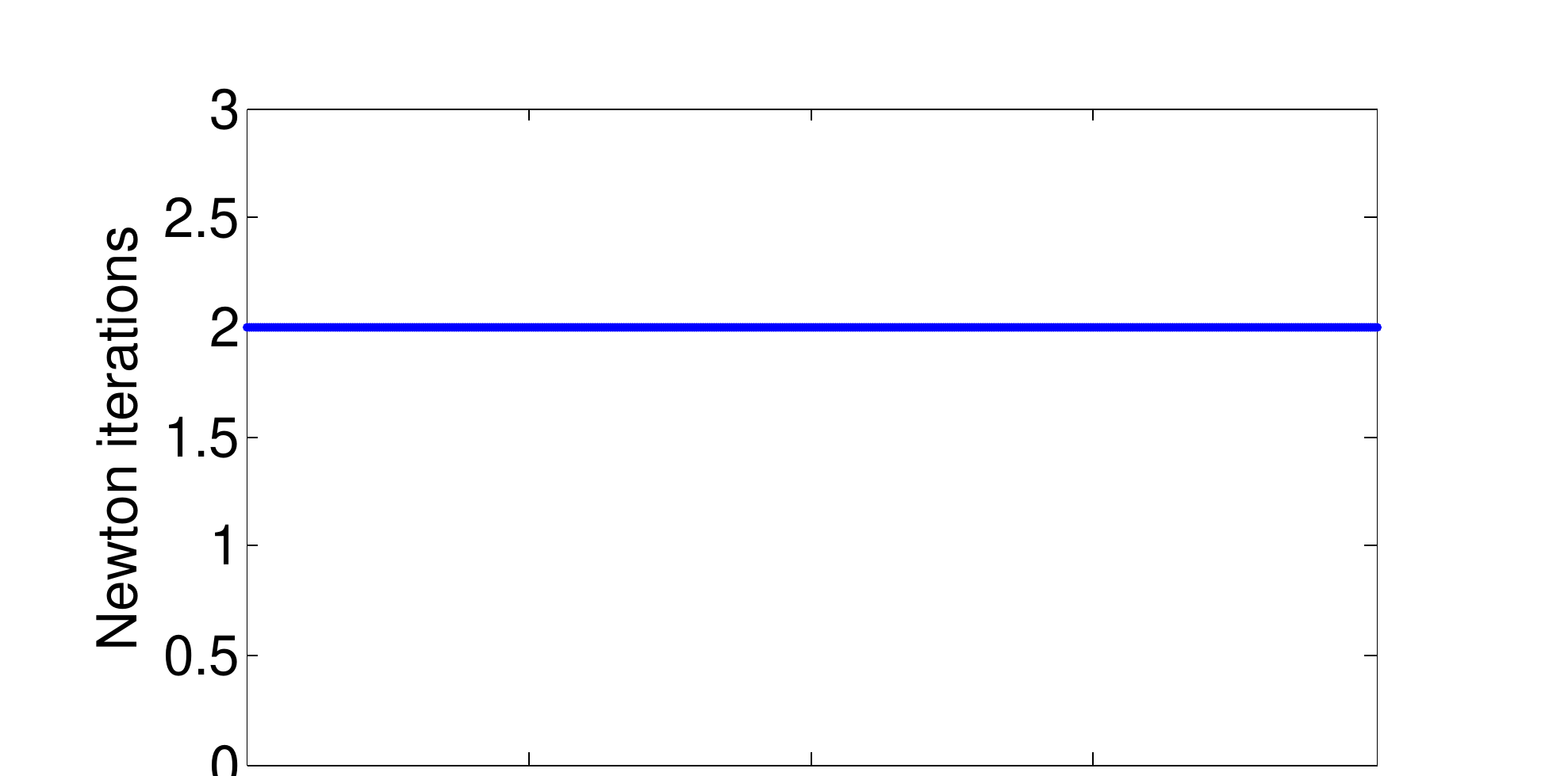}}
 \caption{Number of Newton iterations for fully implicit time-adaptive
  computation and computations with uniform $\cfl$ number}
\label{fig.Newtoniterations}
\end{figure}

In Figure \ref{fig.Newtoniterations} we show the number of Newton iterations in
each timestep for the computations in Section~\ref{sec.numadjoint},
\ref{sec.numadhoc} and \ref{sec.numuniform}. The number of Newton iterations
depend on the $\cfl$ number and is larger, where the solution is nonstationary
and smaller, where the solution is stationary. The time-adaptive computation
gives the smallest total number of Newton iterations, see Table
\ref{table.comp2}, since most timesteps are solved with one Newton iteration.

%\begin{figure}[hbtp]
%\begin{center} {\vspace*{-25mm}
%  {\includegraphics[width=0.6\textwidth]{\FigurePath/diss/5a.pdf}}\vspace*{-60mm}
%  {\includegraphics[width=0.6\textwidth]{\FigurePath/diss/5u.pdf}}\vspace*{-60mm}
%  {\includegraphics[width=0.6\textwidth]{\FigurePath/diss/5u3.pdf}}\vspace*{-60mm}
%  {\includegraphics[width=0.6\textwidth]{\FigurePath/diss/5u10.pdf}}}\vspace*{-35mm}
%  \end{center}
%  \caption{Number of Newton iterations for fully implicit time-adaptive
%  computation and computations with uniform $\cfl$ number, from top to 
%  bottom: adaptive $\cfl$, $\cfl=1$, $\cfl\approx$3.2, $\cfl=10$.
%  \label{fig.Newtoniterations}
%  }
%\end{figure} 

\begin{figure}%
\centering
\subfigure[adaptive $\cfl$]{\includegraphics[width=0.595\textwidth,trim = 7mm 86mm 8mm 86mm,clip]{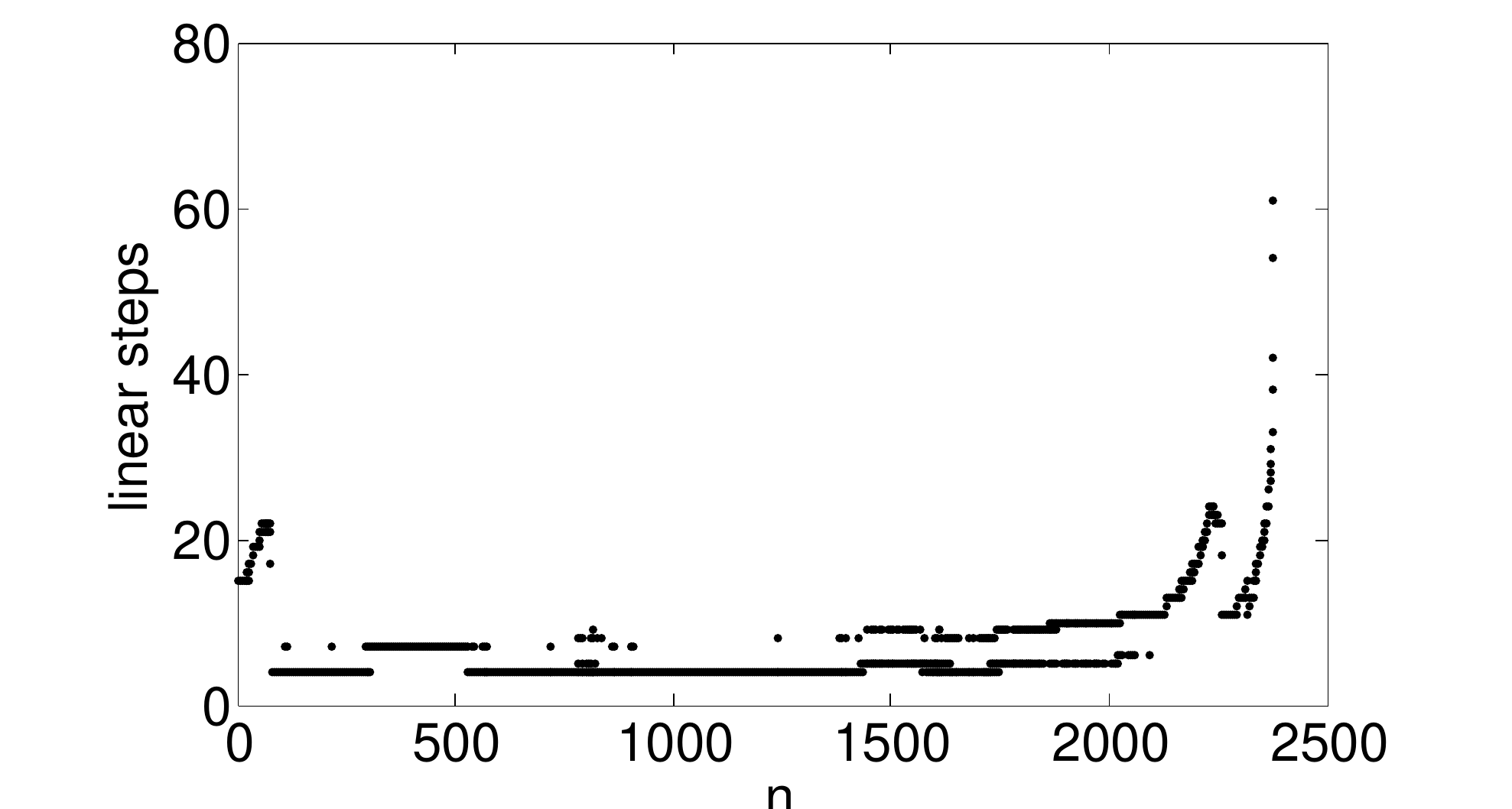}}\\
\subfigure[$\cfl=1$]{\includegraphics[width=0.595\textwidth,trim = 7mm 86mm 8mm 86mm,clip]{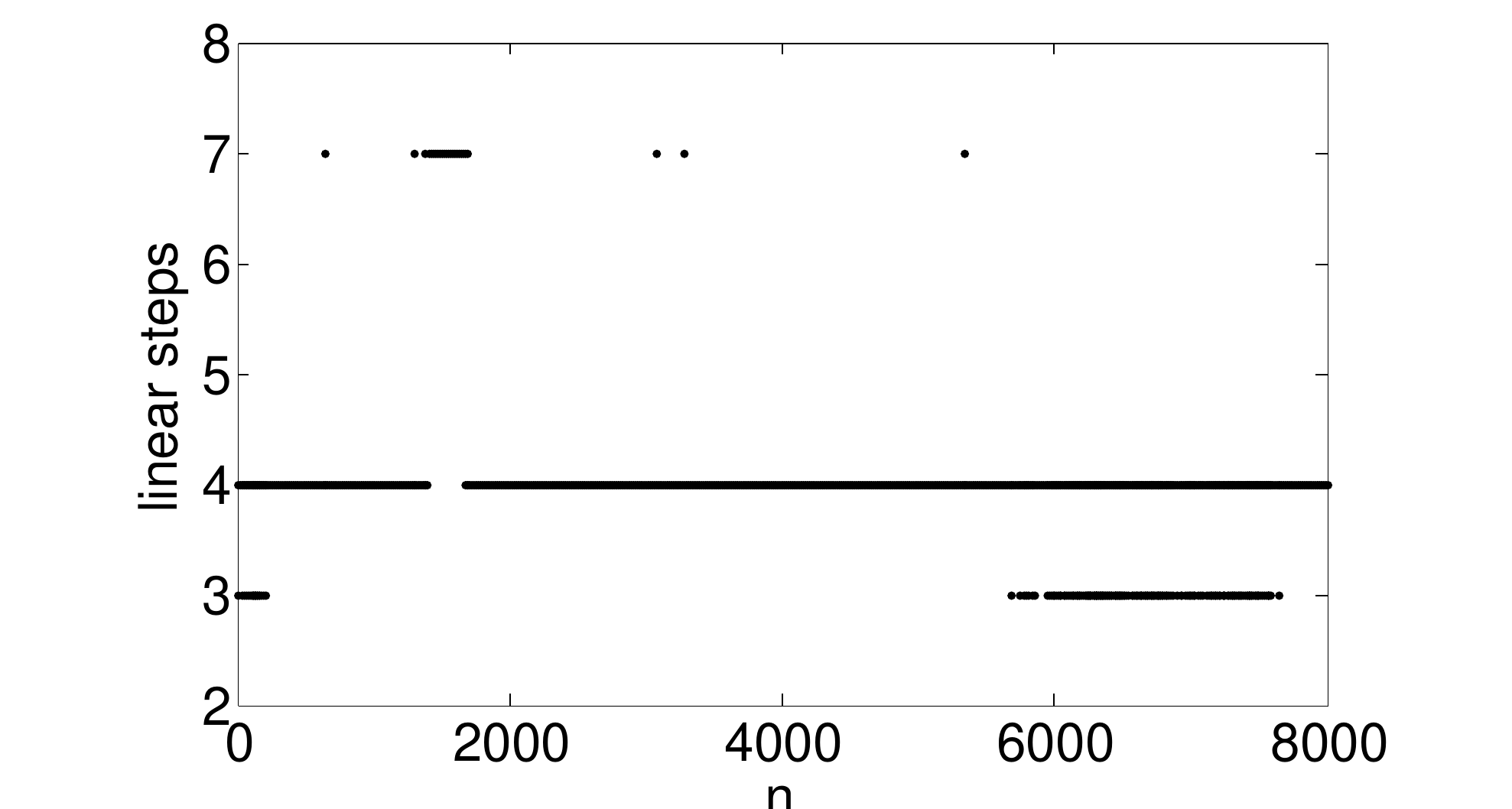}}\\
\subfigure[$\cfl\approx$3.2]{\includegraphics[width=0.595\textwidth,trim = 7mm 86mm 8mm 86mm,clip]{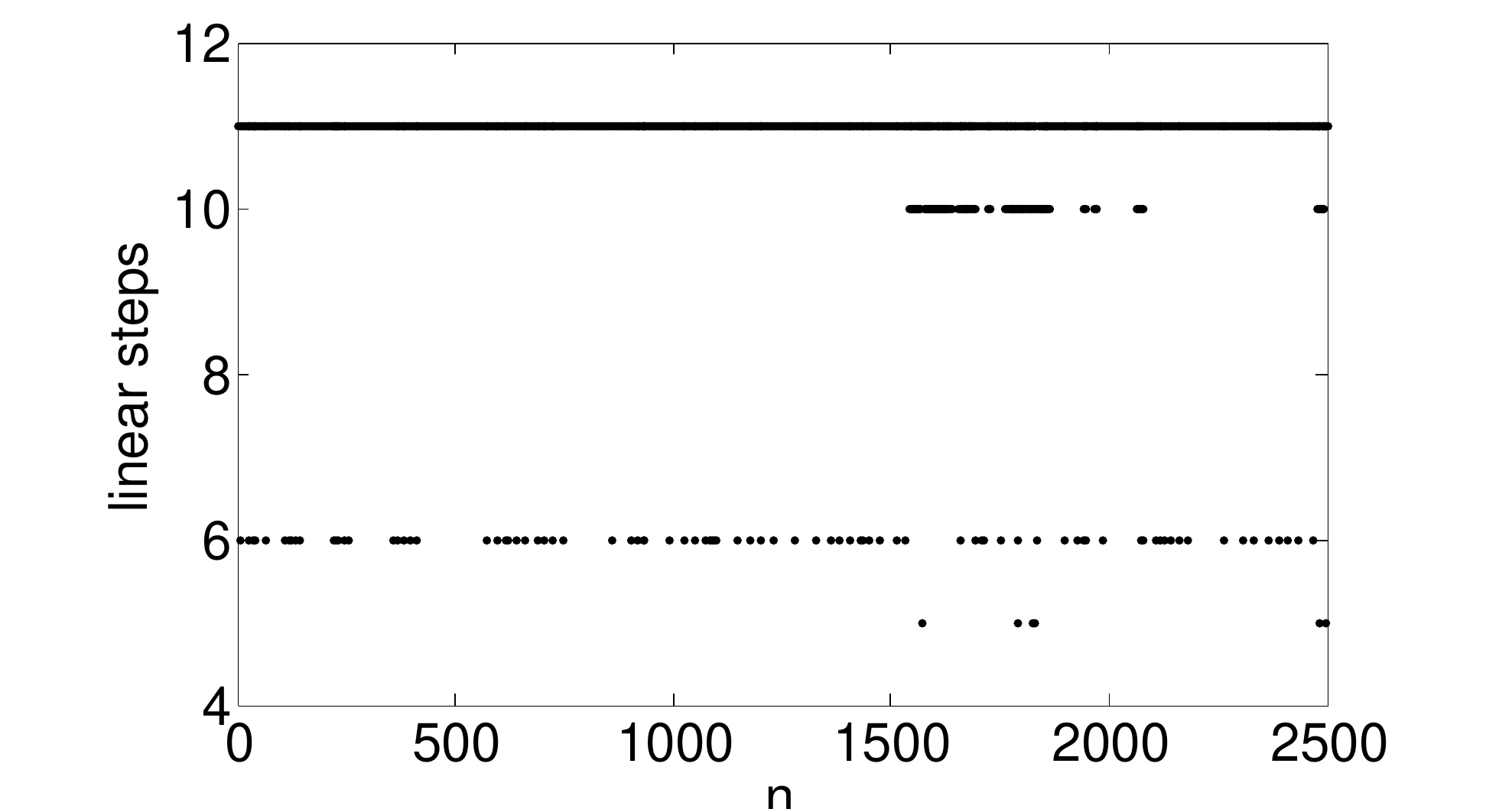}}\\
\subfigure[$\cfl=10$]{\includegraphics[width=0.595\textwidth,trim = 7mm 86mm 8mm 86mm,clip]{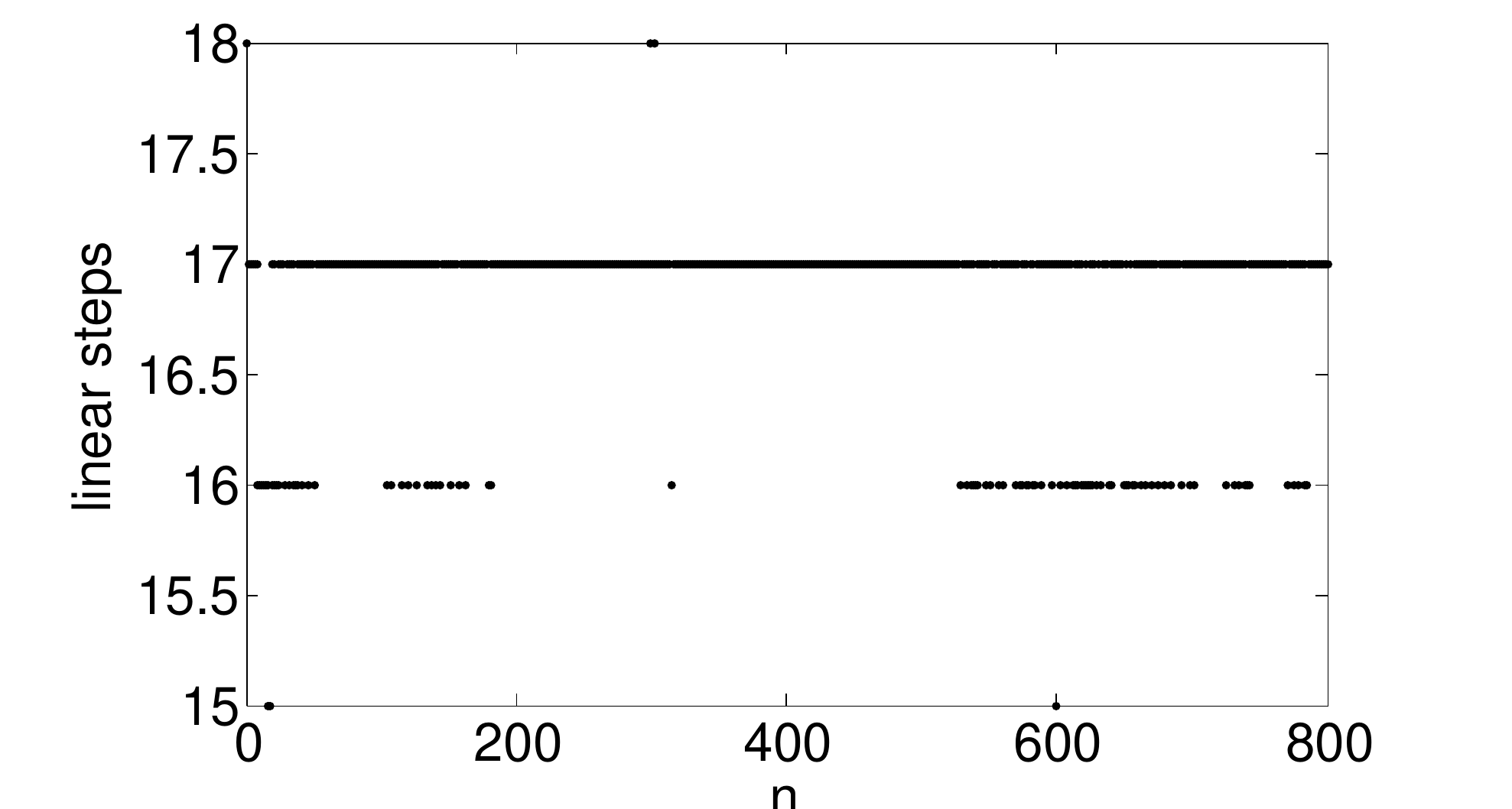}}
 \caption{Number of linear iterations for fully implicit time-adaptive
 computation and computations with uniform $\cfl$ number}
\label{fig.Lineariterations}
\end{figure}

Figure \ref{fig.Lineariterations} shows the total number of linear iterations in
each timestep for the same computations, i.e., the sum of the linear iterations
in each timestep for all Newton iterations. The number of linear iterations also
depends on the $\cfl$ number, for larger $\cfl$ number, we observe more linear
iterations.

%\begin{figure}[hbtp]
%\begin{center} {\vspace*{-25mm}
%  {\includegraphics[width=0.6\textwidth]{\FigurePath/diss/5alin.pdf}}\vspace*{-60mm}
%  {\includegraphics[width=0.6\textwidth]{\FigurePath/diss/5ulin.pdf}}\vspace*{-60mm}
%  {\includegraphics[width=0.6\textwidth]{\FigurePath/diss/5u3lin.pdf}}\vspace*{-60mm}
%  {\includegraphics[width=0.6\textwidth]{\FigurePath/diss/5u10lin.pdf}}}\vspace*{-35mm}
%  \end{center}
%  \caption{Number of linear iterations for fully implicit time-adaptive
%  computation and computations with uniform $\cfl$ number, from top to
%  bottom: adaptive $\cfl$, $\cfl=1$, $\cfl\approx$3.2, $\cfl=10$.
%  \label{fig.Lineariterations}}
%  \vspace*{-5mm}
%\end{figure} 

Table~\ref{table.comp2} gives an overview of the CPU time and the number of
Newton iterations and linear iterations for the computations.  The costs for the
computation of the indicators, both the dual and the ad hoc indicator, on level
$L=2$ are very low compared to the costs of computations on level $L=5$.  An
adaptive computation including the computation of the indicator, i.e., $329 s +
619 s + 9142 s$, is cheaper than the computation using uniform $\cfl$ number,
e.g., $\cfl=1$, that needs 21070$s$. Even the computation with $\cfl=3.2$ is more
expensive than the time-adaptive computation, but leads to worse results, see
Figure \ref{fig.vgl_euler2d_5}.

Table~\ref{table.comp2} shows that the CPU time is roughly proportional to the
number of Newton iterations and not to the number of timesteps or linear solver
steps. The CPU time is about $2.5 s$ per Newton iteration. This means that we
have to minimize the number of Newton iterations in total to accelerate the
computation. This is done very efficiently by the time-adaptive approach. For a
large range of $\cfl$ numbers from 1 to more than 100, it needs only one or two
Newton iterations per timestep, without sacrificing the accuracy.

In this example, the stationary time \replace{regions}{intervalls} are not very large compared to the
overall computation. If the stationary \replace{regions}{time intervalls} were larger, the advantage of the
time adaptive scheme would be even more significant.

\iffalse
\begin{center}
\begin{table}[ht] 
\begin{tabular}{l|c|c}
&adjoint $[s]$ & ad hoc indicator $[s]$ \\ \hline
forward $\cfl=1$& 329&329  \\  
dual problem and error indicator& 619 &- \\ 
error indicator &- & $t\ll 619 ^\ast$\\ \hline
total& 948& $329<t \ll 948 ^\ast$
\end{tabular}\\[1ex] 
\caption{Computational costs for computing a timestep sequence by the adjoint method and an ad hoc indicator on $L=2$, $^\ast$currently the computation is not efficiently implemented.}\label{table.comp1} 
\end{table} 
\end{center}
\fi

\begin{center}
\begin{table}[ht]  
\begin{tabular}{l|c|c|c|c}
&CPU [s]& timesteps & Newton steps & linear steps \\ \hline
adaptive timesteps & 9142& 2379 &3303&16875\\
uniform $\cfl=1$ & 21070&8000&8282&32630\\
uniform $\cfl=3.2$ & 11509&2500 &4904&26895\\
uniform $\cfl=10$ & 3290&800&1600& 13500\\
\end{tabular}\\[1ex] 
\caption{Performance for computations on $L=5$ using different fully implicit
timestepping strategies.}\label{table.comp2}
\end{table} 
\end{center}

\section{Explicit-implicit computational results}
\label{section.numexp_explicit_implicit}

\begin{figure}
  \begin{center}
  \includegraphics[width=0.75\textwidth]{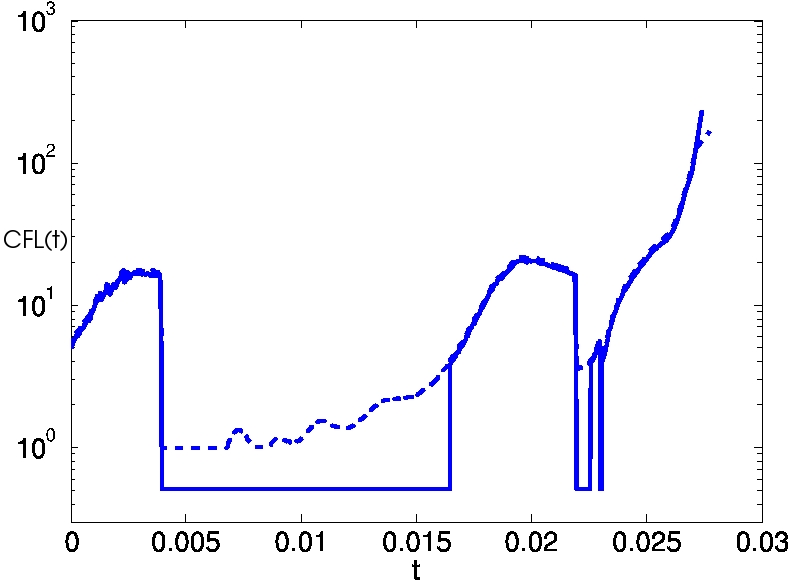}%\vspace*{-10mm}
  \caption{2D Euler equations, comparison of timestep sequence derived from error
  representation for implicit computation (dashed line) and for mixed
  explicit-implicit computation  (bold line).}
  \label{fig.vgl_ex-im_CFL}
  \end{center}
\end{figure}

Now we modify the fully implicit timestepping strategy and introduce a mixed {\em
explicit-implicit} approach. The reason is that implicit
timesteps with $\cfl<5$ are not efficient, since we have to solve
a nonlinear system of equations at each timestep. Thus, for $\cfl<5$,
the new explicit-implicit strategy switches to the cheaper and less dissipative
explicit method with $\cfl=0.5$. The timestep sequence is shown in 
Figure~\ref{fig.vgl_ex-im_CFL}. Of course, we could choose variants of
the thresholds $\cfl=0.5$ and 5.

\begin{figure}
  \begin{center}%\vspace*{-10mm}
  \includegraphics[width=0.85\textwidth]{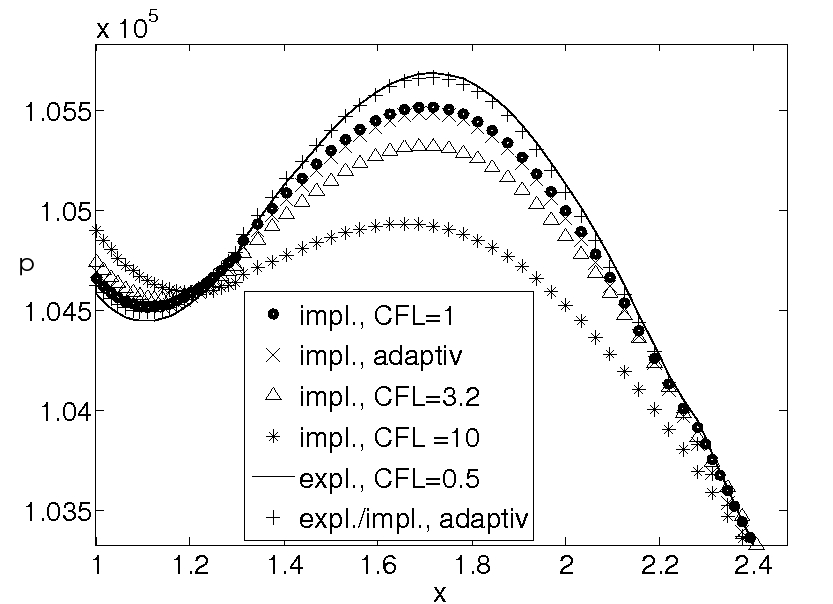}%\vspace*{-40mm}
  \caption{2D Euler equations, comparison of solutions in time on adaptive grid
  with finest level $L=5$, pressure $p$ at the bottom boundary, zoom of the
  perturbation at time $t$=0.011975}
  \label{fig.vgl_ex-im}
  \end{center}
\end{figure}

As we can see in Table~\ref{table.comp3}, the new strategy requires 5802
timesteps, where $95\%$ are explicit. The CPU time of $7730 s$ easily beats
the fully explicit solver ($18702 s$), and is also superior to the fully
implicit adaptive scheme ($ 9142 s$, see Table \ref{table.comp2}). 
The computational results are presented in Figure \ref{fig.vgl_ex-im}. The
results of the combined explicit-implicit strategy are very close to the results
of the fully explicit method, and far superior to all fully implicit methods. 
Note that the explicit scheme serves as reference solution, since it is
well-known that it gives the most accurate solution for an nonstationary problem.

\begin{center}
\begin{table}[ht]  
\begin{tabular}{l|c|c|c}
&CPU [s]& timesteps & timesteps \\
&       & (total)   & (implicit)\\ \hline
explicit $\cfl=0.5$ & 18702&16000&-\\
adaptive expl.-impl. & 7730&5802&259\\
\end{tabular}\\[1ex] 
\caption{Performance for computations on $L=5$ using fully explicit timesteps
and the time-adaptive explicit-implicit strategy.}
\label{table.comp3}
\end{table} 
\end{center}

\section{Conclusion}
\label{section:conclusion}

In this work, explicit and implicit finite volume solvers  on adaptively refined
quadtree meshes have been coupled with adjoint techniques to control the
timestep sizes for the solution of weakly nonstationary compressible inviscid
flow problems.

For the 2D Euler equations we have presented a test case for which the
time-adaptive method  does reach its goals: it separates stationary \replace{regions}{time intervalls} and
perturbations cleanly and chooses just the right timestep for each of them. The
adaptive method leads to considerable savings in CPU time and memory while
reproducing the reference solution almost perfectly.

We have compared the adjoint error representation with several variation-based
indicators. Our prime choice is the adjoint approach, since it has the best
theoretical justification  and needs only half the number of timesteps.

In Theorem~\ref{theorem:error_representation_h} and
Corollary~\ref{cor:error_representation_hk} we state a complete error
representation for nonlinear initial-boundary-value problems with characteristic
boundary conditions for hyperbolic systems of conservation laws, which includes
boundary and linearization errors. Besides building upon well-established
adjoint techniques, we also add a new ingredient which simplifies the
computation of the dual problem \cite{SteinerNoelle}. We show that it is
sufficient to compute the  spatial gradient of the dual solution, $w=\nabla
\vp$, instead of the dual solution $\vp$ itself. This gradient satisfies a
conservation law instead of a transport equation, and it can therefore be
computed with the same algorithm as the forward problem. For discontinuous
transport coefficients, the new conservative algorithm for $w$ is more robust
than transport schemes for $\vp$, see \cite{SteinerNoelle}. Here we also derive
characteristic boundary conditions for the conservative dual problem, which we
use in the numerical examples in Sections~\ref{section.numexp_fully_implicit}
and \ref{section.numexp_explicit_implicit}.

In order to compute the adjoint error representation one needs to compute a
forward and a dual problem and to assemble the space-time scalar product
\eqref{eq.localized_error_indicator}. Together, this costs about three times as
much as the computation of a single forward problem. In our application, the
error representation is computed on a coarse mesh ($L=2$), and therefore it
presents only a minor computational overhead compared with the fine grid
solution ($L=5$). In other applications, the amount of additional storage and
CPU time may become significant. In such cases, checkpointing strategies might
help (see, e.g., \cite{MR2204516}).

We have implemented and tested both a fully implicit and a mixed
explicit-implicit timestepping strategy. The explicit-implicit approach switches
to an explicit timestep with $\cfl=0.5$ in case the adaptive strategy suggests
an implicit timestep with $\cfl<5$. Clearly, the mixed explicit-implicit
strategy is the most accurate and efficient, beating the adaptive fully implicit
in accuracy and efficiency, the implicit approach with fixed $\cfl$ numbers in
accuracy, and the fully explicit approach in efficiency.

Finally, we would like to stress that the computational cost in each timestep is
nearly constant, no matter if the $\cfl$ number is of order 1 or order 100. In
all cases, the solver needs only 1 or 2 Newton iterations per timestep to reach
the (rather strict) break condition, which is related to the multi\replace{scale}{resolution}
\replace{analysis}{decomposition}. This seems to be another major benefit of the adaptive timestep
control.

\bibliographystyle{siam}
\bibliography{christina_lit2}
%\bibliography{christinatlit}

% \include{answer}
\end{document}